\documentclass[reqno]{amsart}
\usepackage{textcomp}

\usepackage[dvipsnames]{xcolor}
\usepackage{amsmath}
\usepackage{amsthm}
\usepackage{hyperref}
\DeclareRobustCommand{\quarticfamilylink}[2]{{\hypersetup{linkcolor=black}\hyperlink{quartic-family-#1}{#2}}}
\newcommand{\quarticfamilytarget}[1]{\raisebox{1.5ex}[0pt][0pt]{\hypertarget{quartic-family-#1}{}}}
\usepackage{amsfonts,graphics,amsthm,amsfonts,amscd,latexsym}
\usepackage{tikz-cd}
\usepackage{epsfig}
\usepackage{flafter}
\usepackage{longtable}
\usepackage{mathtools}
\usepackage{comment}
\usepackage{stmaryrd}

\usepackage{mathabx,epsfig}

\hypersetup{
    colorlinks=true,    
    allcolors=blue,
    linkcolor=blue,          
    citecolor=blue,      
    filecolor=blue,      
    urlcolor=blue           
}
\usepackage{tikz}
\usetikzlibrary{graphs,positioning,arrows,shapes.misc,decorations.pathmorphing}

\tikzset{
    >=stealth,
    every picture/.style={thick},
    graphs/every graph/.style={empty nodes},
}

\tikzstyle{vertex}=[
    draw,
    circle,
    fill=black,
    inner sep=1pt,
    minimum width=5pt,
]
\usepackage[position=top]{subfig}
\usepackage{amssymb}
\usepackage{color}

\calclayout

\usetikzlibrary{decorations.pathmorphing}
\tikzstyle{printersafe}=[decoration={snake,amplitude=0pt}]

\newcommand{\ord}{\operatorname{ord}}

\newcommand{\pp}{\mathbb{P}}

\newcommand{\qq}{\mathbb{Q}}

\newcommand{\rr}{\mathbb{R}}
\newcommand{\cc}{\mathbb{C}}

\newcommand{\G}{\mathbb{G}}

\def\O#1.{\mathcal {O}_{#1}}			
\def\pr #1.{\mathbb P^{#1}}				
\def\af #1.{\mathbb A^{#1}}			
\def\ses#1.#2.#3.{0\to #1\to #2\to #3 \to 0}	
\def\xrar#1.{\xrightarrow{#1}}			
\def\K#1.{K_{#1}}						
\def\bA#1.{\mathbb{A}_{#1}}			
\def\bM#1.{\mathbb{M}_{#1}}				
\def\bL#1.{\mathbb{L}_{#1}}				
\def\bB#1.{\mathbb{B}_{#1}}				
\def\bK#1.{\mathbb{K}_{#1}}			
\def\subs#1.{_{#1}}					
\def\sups#1.{^{#1}}

\usepackage{tikz}
\usetikzlibrary{matrix,arrows,decorations.pathmorphing}

  \newtheorem{theorem}{Theorem}[section]
  \newtheorem{lemma}[theorem]{Lemma}
  \newtheorem{proposition}[theorem]{Proposition}
  \newtheorem{corollary}[theorem]{Corollary}

  \newtheorem{notation}[theorem]{Notation}
  \newtheorem{definition}[theorem]{Definition}
  \newtheorem{example}[theorem]{Example}

\newtheorem{remark}[theorem]{Remark}

\theoremstyle{remark}

\numberwithin{equation}{section}

\usepackage[all]{xy}

\begin{document}

\thanks{J.M. was partially supported by NSF research grant DMS-2443425.}

\title[Cluster type boundaries on $\mathbb{P}^3$]{Cluster type boundaries on $\mathbb{P}^3$} 

\author[J.~Moraga]{Joaqu\'in Moraga}
\address{UCLA Mathematics Department, Box 951555, Los Angeles, CA 90095-1555, USA
}
\email{jmoraga@math.ucla.edu}

\author[T.~Oda]{Tomoki Oda}
\address{UCLA Mathematics Department, Box 951555, Los Angeles, CA 90095-1555, USA
}
\email{tomokioda0723@math.ucla.edu}


\begin{abstract}
We classify cluster type boundaries on the three-dimensional complex projective space. In other words, we classify semi-log canonical quartic surfaces $S\subset \pp^3$ for which there is an open torus embedding $\mathbb{G}_m^3 \hookrightarrow \pp^3\setminus S$. As a consequence, there are precisely $92$ cusp singularities $T_{p,q,r}$ appearing in quartic surfaces $S$ for which the pair $(\pp^3,S)$ is of cluster type. Finally, we provide new examples of toric degenerations of $\pp^3$ for which the whole degeneration space is of cluster type.
\end{abstract} 

\maketitle

\setcounter{tocdepth}{1} 
\tableofcontents

\section{Introduction}
Toric varieties are one of the most ubiquitous classes of algebraic varieties, appearing in essentially every field within algebraic geometry and commutative algebra~\cite[Appendix A]{CLS11}. Recently, there has been substantial research on a new class of algebraic varieties called {\em cluster type} varieties~(see, e.g.,~\cite{Cor25,EFMS25,ELY25,JM24,MZ26}). This class naturally generalizes the notion of toric varieties. A cluster type variety, roughly speaking, is a compactification $X$ of an algebraic torus $\mathbb{G}_m^n$ 
such that the volume form of the algebraic torus has no zeros on $X$ (see~\cite[Definition~2.6]{Mor26}). 
A {\em cluster type boundary} on $X$ is defined by the poles of the volume form of the torus.
Cluster type varieties have importance in mirror symmetry, Fano geometry, and birational geometry~\cite{Cor25,Lutz26,KY23}. The first author together with Alves da Silva, Figueroa and Y\'a\~nez proved that most rational smooth Fano threefolds are cluster type~\cite[Theorem 1.8]{AdSFMY26}. Thus, it is natural to study cluster type Fano varieties and their possible cluster type boundaries. The cluster type boundaries on $\pp^2$ are easy to describe: these are cubic curves whose only singularities are nodes, with at least one node~\cite[Proposition~1.3]{GHK15}. Cluster type boundaries on $\pp^3$ are much harder to understand and some partial results can be found in~\cite{EFM26}. The aim of this paper is to fully classify the cluster type boundaries on $\pp^3$ and use them to study degenerations of $\pp^3$. To understand which pairs $(\pp^3,S)$ are cluster type, one needs to thoroughly understand the birational geometry and singularities of the quartic surface $S$. 
In order to achieve our classification, we will build upon the vast literature on rational quartic surfaces developed by
Noether, Jessop, Umezu, Urabe, Ishii, Nakayama, Mella, Ducat, Alves da Silva, Figueroa, among others (see, e.g.,~\cite{Noe89,Jes16,Ume84, Ura86,IN04,Mel20, Duc24,EFM26,Shah81}).

While cluster type varieties have some rigidity due to the existence of an algebraic torus, they show some interesting behaviour from the perspective of moduli theory. 
For instance, Ji and the first author proved that the cluster type condition is constructible in families of smooth Fano varieties~\cite[Theorem~1.2]{JM24}. On the other hand, Z\'u\~niga and the first author have studied degenerations of cluster type varieties; under some mild conditions the central fiber is a finite quotient of a cluster type variety~\cite[Theorem~1.2]{MZ26}. 
Our first theorem generalizes these two statements and gives a locally constant property for cluster type pairs whenever we control some global invariants of the family and the deeper singularities of the pair.

\begin{theorem}\label{thm:cluster-type-in-families}
Let $\phi\colon (\mathcal{X},\mathcal{B})\rightarrow T$ be a family of pairs satisfying the following conditions:
\begin{enumerate}
\item \label{cond:TF} for $t\in T$ closed the fiber $\mathcal{X}_t$ is a terminal Fano threefold; 
\item \label{cond:CY} for $t\in T$ closed the pair $(\mathcal{X}_t,\mathcal{B}_t)$ is log Calabi--Yau; 
\item \label{cond:CL} the class group restriction ${\rm Cl}(\mathcal{X}/T)\rightarrow {\rm Cl}(X_t)$ is an isomorphism for each $t\in T$ closed; 
\item  \label{cond:simult-q-dlt} the family $\phi$ admits a simultaneous $\qq$-factorial dlt modification; and 
\item \label{cond:normal-connected-base} the variety $T$ is smooth and quasi-projective.
\end{enumerate} 
Then, the cluster type property is constant on the log fibers of $\phi$.
\end{theorem} 

Theorem~\ref{thm:cluster-type-in-families} shows that being of cluster type is
constant in equisingular families of pairs over a smooth base.
Applying Theorem~\ref{thm:cluster-type-in-families} to the trivial family
$\mathbb P^3\times T\to T$, we obtain the following consequence for families of slc quartics.

\begin{corollary}\label{cor:p3-families}
Let $B_t\subset \mathbb{P}^3$ be a family of slc 
quartics parametrized by the closed points $t$ of a  smooth connected variety $T$. Assume that \'etale locally on $T$ the pairs $(\pp^3,B_t)$ admit simultaneous $\qq$-factorial dlt modifications.
Then, the cluster type property of $(\mathbb{P}^3,B_t)$ is independent of $t$.
\end{corollary} 

This makes it possible to study the cluster type boundaries on $\mathbb{P}^3$ one equisingular stratum at a time.
Combining this with the classification of semi-log canonical quartic surfaces,
we aim to complete the classification of the boundaries on $\mathbb{P}^3$ that are cluster type.

The classification of cluster type boundaries on $\mathbb P^3$ is already known when $S$ is reducible. When the boundary $S$ consists of four irreducible components, then the pair $(\mathbb P^3,S)$ is a toric pair, hence a cluster type pair. When the boundary $S$ consists of three irreducible components, then $(\mathbb P^3,S)$ is of cluster type~(\cite[Theorem~6]{ELY25}). When the boundary $S$ consists of two irreducible components, Alves da Silva, Figueroa, and the first author determined exactly when the pair fails to be of cluster type~(\cite[Theorem~5.1 and Proposition~5.2]{EFM26}). We recall one of their main results below.

\begin{theorem}\label{introthm3}
Let $(\mathbb{P}^3, S)$ be a log Calabi--Yau pair of index one and coregularity zero.\footnote{Roughly speaking, coregularity zero means that in some suitable log resolution the boundary has maximal intersection.} Then the following hold.
\begin{enumerate}
\item Assume that $S$ decomposes into two irreducible components. Then, the pair $(\mathbb{P}^3,S)$ is not of cluster type if and only if $S= H+C$, where $H$ is a plane and $C$ is a cubic surface with no singularities along $H$.
\item Assume there is a point $p\in \mathbb{P}^3$ for which $\operatorname{mult}_p(S)\geq 3$. Then, the pair $(\mathbb{P}^3,S)$ is of cluster type.
\end{enumerate}
\end{theorem}

Thus, it remains to study the cluster type property of the pair $(\mathbb P^3,S)$ when $S$ is irreducible and has only multiplicity two singularities.
The classification of strictly log canonical quartics is a classical subject in algebraic geometry. The birational study of such surfaces goes
back to Noether, who analyzed rational quartics in terms of their
singularities~\cite{Noe89}. The classical theory was later summarized by Jessop~\cite{Jes16}. More recently, Mella proved
that every rational quartic surface in $\mathbb P^3$ is Cremona equivalent
to a plane~\cite{Mel20}. Singular quartics
also appear in the study of compactifications
of the moduli space of $K3$ surfaces defined by quartics.
Shah studied GIT description of generations of quartic surfaces (see, e.g.,~\cite{Shah81}).
Urabe and Degtyarev also
studied non-normal and normal quartic surfaces
with cusp singularities, respectively~\cite{Ura86,Deg90}. 
Ducat proved that every quartic pair
$(\mathbb P^3,S)$ of coregularity zero 
is log rational~\cite[Theorem~1.2]{Duc24}.
For a pair, being of cluster type is a stronger condition than being
log rational.\footnote{
We say that $(X,B)$ is \emph{log rational} if it is
crepant birational to the toric pair
$(\mathbb P^n,H_0+\cdots+H_n)$, where the $H_i$
are the coordinate hyperplanes.
}
We aim to determine which of these quartics $S$
give a cluster type pair $(\mathbb{P}^3,S)$.
In order to do so, first we will refine the classical classification of slc quartic surfaces. 
We focus on what kind of log canonical centers happen in the non-normal locus.
In particular, we obtain explicit bounds on the indices of the cusp singularities.  
These bounds depend on the geometry of the surface and its non-normal locus.

\begin{theorem}\label{introthm:non-normal-class}
Let $S\subset\mathbb P^3$ be an irreducible semi-log canonical quartic surface such that the pair $(\mathbb{P}^3,S)$ has coregularity zero.
Assume that every singularity of $S$ has multiplicity two.
Then, the following statements hold:
\begin{enumerate}
\item If $S$ is normal, then it has a unique cusp singularity of type
$T_{2,q,r}$ with $q+r\leq 18$. Moreover, the equality implies that $S$ is smooth 
outside the cusp.
\item If $S$ is non-normal, then its reduced non-normal locus $\Sigma$ is one of the following: a line, a smooth plane conic, a twisted cubic, two intersecting lines, 
a line and a smooth conic meeting transversally at a point, 
three concurrent non-coplanar lines, or a chain of three lines with the first and third skew lines. 
\item If $\Sigma$ is a line or a smooth conic, then $S$ has at
most two degenerate cusps of type $T_{2,q,\infty}$ along $\Sigma$.
Moreover, for each cusp we have 
\[
q\leq 
\begin{cases}
10 & \text{if $\Sigma$ is a line},\\
7  & \text{if $\Sigma$ is a smooth conic}.
\end{cases}
\]
\end{enumerate}
\end{theorem}

We first study cluster type pairs $(\mathbb P^3,S)$
with $S$ a normal quartic surface.
By Wall's classification, these surfaces fall into
four families, denoted by \textup{(A.1)}--\textup{\quarticfamilylink{A4}{(A.4)}}
following Ducat's notation
\cite[Proposition~3.1]{Duc24}.
Surfaces in the \textup{(A.1)} family have a
triple point and hence, by
Theorem~\ref{introthm3}, are of cluster type
when $(\mathbb P^3,S)$ has coregularity zero.
Thus, it remains to consider families
\textup{\quarticfamilylink{A2}{(A.2)}}--\textup{\quarticfamilylink{A4}{(A.4)}}.
We begin with the \textup{\quarticfamilylink{A2}{(A.2)}} family, which contains
quartics with a simple elliptic singularity
$\widetilde E_7$ and its cusp degenerations $T_{2,q,r}$,
where $4\le q\le r$ and $(q,r)\ne(4,4)$.
The natural strategy is to consider the projection from the cusp and study the associated birational double cover $\pi_0\colon S\dashrightarrow \pp^2$. However, for this double cover is hard to describe the singularities of the branch locus in terms of the singularities of $S$. To fix this issue, we perform a birational modification of the ambient space first. We consider a birational map
$\psi\colon \pp^3\dashrightarrow \mathbb P(1,1,1,2)$,
which is toric in the ambient space. Let $T$ be the birational transform
of $S$ and $F$ the divisor contracted by $\psi^{-1}$.
Then $\psi\colon(\pp^3,S)\dashrightarrow
(\mathbb P(1,1,1,2),T+F)$ is a crepant birational map.
We let $E=F|_T$ be the exceptional divisor of
$\psi_T^{-1}\colon T\dashrightarrow S$. In the ambient space $\mathbb{P}(1,1,1,2)$ we can consider the projection out of a torus invariant point (see Subsection~\ref{sec:A3}). Analyzing the induced double cover $T\rightarrow \pp^2$ and the quartic cover of $\pp(1,1,1,2)$ branching along $T$, we deduce the following.

\begin{theorem}\label{thm:intro-A2}
Let $S\subset \pp^3$ be a quartic surface in the \textup{\quarticfamilylink{A2}{(A.2)}} family with a cusp of type $T_{2,q,r}$. 
Then $(\pp^3,S)$ is of cluster type if and only if $T$ has an $A_m$ singularity, with $m\geq 2$, along $E$.
In particular, if $r\leq 6$ then $(\pp^3,S)$ is not cluster type and if $q+r\geq 13$ and $(q,r)\neq (6,7)$ then $(\pp^3,S)$ is of cluster type.
\end{theorem}

This theorem shows that, in the family \textup{\quarticfamilylink{A2}{(A.2)}}, the cusp type alone does not determine the cluster type property when
$(q,r)\in\{(4,7),(4,8),(5,7),(6,7)\}$.
In these cases, one must also examine how the cusp singularity arises.
For families \textup{\quarticfamilylink{A3}{(A.3)}} and \textup{\quarticfamilylink{A4}{(A.4)}}, which we consider next, however, the choice of family together with the singularity type suffices to determine the cluster type property. Although \textup{\quarticfamilylink{A3}{(A.3)}} and \textup{\quarticfamilylink{A4}{(A.4)}} are distinct families, they admit the same singularity types: a simple elliptic singularity $\widetilde E_8$ or a cusp of type $T_{2,3,r}$ with $r\geq7$.
Their global geometry differs: surfaces in \textup{\quarticfamilylink{A3}{(A.3)}} contain a line through the cusp, whereas those in \textup{\quarticfamilylink{A4}{(A.4)}} do not~\cite{Noe89}.
This distinction is also reflected in their GIT behavior, as studied by Shah~\cite{Shah81}.
It affects the cluster type property as well: for the same cusp type $T_{2,3,11}$, the pair $(\mathbb P^3,S)$ is of cluster type in the family \textup{\quarticfamilylink{A3}{(A.3)}}, but not in the family \textup{\quarticfamilylink{A4}{(A.4)}}.
More generally, we obtain the following classification. 

\begin{theorem}\label{them:mainclusternormal}
Let $S\subset\mathbb P^3$ be a normal quartic
surface in the family \textup{\quarticfamilylink{A3}{(A.3)}} or
\textup{\quarticfamilylink{A4}{(A.4)}} with a cusp singularity
of type $T_{2,3,r}$, where $r\geq 7$.
Then the following hold:
\begin{enumerate}
\item If $S$ belongs to the \textup{\quarticfamilylink{A3}{(A.3)}}
family, then $(\mathbb P^3,S)$ is of cluster
type if and only if $r\geq 11$.
\item If $S$ belongs to the \textup{\quarticfamilylink{A4}{(A.4)}}
family, then $(\mathbb P^3,S)$ is of cluster
type if and only if $r\geq 12$.
\end{enumerate}
\end{theorem}

We now turn to the cluster type property of
$(\mathbb P^3,S)$ when $S$ is non-normal.
By Ducat's classification
\cite[Proposition~3.3]{Duc24}, non-normal quartic
surfaces fall, up to degeneration, into four
families \textup{\quarticfamilylink{B1}{(B.1)}}--\textup{(B.4)},
which we refer to as the \emph{Urabe families},
following Urabe's earlier classification~\cite{Ura86}.
The general members of \textup{(B.3)} and
\textup{(B.4)} were already shown to be of
cluster type in \cite[Theorem~6.1]{EFM26}.

For non-normal quartics, Theorem~\ref{introthm:non-normal-class}
describes the possible configurations of the non-normal
locus. We refine Urabe's classification by distinguishing
the corresponding families more precisely.
When the non-normal locus is irreducible, a line gives
rise to two families \textup{\quarticfamilylink{B1-alpha}{(B.1;$\alpha$)}} and
\textup{\quarticfamilylink{B1-beta}{(B.1;$\beta$)}}, distinguished by the decomposition
of the conductor on the normalization.
Although these were treated together in Urabe's
classification, they exhibit different cluster type
properties.
\quarticfamilytarget{B3-alpha}The smooth conic and twisted cubic cases correspond
to \textup{\quarticfamilylink{B2-alpha}{(B.2;$\alpha$)}} and
\textup{\quarticfamilylink{B3-alpha}{(B.3;$\alpha$)}}, respectively.

When the non-normal locus is reducible, we distinguish
four further families:
\textup{\quarticfamilylink{B2-beta}{(B.2;$\beta$)}}, consisting of two intersecting lines;
\textup{\quarticfamilylink{B3-beta}{(B.3;$\beta$)}}, consisting of a smooth conic and
a line meeting it transversely at one point;
\quarticfamilytarget{B4-alpha}\textup{\quarticfamilylink{B4-alpha}{(B.4;$\alpha$)}}, consisting of three concurrent
non-coplanar lines; and
\textup{\quarticfamilylink{B4-beta}{(B.4;$\beta$)}}, consisting of a chain of three lines
whose first and last components are skew.
Using this refinement, we obtain the following
classification of cluster type pairs with non-normal
quartic boundaries.

\begin{theorem}\label{them:mainclusternonnormal}
 Let $S\subset\mathbb P^3$ be an irreducible non-normal slc quartic
such that $(\mathbb P^3,S)$ has coregularity zero. Assume that
every singularity of $S$ has multiplicity at most two.  Let $q$ denote the largest index of a degenerate cusp
of type $T_{2,q,\infty}$ on $S$, allowing $q=\infty$. Then $(\mathbb P^3,S)$ is of cluster type if and only if one of the following holds:
\begin{enumerate}
\item $S$ belongs to one of the families
\textup{\quarticfamilylink{B1-alpha}{(B.1;$\alpha$)}},
\textup{\quarticfamilylink{B3-beta}{(B.3;$\beta$)}} or
\textup{\quarticfamilylink{B4-beta}{(B.4;$\beta$)}};

\item $S$ belongs to \textup{\quarticfamilylink{B1-beta}{(B.1;$\beta$)}}
and $q\geq 7$;

\item $S$ belongs to \textup{\quarticfamilylink{B2-alpha}{(B.2;$\alpha$)}}
and $q\geq 4$;

\item $S$ belongs to \textup{\quarticfamilylink{B2-beta}{(B.2;$\beta$)}} and $(\mathbb P^3,S)$ has two zero dimensional log canonical centers;

\item $S$ belongs to \textup{\quarticfamilylink{B3-alpha}{(B.3;$\alpha$)}}
and its non-normal locus is a twisted cubic.
\end{enumerate}
\end{theorem}

Putting Theorem~\ref{introthm3}, Theorem~\ref{thm:intro-A2}, Theorem~\ref{them:mainclusternormal} and Theorem~\ref{them:mainclusternonnormal} together, we get a complete classification of the cluster type boundaries on $\mathbb{P}^3$.
The following figure is a closed ball representation of the real spectrum of a normal quartic surface \(S\) which gives rise to a cluster type pair $(\pp^3,S)$. 

\begin{figure}[ht]
\centering
\includegraphics[scale=0.4]{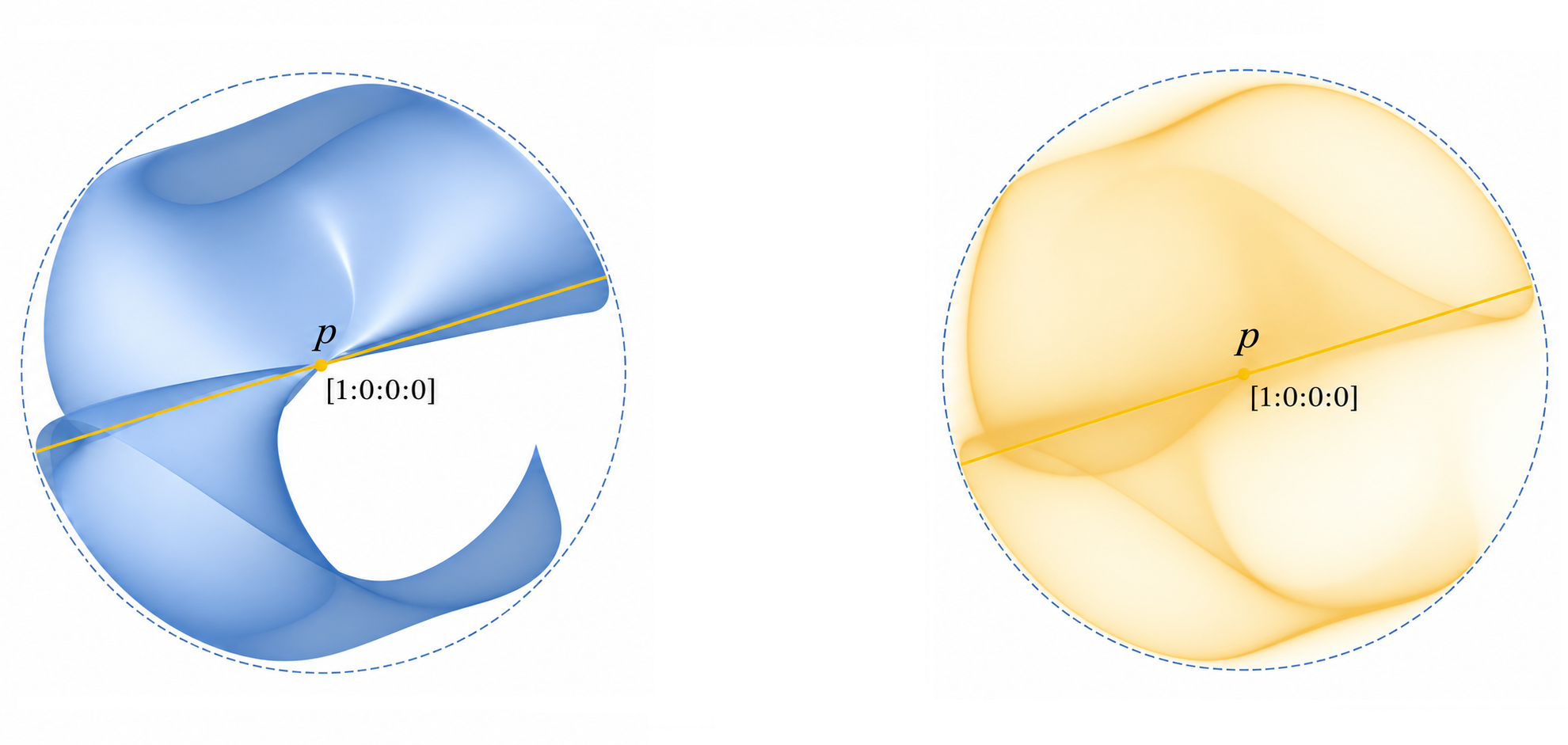}
\caption{Normal quartic surface $V(F)$ and representation of the volume form.}
\end{figure}

The figure on the left is the real locus of the vanishing locus of homogeneous quartic 
\[
F
=
w^2z^2+wy^3+wxyz+2x^3z
+\frac54x^2y^2+yz^3 \in \cc[x,y,z,w]
\]
near $p=[1:0:0:0]\in \pp^3_\rr$. The figure on the right is the level-set representation of the volume form of an algebraic torus $\mathbb{G}_m^3\hookrightarrow \pp^3_\rr \setminus V(F)$. This quartic is studied in detail in Proposition~\ref{prop:A3-low-strata}.(iii). 

The last section of the paper, Section~\ref{sec:cusp-sing-cluster-type-deg}, has two aims.
The first one is to count the number of isolated log canonical centers appearing in normal cluster type boundaries on $\mathbb{P}^3$.
The second one is to describe some new klt degenerations of $\pp^3$ and show how these normal cluster type boundaries on $\pp^3$ induce cluster type boundaries on the klt degenerations.
 
The isolated strictly log canonical singularities appearing on surfaces are called cusp singularities (see Proposition~\ref{prop:slc_local_models}).
Cusp singularities, which below we denote by $T_{p,q,r}$-singualrities, are determined by three positive integers $p,q,r\geq 2$.
There are exactly $101$ analytic type of cusp singularities $T_{p,q,r}$ of coregularity zero.
In Proposition~\ref{prop:101+92}, we will prove that precisely $92$ of these cusp singularities appear in quartic surfaces $S$ for which the pair $(\pp^3,S)$ is of cluster type. Thus, one may say that most cusp singularities of coregularity zero appear in cluster type boundaries of $\pp^3$. 

Finally, the topic of klt degenerations of $\pp^3$ is quite vast. 
The klt degenerations of $\pp^2$ are fully classified by the work of Hacking and Prokhorov in~\cite{HP10}. On the contrast, we are not close to understand a full classification of klt degenerations of $\pp^3$. 
H\"oring and Peternell classified all the degenerations of $\pp^3$ with canonical singularities (see~\cite[Theorem 1.3]{HP24}). On the other hand, DeVleming gave some interesting non-canonical klt degenerations of $\pp^3$ in her Ph.D. thesis (see~\cite{Dev22}). 
There is however a recent heightened interest on klt degenerations of $\pp^3$ including some very recent preprints due to DeVleming, Li, and Torres and Chen and Lee (see, e.g.,~\cite{ADVL23,DVLT26,CL26}). In Section~\ref{sec:cusp-sing-cluster-type-deg}, we will provide some new examples of toric degenerations of $\pp^3$.
We will prove that these toric degenerations admit a cluster type boundary obtained by taking limits of some normal cluster type boundary on $\pp^3$. We hope that this is a good starting point to better understand cluster type klt degenerations of the three-dimensional projective space.

In Appendix~\ref{appendix-A}, we explain how AI was used in this paper. 

\subsection{Strategy of the proof}\label{strat}
Proving that a pair is of cluster type
is somewhat similar to proving rationality of an algebraic variety. Indeed, one needs to explicitly provide the construction of a torus embedding with the property of the volume form. As in irrationality problems, proving that a pair is not of cluster type tends to be more challenging. In this case, one needs to show the existence of an obstruction theory for some suitable modification of the variety. 

If Theorem~\ref{thm:intro-A2}, Theorem~\ref{them:mainclusternormal}  or Theorem~\ref{them:mainclusternonnormal} states that $(\pp^3,S)$ is of cluster type, we explicitly construct a torus embedding.
Our main tool is Lemma~\ref{prop:shear-torus-chart}, which gives a general method for embedding a torus into an affine threefold. Example~\ref{ex:shear-instances} then describes situations in which this lemma applies to the complement of these quartic surfaces.

For the other pairs, we aim to prove that $(\pp^3,S)$ is not of cluster type.
First, for each such family we will prove that the corresponding
marked equisingular strata in the Wall and Urabe coefficient spaces
are smooth algebraic varieties.
We also establish the existence of simultaneous
$\mathbb Q$-factorial dlt modifications over these strata. These statements are proved in
Section~\ref{sec:normal-quartic} and Section~\ref{sec:nonormal}.
Corollary~\ref{cor:p3-families} therefore applies.
Consequently, the cluster type property is constant on each
connected component of these strata.
Hence, to prove that all members of a connected
component are not of cluster type, it suffices to find a single
member of that component which is not of cluster type.
Let $S_0$ be such a quartic surface.
We consider the associated cyclic quartic cover
\[
\pi_0\colon Y_0 \rightarrow \pp^3_{w,x,y,z}
\text{ defined by } 
Y_0:=\{ [w:x:y:z:u] \mid u^4 = F_0(w,x,y,z)\}
\subset \pp^4
\]
where $F_0$ is the defining homogeneous equation of $S_0$. This finite cover is fully ramified along $S_0$. In Lemma~\ref{lem:clustertypeobstruction}.(2), we will prove that
$Y_0$ is rational whenever $(\pp^3,S_0)$ is of cluster type.
Thus, it suffices to prove that the quartic threefold $Y_0$ is irrational.

Furthermore, it is not necessary to study all equisingular strata separately. By Lemma~\ref{lem:clustertypeobstruction}.(1), the variety $Y_0$ is rationally connected. 
Using a theorem by de Fernex and Fusi, if $Y_0$ is irrational, then the quartic cover associated to a very general member of $\mathcal{T}$ is irrational (see~\cite{dFF13} and Proposition~\ref{prop:irrational-covers-in-families}). Thus, it suffices to exhibit a single member with an irrational cyclic quartic cover for some special strata. These strata are the following:
\begin{gather*}
 \mathcal T_{2,4,8}^{\quarticfamilylink{A2}{(A.2)},\mathrm I},\quad
 \mathcal T_{2,4,8}^{\quarticfamilylink{A2}{(A.2)},\mathrm{II}},\quad
 \mathcal T_{2,6,6}^{\quarticfamilylink{A2}{(A.2)},\mathrm I},\quad
 \mathcal T_{2,6,6}^{\quarticfamilylink{A2}{(A.2)},\mathrm{IV}},\quad
 \mathcal T_{2,6,7}^{\quarticfamilylink{A2}{(A.2)},\mathrm I},\quad
 \mathcal T_{2,6,7}^{\quarticfamilylink{A2}{(A.2)},\mathrm{II}},\\
 \mathcal T_{2,3,10}^{\quarticfamilylink{A3}{(A.3)}},\quad
 \mathcal T_{2,3,10}^{\quarticfamilylink{A4}{(A.4)},\mathrm{II}},\quad
 \mathcal T_{2,3,11}^{\quarticfamilylink{A4}{(A.4)}},\quad
 \mathcal T_{2\times(2,4,\infty)}^{\quarticfamilylink{B1-beta}{(B.1;\beta)}},\quad
 \mathcal T_{(2,3,\infty),(2,6,\infty)}^{\quarticfamilylink{B1-beta}{(B.1;\beta)}},\quad
 \mathcal T_{2,6,\infty}^{\quarticfamilylink{B1-beta}{(B.1;\beta)},\mathrm{II}},\quad
 \mathcal T_{2\times(2,3,\infty)}^{\quarticfamilylink{B2-alpha}{(B.2;\alpha)}}
 \quad\text{and}\quad
 \mathcal T_{2,\infty,\infty}^{\quarticfamilylink{B2-beta}{(B.2;\beta)}}.
\end{gather*}
Here, the subscripts record the singularity types that appear in the surfaces, while the superscripts record the Wall and Urabe family. The Roman number denotes the connected components of such families.

For the strata
$\mathcal T_{2,4,8}^{\quarticfamilylink{A2}{(\mathrm{A.2})},\mathrm{I}}$,
$\mathcal T_{2,4,8}^{\quarticfamilylink{A2}{(\mathrm{A.2})},\mathrm{II}}$, and
$\mathcal T_{2,6,6}^{\quarticfamilylink{A2}{(\mathrm{A.2})},\mathrm{I}}$, the cover $Y_0$
is birational to a smooth quartic double solid and hence is
irrational. For the remaining strata, we construct a standard
conic bundle model of $Y_0$.
Some of these conic bundles can be obtained using projections from
singular points of $Y_0$ in $\pp^4$.
We follow the approach of Clemens, Griffiths, Mumford, Beauville, and Shokurov to study the irrationality of
these conic bundles.
For a standard conic bundle over a rational surface with a smooth
irreducible discriminant curve $\Gamma_0$, Shokurov's criterion
provides sufficient conditions for irrationality
(see~\cite{Beau77,Sho83} and Proposition~\ref{Thm:Shokurov}).

In the \textup{\quarticfamilylink{B2-alpha}{(B.2;$\alpha$)}} case, the discriminant
curve $\Gamma_0$ is a cyclic quartic cover of an elliptic curve.
For the other conic bundles, $\Gamma_0$ is a cyclic quartic cover
of $\mathbb P^1$, and its ramification can be read directly from
the defining equation.
Let $B_4$ and $B_2$ denote the numbers of branch points of
$\Gamma_0\to\mathbb P^1$ with ramification indices $4$ and $2$,
respectively. By Proposition~\ref{prop:fractional-kummer-irrationality},
Shokurov's criterion applies if the double cover of $\Gamma_0$
induced by the Brauer residue is connected and \'etale and the
inequalities $B_4\geq 4$ and $3B_4+2B_2\geq 20$ hold.
Thus, it suffices to prove these inequalities in the remaining cases.

\section{Preliminaries}
Throughout this article, we will work over the field of complex numbers $\mathbb{C}$. We may write $(X;x)$ for the germ of an algebraic variety $X$ near a closed point $x\in X$. In this section, we collect some preliminary results regarding dlt modifications (Subsection~\ref{subsec:sim-dlt-mod}), cluster type pairs (Subsection~\ref{subsec:cluster-type}) and some obstructions for a pair to be cluster type (Subsection~\ref{subssec:obstruction-cluster-type}). 
Throughout the paper, we will write $\mathbb{G}_m^n$ for the $n$-dimensional algebraic torus over $\mathbb{C}$. 
We may also write $\mathbb{T}$ or $\mathbb{T}_{x_1,\dots,x_n}$ to represent an algebraic torus or an algebraic torus with coordinates $x_1,\dots,x_n$, i.e., with a fixed isomorphism $\mathbb{T}_{x_1,\dots,x_n}:=\operatorname{Spec}(\mathbb C[x_1^\pm,\dots,x_n^\pm])$. Given an algebraic variety $X$ and a regular function $f$ on $X$, we write $D_X(f)$ for the distinguished open set. We use the same notation for the nonvanishing locus
of a homogeneous coordinate on a projective or
weighted projective variety.

\begin{definition}
{\rm 
Let $f\colon X\rightarrow T$ be a morphism
and $B$ be a boundary divisor on $X$. 
Let $t\in T$ be a point and let $X_t$ be the corresponding fiber. 
Assume that $B$ does not contain the fiber $X_t$. Denote by $B_t$ the pull-back of $B$ to the fiber $X_t$. We call the pair $(X_t,B_t)$ the {\em log fiber} of $f\colon (X,B)\rightarrow T$ over $t\in T$.
If $B$ contains no closed fiber of $f$, we may say that $(X,B)\rightarrow T$ is a {\em family of pairs}.
}
\end{definition}

\subsection{Simultaneous dlt modifications}
\label{subsec:sim-dlt-mod}
In this subsection, we recall the concept of dlt modifications. For the standard definitions of singularities of the MMP we refer the reader to~\cite{KM98} and~\cite{K13}. We introduce the concept of simultaneous $\mathbb Q$-factorial dlt modification. 

\begin{definition}
\em Let $(X,\Delta)$ be a log canonical pair. We say that $(X,\Delta)$ is \emph{divisorial log terminal} (or \emph{dlt}) if there exists an open set $U\subset X$ such that 
\begin{enumerate}
\item $(U,\Delta|_{U})$ is a simple normal crossing pair and
\item every log canonical center of $(X,\Delta)$ intersects $U$ non-trivially and is a stratum of $\lfloor \Delta\rfloor$.
\end{enumerate}

\end{definition}

\begin{definition}
\em Let $(X,\Delta)$ be a log canonical pair. Then a \emph{dlt modification} is a projective crepant birational morphism $\phi\colon (Y,\Delta_Y) \rightarrow (X,\Delta)$ such that $(Y,\Delta_Y)$ is dlt
and every exceptional divisor of $\phi$ is contained in $\lfloor \Delta_Y\rfloor$.
If furthermore $Y$ is $\qq$-factorial, we call it a {\em $\qq$-factorial dlt modification}.
\end{definition}

\begin{definition}
{\em 
Let $\phi\colon (\mathcal{X},\mathcal{B})\rightarrow T$ be a family of log pairs over an irreducible variety $T$. We say that $\phi$ admits a {\em simultaneous dlt modification} if there exists a projective crepant birational morphism $\psi\colon (\mathcal{Y},\mathcal{B}_{\mathcal{Y}})\rightarrow (\mathcal{X},\mathcal{B})$ satisfying the following properties: 
\begin{enumerate} 
\item for each closed $t\in T$, we have that 
$(\mathcal{Y}_t,\mathcal{B}_{\mathcal{Y}_t})\rightarrow (\mathcal{X}_t,\mathcal{B}_t)$ is a dlt modification; and 
\item every log canonical center of $(\mathcal{Y}_t,\mathcal{B}_{\mathcal{Y}_t})$ with codimension $c$ is the restriction to $\mathcal{Y}_t$
of a log canonical center of 
$(\mathcal{Y},\mathcal{B}_{\mathcal{Y}})$ with codimension $c$.
\end{enumerate} 
We further say that it is a {\em simultaneous $\qq$-factorial dlt modification} if each $\mathcal{Y}_t$ is a $\qq$-factorial variety.
}
\end{definition} 

\begin{definition}
{\em 
Let $m\geq 0$ be an integer.
We define by $B_m$ the local isomorphism class of the pair $(\mathbb{A}^3_{x,y,z},D_{1,m}+D_2)$ near the origin, where $D_1:=\{xy+z^m=0\}$ and
$D_2=\{z=0\}$. 
Note that $B_m$ is only when dlt when $m=0$.
}
\end{definition} 

The following local computation will be used repeatedly to produce explicit
$\mathbb{Q}$-factorial dlt modifications.

\begin{lemma}[cf.~{\cite[Theorem~4.57]{KM98}}]
\label{lem:uv-sm-dlt}
For $m\geq 0$, the germ $B_m$ admits an explicit $\qq$-factorial dlt modification obtained as follows. 
\begin{enumerate}
\item $B_{0}$ is log smooth.
 \item $B_1$ is log canonical but not dlt. The weighted blowup of the origin with weights $(1,1,2)$ on $(x,y,z)$ is a $\mathbb{Q}$-factorial dlt modification of $B_1$.
\item For $m\geq 2$, the blow-up of the origin of $B_m$ induces a crepant birational morphism of pairs, and near the non-dlt locus of the exceptional divisor the pull-back pair is of type $B_{m-2}$.
\end{enumerate}
\end{lemma}

The construction of Lemma~\ref{lem:uv-sm-dlt}
also gives the following relative version.

\begin{lemma}
\label{lem:relative-uv-sm-dlt}
Let $(\mathcal X,\mathcal D_1+\mathcal D_2)\to B$
be a family of pairs over a smooth base $B$,
where $\mathcal X$ and $\mathcal D_2$ are smooth
over $B$.
Suppose that, formally along a section of
$\mathcal X\to B$, the boundary is given by
\[
\mathcal D_1=\{uv+y^m=0\},\qquad\text{and}\qquad
\mathcal D_2=\{y=0\},
\]
where $m\ge 1$ is constant.
Then, étale locally over $B$, there is a
projective crepant modification near the section
whose restriction to each fiber is a
$\mathbb Q$-factorial dlt modification.
\end{lemma}

\begin{proof}
Let $\sigma\colon B\to\mathcal X$ be the given section.
For $m\ge2$, we blow up $\sigma(B)$. On the chart
with $u=yu_1$ and $v=yv_1$, the log pull-back of the pair has the same
form with exponent $m-2$ and exceptional divisor $y=0$.
If $m>2$, the relative singular locus of the intersection
of the strict transform of $\mathcal S$ with the exceptional
divisor is the section $u_1=v_1=y=0$. Then, we continue to blow-up this
section until the exponent is zero or one.

If the exponent is one, we take the relative weighted
blow-up of the section with weights $(1,1,2)$. These blow-ups commute with restriction to each
fiber. The calculations in Lemma~\ref{lem:uv-sm-dlt} show
that each step induces a crepant morphism between the corresponding pairs. As a consequence, the resulting
modification on each fiber is $\mathbb{Q}$-factorial and dlt.
\end{proof}

\subsection{Complexity} In this subsection, we recall some concepts and theorems related to the complexity of log pairs and morphisms.
For the most complete forms of the theorems in this direction, we refer the reader to~\cite{BMSZ18,MS21,EFGMSWW25}.

\begin{definition}
{\rm 
Let $f\colon X\rightarrow T$ be a fibration.
Let $(X,B)$ be a log pair and $t\in T$ be a closed point.
The {\em relative complexity} of $(X,B)$ over $t\in T$ is defined to be 
\[
c(X/T,B;t):= \dim(X) + \dim_\qq {\rm Cl}_\qq(X/T) - |B_t|.
\]
Here, the divisor $B_t$ is the largest effective divisor $\leq B$ all whose components intersect $X_t$ non-trivially and $|B_t|$ is the sum of the coefficients of $B_t$.
}
\end{definition} 

\begin{definition}
{\rm 
Let $f\colon X \rightarrow T$ be a fibration.
A pair $(X,B)$ is said to be {\em log Calabi--Yau} over $T$ if it has log canonical singularities and 
$K_X+B\sim_{\qq,T}0$. The {\em index} of $(X,B)$ is the smallest positive integer $i$ for which $i(K_X+B)\sim_T 0$.
}
\end{definition} 

The following theorem is due to the first author and Svaldi (see, e.g.,~\cite[Theorem 1]{MS21}).

\begin{theorem}\label{thm:relative-complexity}
Let $f\colon X \rightarrow T$ be a fibration.
Let $(X,B)$ be a log pair which is log Calabi--Yau over $T$.
Let $t\in T$ be a closed point. 
Then, we have $c(X/T,B;t)\geq 0$.
Furthermore, if $c(X/T,B;t)=0$, then 
$f\colon X\rightarrow T$ is formally toric near $t\in T$.
\end{theorem} 

\begin{remark}
{\rm 
If we set $T={\rm Spec}(\cc)$ in the previous theorem, then we will suppress $T$ and the word {\em relative} from the notation.
In this setting, the theorem says that for a log Calabi--Yau pair $(X,B)$ we have $c(X,B)\geq 0$ and $X$ is toric whenever $c(X,B)=0$.
}
\end{remark}

\subsection{Cluster type pairs}
\label{subsec:cluster-type}
In this subsection, we recall the concept of cluster type varieties and pairs. For a complete survey on cluster type pairs we refer the reader to~\cite{Mor26}.

\begin{definition}[cf.~{\cite[Definition~2.6]{Mor26}}]
{\rm Consider an $n$-dimensional algebraic torus over the complex numbers $\mathbb{T}:={\rm Spec}(\cc[x_1^\pm,\dots,x_n^\pm])$. 
Let $\Omega_n$ be the volume form of the algebraic torus, i.e.,
\[
\Omega_n:=\frac{dx_1\wedge \dots \wedge dx_n}{x_1\cdots x_n}.
\]
An $n$-dimensional normal algebraic variety $X$ is said to be of {\em cluster type} 
if there exists an embedding in codimension one 
$\mathbb{T}\dashrightarrow X$ such that $\Omega_n$ has no zeros on $X$. The effective divisor
$B = \operatorname{Poles}(\Omega_n)$ on $X$ is called a \emph{cluster type boundary}.
The pair $(X,B)$ is called a \emph{cluster type pair}.
}
\end{definition}

The cluster type property behaves well under crepant birational modifications. Furthermore, it can be characterized in terms of dlt modifications and birational contractions to toric pairs. We recall these two statements in the following two lemmas.

\begin{lemma}[cf.~{\cite[Lemma~2.26]{JM24}}]\label{lem:cluster-type-extraction-lcp}
Let $(Y,B_Y)\dashrightarrow (X,B)$ be a crepant birational map that only extracts log canonical places
of the pair $(X,B)$.
Then, we have that $(X,B)$ is of cluster type
if and only if $(Y,B_Y)$ is of cluster type. 
\end{lemma} 

\begin{lemma}[cf.~{\cite[Lemma~2.11]{EFM26}}]\label{lem:cluster-type-MMP}
Let $(X,B)$ be a log Calabi--Yau pair of index one. 
Then, the pair $(X,B)$ is of cluster type
if and only if the following conditions are satisfied: 
\begin{enumerate}
\item there 
exists a $\qq$-factorial
dlt modification $\pi\colon (Y,B_Y) \rightarrow (X,B)$; 
\item there is an effective divisor $E_Y$ not containing any proper dlt strata of $(Y,B_Y)$ and
\item there is a birational contraction 
$\psi\colon Y\dashrightarrow T$ with divisorial exceptional
locus equal to $E_Y$ such that
$T$ is toric and $\psi_*B_Y$ is torus invariant. 
\end{enumerate} 
\end{lemma} 

As a consequence, being of cluster type imposes strong restrictions on the pair $(X, B)$ of index one, and in particular on the singularities of its boundary. The coregularity of a log Calabi--Yau pair is invariant under crepant
birational maps, since the dual complexes of crepant birational dlt pairs are
PL homeomorphic~(cf.~\cite[Proposition~11]{dFKX17}). A complete toric pair has
coregularity zero, because it always has a zero dimensional torus
orbit~(cf.~\cite[Theorem~3.2.6]{CLS11}). Hence coregularity zero is a necessary
condition for a pair to be of cluster type~(cf.~\cite[Theorem~2.7]{Duc24}).
The log Calabi--Yau condition also imposes strong restrictions on the
boundary~$B$. For example, if $(\mathbb P^3,B)$ is a log Calabi--Yau pair of
index one, then $B$ is a reduced quartic surface. By adjunction,
$B$ has at worst semi-log canonical singularities (cf.~\cite[Theorem~4.7]{K13}). The following proposition describes the strictly semi-log
canonical singularities that can appear on an irreducible quartic surface.

\begin{proposition}[cf.~{\cite[Table~1]{SL12}}]
\label{prop:slc_local_models}
Let $S\subset\mathbb P^3$ be an irreducible quartic surface with at worst
semi-log canonical singularities, and let $p\in S$ be a point such that the germ $(S;p)$ is
strictly semi-log canonical and has coregularity $0$ and
multiplicity two. Then one of the following holds.
\begin{enumerate}
  \item $S$ is normal at $p$, and $(S;p)$ is analytically isomorphic
  to the isolated singularity
  \[
    T_{2,q,r}:=\{ (x,y,z) \mid xyz+x^2+y^q+z^r=0\} 
    \quad \text{ with }
    \quad \frac{1}{2}+\frac{1}{q}+\frac{1}{r}<1.
  \]
  \item $S$ is non-normal at $p$, and $(S;p)$ is locally analytically
  equivalent to one of the following:
\begin{align*}
    T_{2,q,\infty}:=&\{(u,v,w)\mid u^2+v^2\bigl(w^2+v^{q-2}\bigr)=0\}
    \text{ with } q\geq 3, &
   T_{2,\infty,\infty}:=&\{(u,v,w)\mid u^2+v^2w^2=0\}.
\end{align*}
\end{enumerate}
\end{proposition}
 
We call a singularity of type $T_{p,q,r}$ a \emph{cusp} and a singularity
of type $T_{2,q,\infty}$ and $T_{2,\infty,\infty}$ a \emph{degenerate cusp}.  Among two-dimensional semi-log canonical hypersurface germs of
multiplicity two, those of coregularity zero are precisely the
cusps and degenerate cusps listed above~(see \cite[\S3.1]{Duc24}).

For a $\mathbb{Q}$-factorial
Fano variety $X$, the cluster type condition of a log Calabi--Yau pair $(X,B)$ can be made more explicit as shown by the following lemma.

\begin{lemma}[cf.~{\cite[Theorem~2.9]{EFM26}}]\label{clustertypeandtorusembedding}
Let $X$ be an $n$-dimensional $\mathbb{Q}$-factorial Fano variety, and let $(X,B)$
be a log Calabi--Yau pair with reduced boundary $B$. Then $(X,B)$ is of cluster type if and only if there
is an open embedding $\mathbb{G}_m^n \hookrightarrow X \setminus B$.
\end{lemma}
 
Thus, for $\mathbb{Q}$-factorial Fano varieties, the cluster type condition can be proved by exhibiting an explicit open torus in the complement of $B$. The following elementary lemma will be used
repeatedly to exhibit algebraic tori in the complement of $B$.

\begin{lemma}\label{prop:shear-torus-chart}
Let $Y$ be an affine variety, set $X=Y\times\mathbb A^1_w$, and let
$a$ and $b$ be regular functions of $Y$.
Set $f:=aw+b$. Then
\[
\Psi:=(\operatorname{pr}_Y,\,f)\colon\
D_X(af)\xrightarrow{\ \sim\ }\ D_Y(a)\times\mathbb G_m
\]
is an isomorphism, with inverse $\Psi^{-1}(s,v)=\bigl(s,(v-b(s))/a(s)\bigr)$.
In particular, if $D_Y(a)\cong\mathbb G_m^2$, then $D_X(af)\cong\mathbb G_m^3$.
\end{lemma}

\begin{proof}
Over $D_Y(a)$, the change of coordinate $\phi (w)=aw+b$
is an automorphism of $D_Y(a)\times\mathbb A^1$, with inverse $w=\frac{f-b}{a}.$ Restricting to $D_{X}(f)$ gives $D_X(af)\cong D_Y(a)\times\G_m.$
\end{proof}

In the setting of Lemma~\ref{prop:shear-torus-chart}, suppose that $Y\cong\mathbb{G}_m^2$ and
$a\in H^0(Y,\mathcal{O}_Y)^{*}$. Then we can find a principal open subset
$D_X(f)\cong\mathbb{G}_m^3$ of $X=\mathbb{G}_m^2\times\mathbb{A}^1$. We set $B:=V(f)$ and we obtain the desired torus embedding in $X\setminus B$. The following
example will be used repeatedly to show the existence of torus embeddings on the complement of quartic surfaces in $\mathbb{P}^3$.

\begin{example}\label{ex:shear-instances}
{\em We explain two particular examples of the previous lemma.
\begin{enumerate}
\item[\textup{(i)}]
Let $X=\mathbb{A}^3$ with variables $x,y,z$.
Let $c\in \cc$ and let $b(y)\in \cc[y]$. 
Set $f:=x^2-(yz+c)^2+b$.
On the open chart $D_X(y)$, we set 
$s=x+yz+c$ and $t=x-yz-c$.
Note that $x=(s+t)/2$ and that $z=(s-t-2c)/(2y)$.
These coordinates identify the open set $D_x(y)$ with 
the product $\mathbb{G}_m\times \mathbb{A}^2_{s,t}$.
In these coordinates we have $f=st+b(y)$. Therefore, 
Lemma~\ref{prop:shear-torus-chart} applies to give
$D_X(fys)\cong\mathbb G_m^3$.

\item[\textup{(ii)}]
Let $X=\mathbb T_{t}\times\mathbb A^2_{x,y}$, let $p\in\cc[t^{\pm1}]^{*}$,
$b\in\cc[t^{\pm1},y]$ and $c\in\cc[t^{\pm1},u]$ with $u$ a new variable, and set
\[
a:=px+b,\qquad\text{and}\qquad f:=ay+c(t,a),
\]
where $c(t,a)$ is obtained from $c(t,u)$ by substituting $a$ for $u$.
Since $p$ is a unit, $a=px+b$ may be used as a coordinate in place of $x$.
Thus $X \cong S\times \mathbb A^1_y$ with $S:=\mathbb T_t\times\mathbb A^1_a$. We put 
$x=\frac{a-b}{p}.$ Since
$D_S(a)=\mathbb T_{t}\times\mathbb T_{a}\cong\mathbb G_m^2$,
Lemma~\ref{prop:shear-torus-chart}, applied with affine coordinate $y$,
coefficient $a$ and constant term $c(t,a)$, yields an isomorphism
\[
\Phi\colon D_X(af)\xrightarrow{\ \sim\ }\mathbb G_m^3 \qquad
(t,x,y)\longmapsto(t,a,f),
\]
whose inverse sends $(t,\alpha,v)$ to the point with
$y=\dfrac{v-c(t,\alpha)}{\alpha}$ and $x=\dfrac{\alpha-b(t,y)}{p(t)}$.
\end{enumerate}
}
\end{example}

\subsection{Obstructions to the cluster type property}
\label{subssec:obstruction-cluster-type}
In this subsection, we detect the failure of the cluster type property using a branched cover. 
For a survey on irrationality of conic bundles, we refer the reader to~\cite{Pro18}. 
For a standard reference on Brauer groups, we refer the reader to~\cite{CTSA21}. 
The following elementary observation is the starting point of this subsection. 

\begin{lemma}\label{lem:clustertypeobstruction}
Let $X$ be a $\mathbb{Q}$-factorial Fano variety and $(X,B)$ a log Calabi--Yau pair of index one. Let $\pi:Y\to X$ be a finite surjective morphism with $Y$ normal and irreducible, such that $\pi$ is étale over $X\setminus B$. Then:
\begin{enumerate}
\item $Y$ is rationally connected, and 
\item if $(X,B)$ is a cluster type pair, then $Y$ is rational.
\end{enumerate}
\end{lemma}

\begin{proof}
Statement (1) is a consequence of the fact that $Y$ is of Fano type~(\cite[Lemma~2.4]{MYY25}). Moreover, all Fano type varieties are rationally connected (cf.~\cite[Theorem~1]{Zha06}). 

For the statement (2), we set $U=X\setminus B$. By the definition of cluster type pair, $U$ contains an algebraic torus $\mathbb{G}_m^n$ as a Zariski dense open subset. Since $\pi$ is \'etale over $X\backslash B$, the induced map $\pi^{-1}(\mathbb G_m^n)\to \mathbb G_m^n$ is also finite \'etale. A connected finite \'etale cover of an algebraic torus is again an algebraic torus~(\cite[Lemma~1]{Miyanishi72}). Therefore, every connected component of $\pi^{-1}(\mathbb G_m^n)$ is isomorphic to $\mathbb{G}_m^n$. Because $\pi^{-1}(\mathbb G_m^n)$ is an open subvariety of the irreducible variety $Y$, any such connected component is dense and open in $Y$. Thus, $Y$ is birational to $\mathbb{G}_m^n$ and, consequently, rational.
\end{proof}

To prove that $(\mathbb{P}^3, S)$ is not of cluster type, it suffices to show
that the cyclic quartic cover of $\mathbb{P}^3$ branched along $S$ is irrational.
We therefore recall an irrationality criterion for threefolds admitting a conic
bundle structure, namely a slightly simplified form of Shokurov's
criterion~\cite{Sho83}. Our formulation differs from the one in the
literature, but it follows immediately from it.

\begin{definition}[cf.~{\cite[Definition~3.1]{Pro18}}]
\em Let $S$ be a smooth projective surface. 
\begin{enumerate}
\item A \emph{conic bundle} is a proper flat morphism $\pi\colon X\to S$ from a smooth projective threefold $X$ such that every scheme theoretic fiber is a plane conic.
\item The \emph{discriminant curve} $\Gamma$ is the reduced divisorial locus over which $\pi$ is not smooth.

\item  The conic bundle $\pi$ is called \emph{standard} if $\pi^{-1}(D)$ is irreducible for every
prime divisor $D\subset S$.
\end{enumerate}
\end{definition}

\begin{proposition}[cf.~{\cite[Theorem~7.5]{Pro18}}]
\label{Thm:Shokurov}
Let $\pi\colon X\to S$ be a standard conic bundle over a smooth projective
rational surface $S$, and assume that its discriminant curve
$\Gamma\subset S$ is smooth and irreducible. If $\Gamma$ is neither
hyperelliptic nor trigonal and is not birational to a plane quintic, then
$X$ is irrational.
\end{proposition}
In the rest of this subsection, we give a criterion for applying
Proposition~\ref{Thm:Shokurov} to the conic bundles that occur below.
In our applications, both the discriminant curve and the associated
double cover can be read directly from the polynomial defining the
generic conic.

\begin{lemma}[{cf.\ \cite[Definition~1.1.11 and Proposition~7.2.4(b)]{CTSA21}}]
\label{lem:conic-brauer-discriminant}
Let $S$ be a rational surface and $K:=\mathbb{C}(S)$ be the function field of $S$.
\begin{enumerate}
\item Let $a,b,c\in K$ and set $q(y,z)=ay^2+2byz+cz^2$.
If $b^2-ac\ne0$, then
\[
Q_\eta:=\{x^2-q(y,z)=0\}\subset\mathbb P^2_K
\]
is a smooth conic. On any open subset of $S$ on which
$a,b,c$ are regular, the family defined by this equation
has singular fibers precisely over $V(b^2-ac)$.

\item For $a,b \in K$, the smooth conic
\[
C:=\{ x^2-a y^2-bz^2=0\} 
\]
over $\mathbb C(S)$ corresponds to the quaternion algebra $(a,b)$, whose class
lies in $\operatorname{Br}(\mathbb C(S))[2]$.
\end{enumerate}
\end{lemma}
Let $\alpha=(\lambda,\mu)$ with $\lambda,\mu\in K^\times$,
and let $D\subset S$ be a prime divisor.
The residue of $\alpha$ along $D$ is
\begin{equation}\label{residueformula}
\partial_D(\alpha)
:=
\left[
(-1)^{\nu_D(\lambda)\nu_D(\mu)}
\frac{\lambda^{\nu_D(\mu)}}{\mu^{\nu_D(\lambda)}}
\right]
\in
\mathbb C(D)^*/\mathbb C(D)^{*2}
\end{equation}
$\nu_D$ is the divisorial valuation of $D$
(see \cite[Definition~1.4.3 and (1.16)--(1.18)]{CTSA21}). 
The expression in brackets is determined modulo $\mathbb C(D)^{*2}$. If $\partial_D(\alpha)$ is nontrivial, then choose $r_D\in \mathbb C(D)$ such that $r_D\equiv \partial_D(\alpha)$ (mod $\mathbb C(D)^{*2}$).  The element $r_D$ determines the quadratic field extension $\mathbb C(D)(\sqrt{r_D})/\mathbb C(D).$
 Thus, if $\Gamma$ is an irreducible component of the discriminant and
$\nu\colon C\to\Gamma$ is its normalization, the pullback of
$\partial_\Gamma(\alpha)$ to $\mathbb C(C)$ determines a double cover $\tau\colon \widetilde C\to C$. We say that this cover is induced by the residue $\partial_\Gamma(\alpha)$.

\begin{lemma}
\label{lem:standardization-after-resolution}
Let $Q$ be a smooth conic over the function field of a smooth
projective rational surface $S_0$. Suppose that its Brauer
class $\alpha$ ramifies along a unique prime divisor $\Gamma$,
and let $C$ be the normalization of $\Gamma$.
If the residue double cover $\tau\colon\widetilde C\to C$
is \'etale, then there exist a sequence of point blow-ups
$S\to S_0$ and a standard conic bundle $\pi\colon X\to S$
whose generic fiber and discriminant curve are isomorphic
to $Q$ and $C$, respectively.
\end{lemma}

\begin{proof}
Choose a sequence of point blow-ups
$\sigma\colon S\to S_0$ such that the strict transform
$\Gamma'$ of $\Gamma$ is smooth. Then $\Gamma'\cong C$.
By \cite[Proposition~3.10]{Pro18}, there exists a standard
conic bundle $\pi\colon X\to S$ whose discriminant is
$\Gamma'$ and whose associated double cover is
$\widetilde C\to\Gamma'$.

Let $\alpha'\in\operatorname{Br}(K)[2]$ be the Brauer
class of the generic fiber of $\pi$, identifying
$\mathbb C(S)$ with $K$ via $\sigma$.
A prime divisor on $S_0$ and its strict transform on $S$
define the same valuation of $K$. Thus the residues of
$\alpha$ and $\alpha'$ on $S_0$ agree: they are trivial
away from $\Gamma$, and along $\Gamma$ they both define
$\tau$. As in the proof of
\cite[Proposition~3.10]{Pro18}, since
$\operatorname{Br}(S_0)=0$, a Brauer class over $K$
is uniquely determined by its residues on $S_0$.
Hence $\alpha=\alpha'$, so the generic fiber of
$\pi$ is isomorphic to $Q$ over $K$.
\end{proof}

The following is the main technical proposition that we will use to prove that conic bundles are irrational.

\begin{proposition}\label{prop:technical-curve-conditions}
Let $R$ and $H\in \mathbb{C}[x]$ be non-constant coprime polynomials of degree $r$ and $h$ respectively. Assume that $r\neq h$ and neither $R$ nor $H$ is a square. Let $c\in \mathbb{C}^*$ and consider the Brauer class
\[
\alpha
=
\bigl(R(x),\,cR(x)u^4+H(x)\bigr)
\in \operatorname{Br}(\mathbb C(x,u))[2].
\]
Put $f(x):=-\frac{H(x)}{cR(x)}.$
Let $\Gamma$ be the closure in
$\mathbb P^1_x\times\mathbb P^1_u$ of the affine curve $u^4=f(x)$ and let $C$ be its normalization. The projection to the first factor
induces a cyclic morphism $\rho\colon C\rightarrow\mathbb P^1_x$ of degree $4$. Let $B_4$ and $B_2$ be the numbers of branch points of $\rho$ with
ramification indices $4$ and $2$ respectively.
Then, we have the following statements
\begin{enumerate}
\item
$\Gamma$ is irreducible and $\partial_\Gamma(\alpha)
= [R]\in\mathbb C(C)^*/\mathbb C(C)^{*2}.$
The corresponding double cover $\tau: \tilde{C}\rightarrow C$ is connected.
\item
The double cover $\tau$ is \'etale if and only if $r\equiv0\pmod2$ or 
$4\nmid(r-h).$
\item
The genus of $C$ is $g(C)=\frac{3B_4+2B_2-6}{2}.$
If $B_4\geq4$ and $g(C)\geq4$, then $C$ is not hyperelliptic.
If $g(C)\geq7$, then $C$ is neither trigonal nor birational to a
plane quintic.
\end{enumerate} 
\end{proposition}

\begin{proof}
Note that $R$ and $H$ are coprime polynomials.
Furthermore, neither of them is a perfect square. 
Hence, the quotient $f=-H/(cR)$ is not a Square. 
Therefore, the curve $\Gamma$ is irreducible 
and $\mathbb{C}(C)\cong \mathbb{C}(x)(\sqrt[4]{f})$ is a cyclic extension of degree $4$ of $\mathbb{C}(x)$.
We turn to compute the residue $\delta_\Gamma(\alpha)$. 
Away from $\Gamma$ the class $\alpha$ is unramified.
Thus, we only need to check the divisors $\{R=0\}$ 
and the two lines at infinity. 
Along $\{u=\infty\}$ the residue is $[R^{-4}]=1$. 
Along $\{x=\infty\}$ the leading coefficient in $x$ is a cosntant or a multiple of $u^4$. This is where the fact that $r\neq h$ is used. 
On $\Gamma$, the polynomial $R$ is a unit and $cRu^4+H$ has a simple zero. Therefore, the formula~\eqref{residueformula} gives $\partial_\Gamma(\alpha)=[R]$.
Moreover, $\widetilde{C}$ is obtained by adjoining $\sqrt{R}$. 
If $\widetilde{C}$ was disconnected, then $\sqrt{R}$ would like in $\mathbb{C}(C)$. 
A cyclic quartic extension has only one quadratic subfield.
Therefore, we would have $\mathbb C(x)(\sqrt R) \cong \mathbb C(x)(\sqrt f)$ and $R/f=-cR^2/H$ is a square.
This would imply that $H$ is a square leading to a contradiction. This implies that $\widetilde{C}$ is connected. This proves (1).

Now, we prove (2).
We next determine when the residue cover is \'etale. For
$Q\in C$, the quadratic extension
$\mathbb C(C)(\sqrt R)/\mathbb C(C)$ is ramified at $Q$ if and only if
$\ord_Q(R)$ is odd. If $Q$ lies over
$p\in\mathbb P^1_x$, then
\[
e(Q/p)
=
\frac{4}{\gcd(4,\ord_p(f))}
\]
and therefore $\ord_Q(R)
=
e(Q/p)\ord_p(R)$
~(see \cite[Proposition~3.7.3]{stichenoth09}). Let $p$ be a finite zero of $R$ of multiplicity $m$.
Since $R$ and $H$ are coprime, $\ord_p(f)=-m.$
Thus $\ord_Q(R)
=
\frac{4m}{\gcd(4,m)},$ which is even. At every other finite point,
$\ord_p(R)=0$. Hence the residue cover is unramified over
$\mathbb A^1_x$.

It remains to consider the points over $x=\{\infty\}$. We have $\ord_\infty(R)=-r$ and  $\ord_\infty(f)=r-h.$ Thus, for every $Q\in C$ lying over infinity, $\ord_Q(R)
=
-\frac{4r}{\gcd(4,r-h)}.$ This integer is odd precisely when $r$ is odd and
$4\mid(r-h)$. Hence the residue cover is \'etale if and only if $r\equiv0\pmod2$ or $4\nmid(r-h).$
This proves~(2).

For~(3), Riemann--Hurwitz gives $2g(C)-2=-8+3B_4+2B_2$ and hence $g(C)=\frac{3B_4+2B_2-6}{2}.$

Suppose that $C$ is hyperelliptic and let
$q\colon C\to\mathbb P^1$ be its degree $2$ map.
If $(\rho,q)\colon C\rightarrow\mathbb P^1\times\mathbb P^1$ is birational onto its image, the Castelnuovo--Severi inequality gives $g(C)\le(4-1)(2-1)=3$ (see \cite[Chapter~VIII, Exercise~C-1]{ACGH85}).  
Thus, if $g(C)\geq4$, then the maps $\rho$ and $q$ factor through a common
curve $C\to D$. Since its degree divides both $4$ and $2$, this map
has degree $2$, and hence $D\cong\mathbb P^1$. The resulting double
cover $D\rightarrow\mathbb P^1_x$ is the quadratic intermediate cover of $\rho$ and is branched exactly at the
$B_4$ points where $\rho$ is totally ramified. Riemann--Hurwitz then
gives $B_4=2$. Therefore $B_4\geq4$ excludes the hyperelliptic case.

If $C$ is trigonal, let $f\colon C\to\mathbb P^1$
be a morphism of degree $3$.
Since $\deg\rho=4$, the map
$(\rho,f)\colon C\to\mathbb P^1\times\mathbb P^1$
is birational onto its image.
Thus, Castelnuovo--Severi gives
$g(C)\le(4-1)(3-1)=6$.
If $C$ is birational to a plane quintic,
the same bound follows from its arithmetic genus.
This proves~(3).
\end{proof}

We will culminate this section by presenting the main proposition that will be used to prove the irrationality of the quartic threefolds considered below.

\begin{proposition}\label{prop:fractional-kummer-irrationality}
Let $R,H\in \mathbb{C}[x]$ be non-constant coprime polynomials of degree $r$ and $h$, respectively. 
Assume that $r\neq h$ and neither $R$ nor $H$ is a square. Let $c\in \mathbb{C}^*$ and consider the Brauer class 
\[
\alpha
=
\bigl(R(x),\,cR(x)u^4+H(x)\bigr)
\in \operatorname{Br}(\mathbb C(x,u))[2].
\]
Put $f(x):= -H(x)/cR(x)$.
Let $\Gamma$ be the closure in $\mathbb{P}^1_x\times \mathbb{P}^1_u$ of the affine curve $u^4=f(x)$ and let $C$ be its normalization. The projection to the first factor induces a finite morphism $\rho\colon C\rightarrow \pp^1_x$ of degree $4$. 
Let $B_2$ and $B_4$ be the number of branch points of $\rho$ with ramification indices $4$ and $2$ respectively. Let $Q_\eta$ be the conic over $\mathbb C(x,u)$ with Brauer
class $\alpha$. Let $X$ be a complex variety with function field $\mathbb C(x,u)(Q_\eta)$. If the following conditions hold
\begin{enumerate}
\item[(i)]
$r\equiv0\pmod2$ or $4\nmid(r-h)$, and
\item[(ii)]
$B_4\geq4$ and $3B_4+2B_2\geq20$,
\end{enumerate}
then $X$ is irrational.
\end{proposition}

\begin{proof}
By the proof of Proposition~\ref{prop:technical-curve-conditions}, we know that $\Gamma$ is the only ramification divisor of $\alpha$ on 
$\pp^1_x\times \pp^1_u$.
By Proposition~\ref{prop:technical-curve-conditions}.(1), Proposition~\ref{prop:technical-curve-conditions}.(2), and Lemma~\ref{lem:standardization-after-resolution}, we know that $X$ is birational to a smooth projective threefold that admits a standard conic bundle to a surface whose discriminant is isomorphic to $C$. 
Conditions (i) and (ii) together with Proposition~\ref{prop:technical-curve-conditions}.(3) then implies that $C$ is neither hyperelliptic nor trigonal nor birational to a plane quintic. 
Thus, we may apply Theorem~\ref{Thm:Shokurov} to conclude that $X$ is irrational.
\end{proof}

\section{Cluster type condition in families} 
In this section, we prove the main theorem regarding the cluster type condition in families of Calabi--Yau pairs on terminal Fano varieties. We will use the following version of invariance of plurigenera. 

\begin{theorem}\label{thm:invariance-better}
Let $f\colon (\mathcal{Z},\Theta)\rightarrow C$ be a projective family of pairs over a smooth curve $C$.
Let $c\in C$ be a closed point.
Assume that $f$ has $\qq$-factorial terminal log fibers $(\mathcal{Z}_t,\Theta_t)$ for $t$ near $c\in C$. 
Assume that $\Theta$ is big
and horizontal over $C$.
Assume that no component of $\Theta_c$ 
is contained in the stable base locus ${\rm SBs}(K_{\mathcal{Z}_c}+\Theta_c)$. 
Then, sections of sufficiently divisible multiples of 
$K_{\mathcal{Z}_c}+\Theta_c$
extend to sections of the corresponding multiples of
$K_{\mathcal{Z}}+\Theta$ near $c\in C$.
\end{theorem}

The proof of the previous theorem is identical to that of~\cite[Theorem 4.2]{HMX13}. In that statement, it is assumed that the restrictions $\Theta_c$ of the components of $\Theta$ to fibers near $c\in C$ remain irreducible. In such proof, the authors extend to $\mathcal{Z}$ the boundary obtained from substracting ${\rm SBs}(K_{\mathcal{Z}_c}+\Theta_c)$ from $\Theta_c$.
Our assumption makes this overlap zero, so the original boundary already supplies the required extension. The remaining MMP and vanishing arguments from~\cite[Theorem 4.2]{HMX13} apply verbatim. 

\begin{proof}[Proof of Theorem~\ref{thm:cluster-type-in-families}]
We let $\{0\}\in T$ and $\{1\}\in T$ be two closed points. We want to prove that $(\mathcal{X}_0,\mathcal{B}_0)$
is of cluster type if and only if $(\mathcal{X}_1,\mathcal{B}_1)$ is of cluster type.
From now on, we assume that $(\mathcal{X}_1,\mathcal{B}_1)$ is a cluster type pair. 
First, we reduce to the case in which the base is a smooth curve. Indeed, if $\{0\}\in T$ and $\{1\}\in T$ are two closed points, then they can be connected via a smooth curve $C$. 
This follows from Bertini's theorem and Condition~\eqref{cond:normal-connected-base} of the theorem.
The pull-back family $(\mathcal{X}_C,\mathcal{B}_C)\rightarrow C$ satisfies the same conditions as stated in the theorem; its fibers are terminal Fano varieties, 
the log fibers are log Calabi--Yau of index one, and it
admits a simultaneous $\mathbb{Q}$-factorial dlt modification.
We argue that $(\mathcal{X}_C,\mathcal{B}_C)\rightarrow C$ has isomorphic class group restrictions. Indeed, the structure sheaf is acyclic on the fibers of this family, so the relative Picard functor is discrete and separated. Henceforth, a Cartier class which is trivial on one fiber is therefore relative trivial over the connected curve. 
Therefore, if we prove the statement for the family over the curve, then $(\mathcal{X}_{0},\mathcal{B}_{0})$ is of cluster type
if and only if 
$(\mathcal{X}_{1},\mathcal{B}_{1})$ is of cluster type.  
Thus, it suffices to prove the statement for $T$ a smooth curve.
From now on, we assume that $T$ is itself a smooth curve.

Second, we argue that the cluster type condition 
for the log pairs $(\mathcal{X}_t,\mathcal{B}_t)$ is constant over the curve $T$.
Let $\{0\}$ and $\{1\}$ be two closed points of $T$.
We aim to show that the pair $(\mathcal{X}_{0},\mathcal{B}_{0})$ is also of cluster type. 
Since $(\mathcal{X},\mathcal{B})$ admits a simultaneous $\qq$-factorial dlt modification (see Condition~\eqref{cond:simult-q-dlt}), 
we know that every log canonical center of
$(\mathcal{X},\mathcal{B})$ is horizontal over $T$. 
In particular, the pair
$(\mathcal{X},\mathcal{X}_{0}+\mathcal{X}_{1})$ is a purely terminal Fano pair over $T$ by inversion of adjunction. Here, we are using Condition~\eqref{cond:TF}.
Thus, we may find a boundary divisor $\Delta$ such that $(\mathcal{X},\mathcal{X}_{0}+\mathcal{X}_{1}+\Delta)$
is a purely terminal Calabi--Yau pair over $T$. 
By construction, the effective divisor $\Delta$ is ample over $T$. 
Let $\psi\colon (\mathcal{Y},\mathcal{B}_{\mathcal{Y}})\rightarrow 
(\mathcal{X},\mathcal{B})$ be a simultaneous $\qq$-factorial dlt modification. We know that $K_{\mathcal{X}_1}+\mathcal{B}_1\sim 0$ and so every component of $\mathcal{B}_{\mathcal{Y}_1}$ is reduced. 
By the simultaneous dlt modification condition, every component of $\mathcal{B}_{\mathcal{Y}_1}$ is the restriction of a component of $\mathcal{B}_{\mathcal{Y}}$. Therefore, $\mathcal{B}$ and $\mathcal{B}_{\mathcal{Y}}$ are integral. Condition~\eqref{cond:CL} implies that $K_{\mathcal{X}/T}+\mathcal{B}\sim_T 0$ and $K_{\mathcal{Y}/T}+\mathcal{B}_{\mathcal{Y}}\sim_T 0$.
We can write 
\[
\psi^*(K_{\mathcal{X}}+\mathcal{B}+\mathcal{X}_{0}+\mathcal{X}_{1})=
K_{\mathcal{Y}}+\mathcal{B}_{\mathcal{Y}}+\mathcal{Y}_{0}+\mathcal{Y}_{1}, 
\]
where $\mathcal{Y}_{i}$ is the strict transform
in $\mathcal{Y}$ of $\mathcal{X}_{i}$, for $i\in \{0,1\}$. 
Note that for each $t\in T$, the crepant birational map 
\[
\psi_t\colon (\mathcal{Y}_t,\mathcal{B}_{\mathcal{Y}_t}) 
\rightarrow 
(\mathcal{X}_t,\mathcal{B}_t) 
\]
only extracts log canonical places of the pair
$(\mathcal{X}_t,\mathcal{B}_t)$.
Therefore, for $\{1\}\in T$, the 
pair $(\mathcal{Y}_{1},\mathcal{B}_{\mathcal{Y}_{1}})$
is of cluster type due to Lemma~\ref{lem:cluster-type-extraction-lcp}.
Therefore, by Lemma~\ref{lem:cluster-type-MMP}, we have: 
\begin{itemize} 
\item a $\qq$-factorial dlt modification
$\alpha_{1}\colon (\mathcal{Z}_{1},\mathcal{B}_{\mathcal{Z}_{1}}) 
\rightarrow (\mathcal{Y}_{1},\mathcal{B}_{\mathcal{Y}_{1}})$; 
\item an effective divisor $E_{1}$ in $\mathcal{Z}_{1}$
which contains no proper dlt strata
of $(\mathcal{Z}_{1},\mathcal{B}_{\mathcal{Z}_{1}})$; and 
\item a birational contraction $\beta_{1}\colon \mathcal{Z}_{1}\dashrightarrow T_{1}$ with divisorial exceptional locus equal to $E_{1}$ such that 
$T_1$ is toric and ${\beta_{1}}_*\mathcal{B}_{\mathcal{Z}_{1}}$ is torus invariant.
\end{itemize} 
After passing to a higher dlt modification, we may assume that $\alpha_1$ is obtained, up to a small birational modification, by finitely many formally toric weighted blowups of strata.
We can find a simultaneous $\qq$-factorial dlt modification 
$\alpha\colon (\mathcal{Z},\mathcal{B}_{\mathcal{Z}})\rightarrow (\mathcal{Y},\mathcal{B}_{\mathcal{Y}})$ that restricts to $\alpha_{1}$, in codimension one, for the closed point $\{1\}\in T$. 
Indeed, any dlt strata of $(\mathcal{Y}_{1},\mathcal{B}_{\mathcal{Y}_{1}})$ is the restriction to $\mathcal{Y}_{1}$ of a dlt strata of 
$(\mathcal{Y},\mathcal{B}_{\mathcal{Y}})$.
Thus, the extension statement follows from~\cite[Lemma 36 and Proposition 37]{dFKX17}.
We write $\rho$ for the induced morphism from $\mathcal{Z}$ to $\mathcal{X}$.
Moreover Condition~\eqref{cond:CL} is preserved under simultaneous $\mathbb Q$-factorial dlt modifications. This follows from the restriction exact sequence applied to the open subset on which the blow-up is an isomorphism.
Hence for every closed point $t\in T$
the homomorphism $\operatorname{Cl}(\mathcal{Z}/T)\rightarrow\operatorname{Cl}(\mathcal{Z}_t)$ is an isomorphism.
Consider the boundary divisor 
\[
\Delta_{\epsilon}:=\epsilon  \Delta + (1-\epsilon) \mathcal{B} 
\]
where $\epsilon \in (0,1)$.
By Condition~\eqref{cond:CY}, every exceptional non-terminal valuation of
the pair $(\mathcal{X},\Delta_\epsilon)$ 
is a log canonical place of the pair $(\mathcal{X},\mathcal{B})$, 
so they are horizontal over $T$.
On the simultaneous dlt modification these places are monomial along strata, with discrepancies computed by the same horizontal boundary.
For $\epsilon_0 \in (0,1)$ small enough, only finitely many monomial valuations have discrepancy at most one. 
Thus, for $\epsilon_0 \in (0,1)$ small enough, 
possibly replacing $(\mathcal{Z},\mathcal{B}_{\mathcal{Z}})$ with a higher dlt model, we may assume that every exceptional non-terminal valuation of $(\mathcal{X},\Delta_{\epsilon_0})$ is extracted by $\alpha \circ \psi$. 
Furthermore, perturbing $\epsilon_0$ if possible, we may assume that every divisor extracted by $\psi$ is a exceptional non-terminal valuation of $(\mathcal{X},\Delta_{\epsilon_0})$. 
Therefore, if we write 
\[
\rho^*(K_\mathcal{X}+\mathcal{X}_{0}+\mathcal{X}_{1}+\Delta_{\epsilon_0})=
K_{\mathcal{Z}}+\mathcal{Z}_{0}+\mathcal{Z}_{1}+\Delta_{\mathcal{Z},\epsilon_0}, 
\]
then the following properties are satisfied: 
\begin{enumerate}
\item[(i)]\label{cond:boundary} the divisor $\Delta_{\mathcal{Z},\epsilon_0}$ is effective
and $\Delta_{\mathcal{Z},\epsilon_0}\geq c\rho^*\Delta$ for some $c>0$; 
\item[(ii)]\label{cond:terminal} the pair $(\mathcal{Z},\Delta_{\mathcal{Z},\epsilon_0})$ is a terminal Calabi--Yau pair over $T$; and 
\item[(iii)]\label{cond:purely-terminal} the pair 
$(\mathcal{Z},\mathcal{Z}_{0}+\mathcal{Z}_{1}+\Delta_{\mathcal{Z},\epsilon_0})$ is a purely terminal Calabi--Yau pair over $T$.
\end{enumerate} 
Let $P_{1} \leq E_{1}$ be a prime component of $E_{1}$.
From the isomorphism condition
${\rm Cl}(\mathcal{Z}/T)\rightarrow {\rm Cl}(\mathcal{Z}_{1})$ we conclude that there exists a Weil divisor $D$ on $\mathcal{Z}$ such that $D$ does not contain $\mathcal{Z}_{1}$ on its support 
and $D|_{\mathcal{Z}_{1}}\sim_\qq E_{1}$. 
Note that the divisor $D$, a priori, is not $\qq$-Cartier and  need not be effective.

However, by~\cite[Corollary 3.2]{dfH09b}, we conclude that $\mathcal{Z}$ is $\qq$-factorial. 
Therefore, $D$ is a $\qq$-Cartier divisor.
Nevertheless it can have some negative coefficients.
Replacing $P_1$ with a multiple, we may assume that $D$ itself is Cartier and we can write $D\sim H_1-H_2$ where both $H_1$ and $H_2$ are ample. Choose $G\in |N\rho^*\Delta-H_2|$ where $N\gg 0$ general whose restriction to $\mathcal{Z}_1$ does not contain the support of $P_1$. Pick $0<\lambda\ll 1$ and $0< \delta \ll \lambda cN^{-1}$ and define
\[
\Theta:=(\Delta_{\mathcal{Z},\epsilon_0}-c\rho^*\Delta)
+(1-\lambda)c\rho^*\Delta + \frac{\lambda c}{N}(H_2+G)-\delta(H_2-H_1).
\]
By construction, the pair $(\mathcal{Z},\Theta)$ is terminal, $\Theta$ is big, $K_{\mathcal{Z}}+\Theta \sim_\qq \delta D$,
and $\Theta_1$ does not contain $P_1$. 
Note that ${\rm SBs}(P_1)=P_1$ as $P_1$ is exceptional for a birational contraction. Therefore, we can apply Theorem~\ref{thm:invariance-better}, to conclude that
\[
H^0(\mathcal{Z},MD) \rightarrow H^0(\mathcal{Z}_{t_1},MD|_{\mathcal{Z}_{t_1}})
\]
is surjective for $M$ divisible enough.
Thus, we can find an effective $\qq$-divisor $P_{\mathcal{Z}}$ on $\mathcal{Z}$ that restricts to $P_{1}$ on $\mathcal{Z}_{1}$.
Proceeding component by component, we may assume that there exists an effective $\qq$-divisor
$E_{\mathcal{Z}}$ in $\mathcal{Z}$ that restricts to $E_{1}$. 
Furthermore, we may assume that the number of components of $E_{\mathcal{Z}}$ is at least the number of components of $E_{1}$.
Applying invariance of plurigenera again we can see that $E_{\mathcal{Z}}$ must have relative Kodaira dimension zero over $T$.
Observe that $\mathcal{Z}/T$ is of Fano type by~\cite[Lemma 2.17]{Mor24a}.
Therefore, we can run a $E_{\mathcal{Z}}$-MMP over $T$ that terminates after contracting all the components of $E_{\mathcal{Z}}$.
Let $\beta\colon \mathcal{Z}\dashrightarrow \mathcal{W}$ be such birational contraction. 
We let $\mathcal{B}_{\mathcal{W}}$ be the push-forward of $\mathcal{B}_{\mathcal{Z}}$ to $\mathcal{W}$. 
Note that the number of components $|E_{1}|$ of 
$E_{1}$ equals $c(\mathcal{Z}_{1},\mathcal{B}_{\mathcal{Z}_{1}})$. 
In particular, we conclude that 
\begin{equation}\label{inequality:comp} 
|E_{\mathcal{Z}}|\geq 
\dim \mathcal{Z} + \rho(\mathcal{Z})-|\mathcal{B}_{\mathcal{Z}}|-1. 
\end{equation} 

The previous inequality implies that the complexity of $(\mathcal{Z},\mathcal{B}_{\mathcal{Z}})$ over any closed point of $T$ is at most $|E_{\mathcal{Z}}|+1$. Note that by~\eqref{inequality:comp}
the number of divisors contracted by $\mathcal{Z}\rightarrow \mathcal{W}$ is at least $\dim \mathcal{Z}+\rho(\mathcal{Z})-|\mathcal{B}_{\mathcal{Z}}|-1$. 
On the other hand, the divisor $E_{\mathcal{Z}}$ contains no component of $\mathcal{B}_{\mathcal{Z}}$, indeed this happens on a closed fiber and so it happens over an open subset of the base.
Therefore, the crepant birational map $(\mathcal{Z},\mathcal{B}_{\mathcal{Z}})\dashrightarrow (\mathcal{W},\mathcal{B}_{\mathcal{W}})$ decreases the complexity of the underlying pairs by at least $|E_{\mathcal{Z}}|$.
Thus, the pair $(\mathcal{W},\mathcal{B}_{\mathcal{W}})$ which is log Calabi--Yau over $T$ has complexity at most one over any closed point. Therefore, we deduce that 
the pair 
\begin{equation}\label{w-pair}
(\mathcal{W},\mathcal{B}_{\mathcal{W}}+\mathcal{W}_{0}+\mathcal{W}_{1})
\end{equation} 
is log Calabi--Yau over $T$
and it has complexity zero near $\{0\}$ and $\{1\}$.
We conclude that the log Calabi--Yau pairs 
$(\mathcal{W}_{i},\mathcal{B}_{\mathcal{W}_{i}})$,
induced by adjunction of~\eqref{w-pair}
to the fibers over $i\in \{0,1\}$
are normal toric varieties. 
Indeed, this follows from Theorem~\ref{thm:relative-complexity}.
Therefore, we have a diagram 
\[
\xymatrix{
(\mathcal{Z}_{0},\mathcal{B}_{\mathcal{Z}_{0}})\ar[d]_-{\rho_0}\ar@{-->}[r]^-{\beta_{0}} & (\mathcal{W}_{0},\mathcal{B}_{\mathcal{W}_{0}})\\ 
(\mathcal{X}_{0},\mathcal{B}_{0})} 
\]
where $\rho_{0}$ is a $\qq$-factorial dlt modification and $\beta_{0}$ is a crepant birational map onto a toric log Calabi--Yau pair. 
By Lemma~\ref{lem:cluster-type-MMP}, in order to conclude that $(\mathcal{X}_{0},\mathcal{B}_{0})$ is of cluster type, it suffices to prove that $\beta_{0}$ is a birational contraction, i.e., it does not extract divisors. 

By Condition~\eqref{cond:purely-terminal}, the pair 
$(\mathcal{Z},\mathcal{Z}_{0}+\mathcal{Z}_{1}+\Delta_{\mathcal{Z},\epsilon_0})$ is purely terminal and log Calabi--Yau over $T$.
Therefore, the birational contraction $\beta$ induces a crepant birational map 
\[
\beta\colon 
(\mathcal{Z},\mathcal{Z}_{0}+\mathcal{Z}_{1}+\Delta_{\mathcal{Z},\epsilon_0})
\dashrightarrow 
(\mathcal{W},\mathcal{W}_{0}+\mathcal{W}_{1}+\Delta_{\mathcal{W},\epsilon_0})
\]
where $\Delta_{\mathcal{W},\epsilon_0}$ is the push-forward of $\Delta_{\mathcal{Z},\epsilon_0}$ to $\mathcal{W}$. 
By adjunction, we get a crepant birational map
between Calabi--Yau pairs
\[
\beta_{0}\colon 
(\mathcal{Z}_{0},\Delta_{\mathcal{Z}_{0},\epsilon_0})\dashrightarrow  
(\mathcal{W}_{0},\Delta_{\mathcal{W}_{0},\epsilon_0}),
\]
such that the pair on the domain is terminal.
Moreover, both $\Delta_{\mathcal{Z}_{0},\epsilon_0}$ and $\Delta_{\mathcal{W}_{0},\epsilon_0}$ are effective divisors as they are defined via adjunction. 
Therefore, we conclude that $\beta_{0}$ does not extract divisors, otherwise $\Delta_{\mathcal{W}_{0},\epsilon_0}$ would have some negative coefficients. 
We conclude that $\beta_{0}$ is a birational contraction and so by Lemma~\ref{lem:cluster-type-MMP}, we conclude that the pair $(\mathcal{X}_{0},\mathcal{B}_{0})$ is of cluster type. 
This finishes the proof of the theorem.
\end{proof}

\begin{remark}\label{rem:analytic-base}
{\em 
Theorem~\ref{thm:cluster-type-in-families}, and hence Corollary~\ref{cor:p3-families}, are also valid for  projective flat families of complex analytic pairs over smooth irreducible complex bases. Indeed, the proof is local and reduces to projective families over complex disks. After shrinking to suitable Stein neighborhoods, Fujino's analytic minimal model program and vanishing theorems provide the analytic counterparts of the algebraic theorems used in the proof. See references for MMP~\cite[Theorem 1.7, Theorem 1.8, and Theorem 5.2]{Fuj22} and references for vanishing ~\cite[Theorem 1.1.(ii) and Theorem 1.2]{Fuj25}).
}
\end{remark} 
\begin{proof}[Proof of Corollary~\ref{cor:p3-families}]
As $T$ is smooth, we can cover it with complex disks $\{ \mathbb{D}_i\}_{i \in I}$ where $I$ is a countable set.
It suffices to show that the cluster type property is constant over each disk as $T$ is connected. 
Without loss of generality, we may replace $T$ with one of these disks.
Thus, we may assume that $(\pp^3,B_t)$ admit a simultaneous $\qq$-factorial dlt modification over $T$.
We consider the constant family $\mathcal{X}:=\pp^3\times T\rightarrow T$ with the boundary $\mathcal{B}$ that satisfies $\mathcal{B}_t:=B_t\times \{t\}$ for every $t\in T$.
Note that $\mathcal{X}_t$ is a terminal Fano threefold
and $(\mathcal{X}_t,\mathcal{B}_t)\cong (\pp^3,B_t)$ is a log Calabi--Yau pair of index one. 
Condition~\eqref{cond:CL} of Theorem~\ref{thm:cluster-type-in-families} holds trivially as we are considering a product family. 
By the reduction, the pairs $(\pp^3,B_t)$ admit a simultaneous $\mathbb Q$-factorial dlt modification. By Theorem~\ref{thm:cluster-type-in-families} and Remark~\ref{rem:analytic-base}, we conclude that
the cluster type property of $(\pp^3,B_t)$ is independent of $t\in T$.
\end{proof}

For three dimensional Calabi--Yau pairs, the irrationality obstruction in Lemma~\ref{lem:clustertypeobstruction} extends to a very general fiber. This follows from
the deformation result of de Fernex and Fusi~\cite[Theorem~1.2]{dFF13}. 
We record this in the following proposition.

\begin{proposition}
\label{prop:irrational-covers-in-families}
Let $C$ be a smooth connected curve with $0\in C$, and let $f\colon (\mathcal X,\mathcal B)\to C$ be a projective flat family such that, for every $t\in C$,
$\mathcal X_t$ is a klt Fano threefold and
$(\mathcal X_t,\mathcal B_t)$ is a log Calabi--Yau pair with
$\mathcal B_t$ reduced.

Suppose there is a finite surjective $C$-morphism $\pi\colon \mathcal Y\to\mathcal X$
such that $\mathcal Y\to C$ is flat with normal irreducible fibers and $\pi_t\colon \mathcal Y_t\to\mathcal X_t$ is \'etale over $\mathcal X_t\setminus\mathcal B_t$ for every $t\in C$.
If $\mathcal Y_0$ is irrational, then
$(\mathcal X_t,\mathcal B_t)$ is not of cluster type for very general
$t\in C$.
\end{proposition}

\begin{proof}
By Lemma~\ref{lem:clustertypeobstruction}.(1), the fibers $\mathcal Y_t$ are
rationally connected. The locus of
$t\in C$ for which $\mathcal Y_t$ is rational is a countable union of closed
subsets of $C$ by \cite[Theorem~1.2]{dFF13}. Since $\mathcal{Y}_0$ is irrational, all of these subsets are
proper. Hence $\mathcal{Y}_t$ is irrational for very general $t\in C$. Therefore
Lemma~\ref{lem:clustertypeobstruction}.(2) gives that
$(X_t,B_t)$ is not of cluster type.
\end{proof}

\section{Normal quartic surfaces with a cusp singularity}
\label{sec:normal-quartic}
In this section,  we classify cluster type pairs $(\pp^3,S)$ where
$S$ is a normal quartic surface.
First, we recall Wall's parameterization of normal log canonical quartics with a $T_{2,q,r}$-singularity.

 \begin{proposition}[cf.~{\cite[Proposition~3.1 and \S5.2--5.3]{Duc24}}]
\label{prop:ducat-normal-forms}
Let $S\subset\mathbb P^3$ be a log canonical quartic surface with a
$T_{2,q,r}$-singularity at $p:=[1:0:0:0]$. Then, after a projective change of coordinates fixing $p$, the equation defining $S$ has one of the following three forms.
\begin{align}
\quarticfamilytarget{A2}\textup{\quarticfamilylink{A2}{(A.2)}}\qquad
F:={}& a_0w^2x^2+b_2\,wx+c_4,
\label{eq:A2-normal}\\[1mm]
\quarticfamilytarget{A3}\textup{\quarticfamilylink{A3}{(A.3)}}\qquad
F:={}&a_0w^2z^2+b_0wy^3+c_1wxz+d_2wz
+e_0x^3z+f_2x^2+g_3x+h_4,
\label{eq:A3-normal}\\[1mm]
\quarticfamilytarget{A4}\textup{\quarticfamilylink{A4}{(A.4)}}\qquad
F:={}&Q^2+wy^3+(2a_1x+b_2)Q
+c_2x^2+d_3x+e_4,
\quad\text{with}\quad
Q:=wz+x^2.
\label{eq:A4-normal}
\end{align}
In case \textup{\quarticfamilylink{A2}{(A.2)}}, $a_0,b_2$ and $c_4$ are homogeneous polynomials in $x,y,z$ of degree indicated by the subscript. In cases of \textup{\quarticfamilylink{A3}{(A.3)}} and \textup{\quarticfamilylink{A4}{(A.4)}} the
coefficients are homogeneous polynomials in $y,z$ of degree indicated
by the subscript.
\end{proposition}

We call ~\eqref{eq:A2-normal}--\eqref{eq:A4-normal} the \emph{Wall normal forms}. In case \quarticfamilylink{A2}{(A.2)} the coefficients $a_0,b_2$ and $c_4$ involve $1,6,$ and $15$ parameters, respectively. Thus equation~\eqref{eq:A2-normal} is determined by $22$ complex coefficients. We denote the resulting affine space parametrizing these coefficients by $\mathcal W^{\quarticfamilylink{A2}{(A.2)}}$. We define $\mathcal W^{\quarticfamilylink{A3}{(A.3)}}$ and $\mathcal W^{\quarticfamilylink{A4}{(A.4)}}$ similarly. They have complex dimensions $20$ and $17$, respectively.

\begin{definition}
\label{def:wall-parameter-space}
For $i\in\{\quarticfamilylink{A2}{2},\quarticfamilylink{A3}{3},\quarticfamilylink{A4}{4}\}$ and $t\in\mathcal W^{\quarticfamilylink{A2}{(A.i)}}$,
let $S_t:=V(F_t)\subset\mathbb P^3$.  The \emph{Wall coefficient space of type \textup{\quarticfamilylink{A2}{(A.$i$)}}} is
\[
\mathcal N^{\quarticfamilylink{A2}{(A.i)}}
:=
\left\{
t\in\mathcal W^{\quarticfamilylink{A2}{(A.i)}}
\ \middle|\
\begin{array}{l}
S_t \text{ is a normal log canonical quartic surface},\\
{
\operatorname{mult}_q(S_t)\le2
\text{ for every }q\in S_t
}
\end{array}
\right\}.
\]
\end{definition}
On the locus of nonzero quartic forms, the conditions that
$S_t$ be normal and log canonical are open
(see~\cite[Th\'eor\`eme~12.2.4]{EGA4}
and~\cite[Example~9.5.41]{Laz04b}). The locus of multiplicity at least three is closed (see~\cite[Exercise~4.2]{Har92}). Thus $\mathcal N^{\quarticfamilylink{A2}{(A.i)}}$ is open in $\mathcal W^{\quarticfamilylink{A2}{(A.i)}}$
and hence smooth and irreducible.

\begin{proposition}[cf.~{\cite[Proposition~3.1]{Duc24}}]
\label{uniquecuspsurfaces}
Let $S\subset \mathbb P^3$ be a normal strictly log canonical quartic
surface of coregularity $0$. Then $S$ is rational and has at most one cusp singularity. Furthermore, all
its other singularities are Du Val.
\end{proposition}

Since $S$ is strictly log canonical, Proposition~\ref{uniquecuspsurfaces}
implies that its lc center is unique. In each Wall normal form, this lc center is
fixed at $p=[1:0:0:0].$ Thus the corresponding Wall family carries a natural section obtained by
marking the lc center. This motivates the following definition.

\begin{definition}
\em \label{def:marked-wall-family}
Fix $i\in\{\quarticfamilylink{A2}{2},\quarticfamilylink{A3}{3},\quarticfamilylink{A4}{4}\}$. For $t\in\mathcal N^{\quarticfamilylink{A2}{(A.i)}}$
we let $F_t$ denote the corresponding Wall normal form
and set $S_t:=\{F_t=0\}\subset\mathbb P^3$.
The \emph{tautological Wall family of type
\textup{\quarticfamilylink{A2}{(A.$i$)}}} is
\[
 \mathcal S^{\quarticfamilylink{A2}{(A.i)}}
 :=\{(x,t)\in\mathbb P^3\times\mathcal N^{\quarticfamilylink{A2}{(A.i)}}
      \mid F_t(x)=0\}.
\]
\begin{enumerate}
\item Since $p=[1:0:0:0]$ lies on every fiber,
the tautological family has a natural section
$\sigma^{\quarticfamilylink{A2}{(A.i)}}\colon
\mathcal N^{\quarticfamilylink{A2}{(A.i)}}\rightarrow\mathcal S^{\quarticfamilylink{A2}{(A.i)}}$
with $t\rightarrow(p,t)$.
We call $\sigma^{\quarticfamilylink{A2}{(A.i)}}$ the \emph{marked section}.

\item
For a morphism $f\colon B\rightarrow\mathcal N^{\quarticfamilylink{A2}{(A.i)}}$
we set
\[
 \mathcal S_B
 :=\mathcal S^{\quarticfamilylink{A2}{(A.i)}}
   \times_{\mathcal N^{\quarticfamilylink{A2}{(A.i)}}}B
 \subset \mathbb P^3_B
\]
and let $\sigma_B\colon B\rightarrow\mathcal S_B$
be the section induced by $\sigma^{\quarticfamilylink{A2}{(A.i)}}$.
We call $(\mathcal S_B\to B,\sigma_B)$
the \emph{marked Wall family of type
\textup{\quarticfamilylink{A2}{(A.$i$)}} over $B$}.

\item
If for every $b\in B$ the marked point
$\sigma_B(b)\in(\mathcal S_B)_b$
is a singularity of type $T_{2,q,r}$, then we call
$(\mathcal S_B\to B,\sigma_B)$
an \emph{equisingular $T_{2,q,r}$-family of type
\textup{\quarticfamilylink{A2}{(A.$i$)}}}.

For a fixed marked Wall family over $B$ and fixed $q,r$,
we let $\mathcal T^{\quarticfamilylink{A2}{(A.i)}}_{2,q,r}\subset B$
denote the locus of points $b$ for which
$((\mathcal S_B)_b,\sigma_B(b))$ has type $T_{2,q,r}$. If it has several irreducible components,
we denote them by
\[
 \mathcal T^{\quarticfamilylink{A2}{(A.i)},I}_{2,q,r},
 \ldots,
 \mathcal T^{\quarticfamilylink{A2}{(A.i)},m}_{2,q,r}.
\]
\end{enumerate}
\end{definition}

We will occasionally identify a parameter
$t\in\mathcal N^{\quarticfamilylink{A2}{(A.i)}}$ with the corresponding
quartic surface $S_t$.
As an abuse of notation, we may also write
$S\in\mathcal N^{\quarticfamilylink{A2}{(A.i)}}$ or
$S\in\mathcal T_{2,q,r}^{\quarticfamilylink{A2}{(A.i)}}$.
Next, for a stratum 
$\mathcal T^{\quarticfamilylink{A2}{(A.i)}}_{2,q,r}\subset B$,
we establish bounds on $q$ and $r$.
Afterwards, we show that
an equisingular Wall family over a smooth irreducible
base admits a simultaneous $\mathbb Q$-factorial dlt
modification \'etale locally over the base.

\begin{proposition}\label{Theorem:q+rleq18}
Let $S\subset \mathbb P^3$ be a normal log canonical quartic surface with a unique cusp singularity of type $T_{2,q,r}$. Then $q+r \leq 18$. Furthermore, if $q+r=18$, then the surface $S$ has no Du Val singularities.
\end{proposition}

\begin{proof}
Since $S$ has a cusp, $(\mathbb P^3,S)$ has coregularity
zero. By Proposition~\ref{uniquecuspsurfaces}, $S$ is rational and all its
other singularities are Du Val. Let $\widetilde{S}\to S$ be the minimal resolution. Let $n$ be the number of irreducible components of the exceptional divisors over the cusp. Let $k$ be the total number of exceptional curves over the remaining Du Val singularities.

By the classification of normal quartics with an isolated strictly log canonical singularity, we have 
\[
\rho(\widetilde{S})=
\begin{cases}
12, &\quad\text{if}\quad q\geq4,\\
11, &\quad\text{if}\quad q=3.
\end{cases}
\] 
(see \cite[Proposition~3.1 and Table~2]{Duc24}). The exceptional cycle of the $T_{2,q,r}$ singularity has
\[
n=
\begin{cases}
q+r-8, &\quad\text{if}\quad q\geq4,\\
r-6, &\quad\text{if}\quad q=3
\end{cases}
\]
irreducible components (see \cite[Tables~1 and 2]{Lauf77}). Since the pullback of a hyperplane section and the $n+k$ exceptional curves are linearly independent in $\operatorname{NS}(\widetilde{S})$, we have $1+n+k\le\rho(\widetilde{S}).$

In either case, this yields $q+r+k\le 19$.
Let us exclude the case $q+r=19$.
Let $d$ be the absolute value of the determinant of the
intersection matrix of the exceptional curves over the cusp.
Then $d=(q-2)(r-2)-4$ (cf.~\cite[Table~1]{FP84}). Since $\widetilde S$ is a smooth rational surface,
we have
$\operatorname{NS}(\widetilde S)
\cong H^2(\widetilde S,\mathbb Z)$,
whose intersection pairing is unimodular
by Poincar\'e duality.

Suppose that $q+r=19$. Then $k=0$ and
$n+1=\rho(\widetilde S)$. Since $H^2=4$
and $H\cdot E_i=0$, the classes
$H,E_1,\ldots,E_n$ span a full-rank
sublattice $L\subset\operatorname{NS}(\widetilde S)$
with $|\det L|=4d.$
As $\widetilde S$ is rational, its
Neron--Severi lattice is unimodular. Hence $4d=[\operatorname{NS}(\widetilde S):L]^2,$ so $d$ must be a perfect square.
This is impossible for $3\le q\le r$
and $q+r=19$. Thus $q+r\le18$. Finally, if $q+r=18$ and $k>0$, then
$k=1$, so the additional singularity is
of type $A_1$. Its exceptional curve
has self-intersection $-2$ and is
orthogonal to $H,E_1,\ldots,E_n$.
The resulting full-rank sublattice
therefore has discriminant $8d$,
which must be a perfect square.
However, this is impossible when
$3\le q\le r$ and $q+r=18$.
Thus $k=0$.
\end{proof}

\begin{notation}
\label{not:weighted-blowup}
\em Let $(S;p)\cong T_{2,q,r}$ with $q\le r$, and let
$F\in\mathbb C[w,x,y,z]$ be a defining equation of $S$.
We use \textup{\quarticfamilylink{A2}{(A.2)}} for $q\ge4$ and \textup{\quarticfamilylink{A3}{(A.3)}}
or \textup{\quarticfamilylink{A4}{(A.4)}} for $q=3$. If $b_0=0$ in the
\textup{\quarticfamilylink{A3}{(A.3)}} normal form, exchanging $x$ and $z$
gives the \textup{\quarticfamilylink{A2}{(A.2)}} form.

On $D_{\mathbb P^3}(w)$, we use the following local
coordinates and weights. The weighted order of $F$
is computed after setting $w=1$.
\[
\begin{aligned}
\quarticfamilylink{A2}{(\mathrm{A}.2)}&\quad
(x,y,z), && \operatorname{wt}(x,y,z)=(2,1,1),
&& \operatorname{ord}_{\mathrm{wt}}(F)=4,\\
\quarticfamilylink{A3}{(\mathrm{A}.3)}&\quad
(x,y,z), && \operatorname{wt}(x,y,z)=(1,2,3),
&& \operatorname{ord}_{\mathrm{wt}}(F)=6,\\
\quarticfamilylink{A4}{(\mathrm{A}.4)}&\quad
(x,y,Q),\quad Q=wz+x^2,
&& \operatorname{wt}(x,y,Q)=(1,2,3),
&& \operatorname{ord}_{\mathrm{wt}}(F)=6.
\end{aligned}
\]
Let $\pi\colon X\to\mathbb P^3$ be the corresponding
weighted blow-up at $p$, with exceptional divisor $E_1$.
Let $\widetilde S$ be the strict transform of $S$ and set
$C:=\widetilde S\cap E_1$. We write $D_X(x)$ for the affine
chart corresponding to $x$ and $\widetilde F|_{x=1}$ for
a defining equation of $\widetilde S$ on that chart.
We use the same notation for the other local coordinates
in the corresponding row of the table.
\end{notation}
 
\begin{proposition}\label{prop:constancy-of-m}
Let $S$ be a normal quartic surface
in one of the Wall families 
\textup{\quarticfamilylink{A2}{(A.2)}}-\textup{\quarticfamilylink{A4}{(A.4)}} with a $T_{2,q,r}$ singularity.
Let $\pi\colon X\to\pp^3$ be the weighted blow-up described
in Notation~\ref{not:weighted-blowup}.
Then $C$ is nodal, and every singular point of $\widetilde S$
along $E_1$ is a node of $C$. At each node of $C$, the
surface $\widetilde S$ has type $A_m$. The value
$m$ only depends on the analytic type of the cusp.
\end{proposition}

\begin{proof} 
Throughout this proof, we use the local analytic equations given by Proposition~\ref{prop:slc_local_models}.
Suppose first that $S$ belongs to \textup{\quarticfamilylink{A2}{(A.2)}}.
On the chart $D_X(z)$, write $x=s^2x_1$, $y=sy_1$, and $z=s$.
The strict transform is given by
\begin{equation}\label{eq.A.2-blow-up}
\{x_1^2+x_1y_1+s^{q-4}y_1^q+s^{r-4}=0\}
\subset\mathbb A^3_{x_1,y_1,s}.
\end{equation}
Setting $s=0$ gives a node of $C$ at the origin, where
$\widetilde S$ has type $A_{r-5}$. If $q\ge5$, the same
computation on $D_X(y)$ gives type $A_{q-5}$ at
$[0:1:0]\in E_1$.

For \textup{\quarticfamilylink{A3}{(A.3)}}, on the chart $D_X(x)$ write $x=s$, $y=s^2y_1$,
and $z=s^3z_1$. The strict transform is given by
\[
\{z_1^2+y_1z_1+y_1^3+s^{r-6}=0\}
\subset\mathbb A^3_{y_1,z_1,s}.
\]
At the origin, $C$ has a node and $\widetilde S$ has type
$A_{r-7}$. The same computation applies to \textup{\quarticfamilylink{A4}{(A.4)}}
in the coordinates $(x,y,Q)$. Consequently,
\begin{equation}\label{computationofDuvalsingularities}
m=
\begin{cases}
r-5, & \textup{\quarticfamilylink{A2}{(A.2)}},\ q=4,\\
q-5 \text{ and } r-5, & \textup{\quarticfamilylink{A2}{(A.2)}},\ q\ge5,\\
r-7, & \textup{\quarticfamilylink{A3}{(A.3)}},\ \textup{\quarticfamilylink{A4}{(A.4)}}.
\end{cases}
\end{equation}

By computing the local coordinate along the open sets $D_X(x)$ and $D_X(z)$, we can see $C$ is smooth away node. Near each node, there is an analytic isomorphism of
$(X,\widetilde S+E_1)$ with $B_m$
\end{proof}

\begin{proposition}
\label{prop:simultaneous-dlt}
Let $B$ be a smooth and irreducible variety. 
Let $\mathcal{S}\rightarrow B$ be an equisingular family of type \textup{\quarticfamilylink{A2}{(A.i)}} with $i\in \{\quarticfamilylink{A2}{2},\quarticfamilylink{A3}{3},\quarticfamilylink{A4}{4}\}$.
Assume that $\mathcal{S}$ has an $T_{2,q,r}$ singularity on each fiber. 
Then, \'etale locally over $B$, the family 
$(\pp^3\times B,\mathcal{S})\rightarrow B$ 
admits a simultaneous $\qq$-factorial dlt modification.
\end{proposition}

\begin{proof}
Our family of surfaces has a marked section on which every fiber has a cusp singularity.
Away from this marked section
every fiber $S_b$ has only Du Val singularities.
Threfore, the pairs $(\pp^3,S_b)$ are plt for each $b\in B$ by inversion of adjunction. It therefore suffices to construct the modification near the log canonical center $p\times B$.
Perform the weighted blow-up described in Notation~\ref{not:weighted-blowup} along $p\times B$.
After doing so, \'etale locally near any point of the blow-up, we are in the situation of Lemma~\ref{lem:relative-uv-sm-dlt}. Therefore, such lemma implies the existence of the simultaneous $\qq$-factorial dlt modification \'etale locally over $B$.
\end{proof}

\subsection{Quartics of type \quarticfamilylink{A2}{(\textup{A.2})}}

In this subsection, we study the cluster type property of normal quartic surfaces
$S$ in the \textup{\quarticfamilylink{A2}{(A.2)}} family using the double cover
associated with projection from the cusp.
We first construct a crepant birational model $(T,C)$
of $(S,0)$. We then prove that $(\pp^3,S)$ is of
cluster type if $T$ has an $A_m$ singularity,
with $m\ge2$, at a marked node of $C$.

\begin{notation}\label{not:T}
{\em 
Let $S\subset\pp^3$ be a normal quartic surface in the
\textup{\quarticfamilylink{A2}{(A.2)}} family with a cusp $p$ of type $T_{2,q,r}$,
where $4\le q\le r$ and $(q,r)\ne(4,4)$.
Choose coordinates $[w:x:y:z]$ such that $p=[1:0:0:0]$
and write
\begin{equation}\label{eq:A2-square-form}
 F=x^2w^2+xwb_2+c_4,\qquad
 4F=(2wx+b_2)^2-\Delta,\qquad\text{and}\qquad
 \Delta=b_2^2-4c_4,
\end{equation}
where $b_2,c_4\in\mathbb C[x,y,z]$ are homogeneous of
degrees two and four, respectively.
Define $T=\{v^2=\Delta(x,y,z)\}\subset \pp(1,1,1,2)_{x,y,z,v}$. 
Consider the induced double cover $\pi\colon T 
 \rightarrow\pp^2_{x,y,z}.$ Let $D$ and $L$ be the vanishing loci of $\Delta$ and $x$ on $\pp^2$, respectively.
The substitution $v=2wx+b_2$ identifies
$S\cap\{x\ne0\}$ with $T\setminus\pi^{-1}(L)$.
We also set
\[
 C=\pi^{-1}(L),\qquad
 \Theta=V_T(v-b_2),\qquad\text{and}\qquad
 Z=C\cap\Theta=V_T(x,v-b_2).
\]
The birational map
$(S,0)\dashrightarrow(T,C)$ is crepant.
Let $\widetilde S$ be the strict transform of $S$
under the weighted blow-up of $\pp^3$ at marked cusp $p$
with weights $(2,1,1)$.
The induced map $\widetilde S\to T$ is a morphism. After changing coordinates, we may write
\begin{equation}\label{eq:A2-normalized-restriction}
 \Delta|_L=y^2z^2-\delta_{q,4}y^4,
\end{equation}
where $\delta_{q,4}$ is the Kronecker delta. We call the nodes of $C$ the marked nodes.
These are
$p_z=[0:0:1:0]$ and, when $q\ge5$,
$p_y=[0:1:0:0]$.
The indices of $p_y$ and $p_z$ are $q$ and $r$,
respectively.
For each marked node $p_i$, write $\bar p_i=\pi(p_i)$.}
\end{notation}

First, we will prove that whenever $T$ has $A_m$-singularities, with $m\geq 2$, along a marked node, the pair $(\pp^3,S)$ is of cluster type.

\begin{proposition}\label{thm:worse-nodal-normal-form}
We use Notation~\ref{not:T}.
Let $S\subset \pp^3$ be a normal quartic surface in the \textup{\quarticfamilylink{A2}{(A.2)}} family.
If $T$ has an $A_m$-singularity, with $m\ge2$ at a
marked node of $C$, then $(\pp^3,S)$ is of cluster type.
\end{proposition}

\begin{proof}
Using
\eqref{eq:A2-normalized-restriction}, write
\[
 \Delta=z^2(y^2+a_2xy+a_1x^2+a_0xz)
        -\delta_{4,q} y^4+xR,
 \qquad\text{with}\qquad \deg_zR\le1.
\]
Since $T$ is singular at $p_z$, we have $a_0=0$.
The quadratic form $y^2+a_2xy+a_1x^2$ must be a square.
Otherwise the quadratic part of $v^2-\Delta(x,y,1)$
would be nondegenerate, giving an $A_1$-singularity. Set $y':=y+a_2x/2$ and relabel $y'$ as $y$.
We then obtain
\[
 \Delta(1,y,z)=y^2z^2+A(y)z+B(y)
\]
for some $A,B\in\mathbb C[y]$.
Set $t=2w+b_2(1,y,z)-yz$ and $a=2yt-A(y)$.
Then equation~\eqref{eq:A2-square-form} becomes
\[
 4F|_{x=1}=az+t^2-B(y).
\]
By Example~\ref{ex:shear-instances}, there is an open
embedding $\mathbb G_m^3\hookrightarrow\pp^3\setminus S$.
Therefore, Lemma~\ref{clustertypeandtorusembedding} proves the claim.
\end{proof}

We next relate the cusp indices to the singularities of $T$ and the intersection $Z$.

\begin{lemma}\label{lem:A2-singularity-origins}
We use Notation~\ref{not:T}.
Let $S\subset \pp^3$ be a quartic surface in the \textup{\quarticfamilylink{A2}{(A.2)}} family.
Let $p_i\in C\subset T$ be a marked node of index $n$. Then,
\item \begin{equation}\label{eq:A2-local-index}
 n=
 \begin{cases}
  5+\mu(T;p_i),&p_i\notin Z,\\
  4+i_{p_i}(C,\Theta),&p_i\in Z.
 \end{cases}
\end{equation} 
If $p_i\notin Z$, then $(T;p_i)$ has type $A_{n-5}$ where $n=q$ or $n=r$.
If $p_i\in Z$, then $T$ and 
$\Theta$ are smooth at $p_i$.
\end{lemma}

\begin{proof}
We first describe the morphism
$\widetilde S\to T$ locally near a zero dimensional
log canonical place of $S$.
On the open chart $D_T(y)$, put $u=\frac{x}{y},s=\frac{z}{y}$, and $a=\frac{v-b_2}{2y^2}$. On the chart $D_S(y)$, the weighted blow-up is given by $w=1, x=t^2a,
y=t$ and $z=ts$. 
Thus $\widetilde S\to T$ is locally defined by
$u=at$. The same description holds on $D_T(z)$.

If $p_i\notin Z$, then $a(p_i)\ne0$, so
$t=u/a$ and $\widetilde S\to T$ is an
isomorphism near $p_i$.
Thus $(T;p_i)$ is smooth if $n=5$ and has
type $A_{n-5}$ otherwise. Since
$\mu(A_{n-5})=n-5$, we obtain the first
part of the statement. Suppose that $p_i\in Z$, and put
$\bar p_i=\pi(p_i)$.
Since $p_i$ is a node of $C$, we have
$\Delta(\bar p_i)=0$ and
\[
\partial_y\Delta(\bar p_i)
=\partial_z\Delta(\bar p_i)=0.
\]
Moreover, $b_2(\bar p_i)=0$, so
$c_4(\bar p_i)=0$ and
\[
4(\partial_xc_4,\partial_yc_4,\partial_zc_4)
(\bar p_i)
=-(\partial_x\Delta,0,0)(\bar p_i).
\]
If $\partial_x\Delta(\bar p_i)=0$,
then all partial derivatives of $c_4$
vanish at $\bar p_i$.
Since $b_2(\bar p_i)=c_4(\bar p_i)=0$,
the line joining the cusp and $\bar p_i$
would be contained in $\operatorname{Sing}(S)$,
contradicting normality.
Thus $\partial_x\Delta(\bar p_i)\ne0$.
Since $T=\{v^2=\Delta\}$ and
$\Theta\simeq V(c_4)$, both are smooth
at $p_i$.

Put $k=i_{p_i}(C,\Theta)$.
Since $\Theta=V_T(a)$ is smooth at $p_i$,
we consider $(a,s)$ as local analytic
coordinates on $T$ centered at $p_i$. As $C=V_T(u)$, the restriction of $u$
to $a=0$ has order $k$.
Thus we may write
\[
u=ag(a,s)+s^kh(s),\qquad h(0)\ne0.
\]
The local equation of $\widetilde S$ over $p_i$ is therefore $a(t-g(a,s))=s^kh(s),$
which defines an $A_{k-1}$-singularity.
Comparing with the $A_{n-5}$-singularity
of the weighted blow-up of the cusp gives
$k=n-4$.
\end{proof}

Now using Lemma~\ref{lem:A2-singularity-origins}, we give bounds for the cusp indices
in the \textup{\quarticfamilylink{A2}{(A.2)}}. We give bounds in a range which is complementary to such of Proposition~\ref{thm:worse-nodal-normal-form}.

\begin{corollary}\label{eq:A2-index-bound}
We use Notation~\ref{not:T}.
Let $S\subset \pp^3$ be a quartic surface in the \textup{\quarticfamilylink{A2}{(A.2)}} family.
If $T$ has at most $A_1$-singularities at the 
marked nodes, then we have $q+r\leq 12$ or $(q,r)=(6,7)$.
\end{corollary}

\begin{proof}
Suppose that $p_i\in Z$ is a marked node of index $n$.
By Lemma~~\ref{lem:A2-singularity-origins}, we have $k:=i_{p_i}(C,\Theta)=n-4.$
Since $2\le k\le\operatorname{length}Z=4$,
we deduce that $6\le n\le8$.
If $C$ is reducible, and we denote $C=C_1+C_2$, each component $C_i$ maps isomorphically onto $L$ under the projection $\pi$. Since $\Theta$ is defined by a polynomial of
weighted degree two, we have
$\Theta\cdot C_j=2$ for $j\in \{1,2\}$.
As $\Theta$ is smooth, then it is transversal to at least one branch of $C$ at $p_i$. Hence, we get an inequality 
\begin{equation}\label{ineq:A2}
k=i_{p_i}(C_1,\Theta) + i_{p_i}(C_2,\Theta)\leq 3.
\end{equation}
If both nodes lie in $Z$,
then their local lengths are at least two and add up to four. 
Therefore, both lengths are two, showing that $q=r=6$.

Finally, suppose that $T$ is smooth or of type $A_1$
at every marked node. By equality~\eqref{eq:A2-local-index},
a marked node outside $Z$ has index five or six.
If $q=4$, then $r\le8$, so $q+r\le12$.
If $q\ge5$, then $C$ has two irreducible components $C=C_1\cup C_2$ where $C_1=V(x,v-yz)$ and $C_2=V(x,v+yz)$. 
Moreover, in this case, we have $q,r\le7$.
When $r=7$, we have
$\operatorname{length}_{p_z}Z=3$, so
$\operatorname{length}Z=4$ forces $p_y\notin Z$
and hence $q\le6$.
Thus $q+r\le12$ unless $(q,r)=(6,7)$.
\end{proof}

 Let $B$ be the space of polynomials $(b_2,c_4)$ in
equations~\eqref{eq:A2-square-form} for which the corresponding quartic 
belongs to $\mathcal N^{\quarticfamilylink{A2}{(A.2)}}$ and equation~\eqref{eq:A2-normalized-restriction} holds.
We write $\mathcal T_{2,q,r}^{\quarticfamilylink{A2}{(A.2)}}\subset B$ for the
corresponding cusp stratum, with the marked nodes ordered
as above.

\begin{proposition}\label{prop:A2-degenerations}
We use Notation~\ref{not:T}. For $4\le q\le r$ with $(q,r)\ne(4,4)$, consider the
locus in $\mathcal T_{2,q,r}^{\quarticfamilylink{A2}{(A.2)}}$ where $T$ 
has at most $A_1$-singularities at the marked nodes.
Then, the following statements hold.
\begin{enumerate}
\item Each connected component is smooth, irreducible,
and quasi-projective.
\item
For each connected component of an equisingular
locus, the closure in $\mathcal N^{\quarticfamilylink{A2}{(A.2)}}$
of the locus of quartics projectively equivalent
to its members contains one of the following strata.\[
 \mathcal T_{2,4,8}^{\quarticfamilylink{A2}{(A.2)},\mathrm I},\quad
 \mathcal T_{2,4,8}^{\quarticfamilylink{A2}{(A.2)},\mathrm{II}},\quad
 \mathcal T_{2,6,6}^{\quarticfamilylink{A2}{(A.2)},\mathrm I},\quad
 \mathcal T_{2,6,6}^{\quarticfamilylink{A2}{(A.2)},\mathrm{IV}},\quad
 \mathcal T_{2,6,7}^{\quarticfamilylink{A2}{(A.2)},\mathrm I},\quad\text{or}\quad
 \mathcal T_{2,6,7}^{\quarticfamilylink{A2}{(A.2)},\mathrm{II}}.
\]
\end{enumerate}
\end{proposition}

\begin{proof}
Write
\[
 \Delta=y^2z^2-\delta_{q,4} y^4+xH_3,\qquad
 b_2|_L=\alpha y^2+\beta yz+\gamma z^2,\qquad
 h_i=H_3(\bar p_i).
\]
At a marked node $p_i$, the surface $T$ is smooth
if $h_i\ne0$. If $h_i=0$, then $T$ has an $A_1$-singularity at $p_i$ precisely when the quadratic part
of $\Delta$ at $\bar p_i$ is nondegenerate. We have $p_z\in Z$ if and only if $\gamma=0$.
If $\gamma\ne0$, then \eqref{eq:A2-local-index} gives
$r=5$ when $h_z\ne0$ and $r=6$ when $h_z=0$.
If $\gamma=0$, Lemma~\ref{lem:A2-singularity-origins}
gives $h_z\ne0$, and
\begin{equation}\label{eq:A2-contact-orders}
 4c_4(0,y,1)
 = (\beta^2-1)y^2+2\alpha\beta y^3
   +(\alpha^2+\delta_{4,q})y^4.
\end{equation}
The order of vanishing of $c_4|_{x=0}$ at $y=0$ is $r-4$, hence 
\[ r=
 \begin{cases}
  6,&\beta^2\ne1,\\
  7,&\beta=\pm1,\quad\alpha\ne0,\\
  8,&\beta=\pm1,\quad\alpha=0,\quad \delta_{q,4}=1.
 \end{cases}\]
The remaining case has $c_4|_L=0$ and is excluded by
irreducibility. For $q\ge5$, the conditions at $p_y$
are obtained by exchanging $\alpha,\gamma$ and $h_y,h_z$.
Together with Lemma~\ref{lem:A2-singularity-origins} these conditions determine the strata. 

Let $\mu_C(T)$ denote the number of
$A_1$-singularities of $T$ on $C$. We label the strata as follows.
The strata at $(4,5)$ and $(5,5)$ are irreducible.
At $(4,6)$ and $(5,6)$, the irreducible components with 
$\mu_C(T)=0,1$ are labelled $\mathrm I,\mathrm{II}$,
respectively. At $(6,6)$, the labels
$\mathrm I,\mathrm{II},\mathrm{III},\mathrm{IV}$
correspond, respectively, to
\[
 (\mu(T;p_y),\mu(T;p_z))
 =(0,0),(1,0),(0,1),(1,1).
\]
At $(4,7)$, $(4,8)$, $(5,7)$, and $(6,7)$,
the labels $\mathrm I,\mathrm{II}$ correspond to
$\beta=1,-1$.

Each stratum is open in an affine linear space. The sign of $\beta$ is locally constant
when $\beta^2=1$.
Thus, the strata are smooth, irreducible, quasi-projective varieties and are connected components
of the corresponding equisingular loci,
proving \textup{(1)}.

For \textup{(2)}, we record the numerical data of $T$ by the triple
$(q,r,\mu_C(T))$. We can specialize $(5,5,0)$ and
$(5,6,0)$ to $(6,6,0)$ by setting
$\alpha=\gamma=0$ with
$h_yh_z(\beta^2-1)\ne0$.
Similarly, $(4,r,0)$ for $5\le r\le7$
specializes to $(4,8,0)$ by setting
$\alpha=\gamma=0$ and $\beta=\pm1$,
with $h_z\ne0$.
Finally, $(6,6,1)$, $(5,6,1)$, and
$(5,7,0)$ specialize to $(6,7,1)$
by setting $\gamma=h_y=0$ and
$\beta=\pm1$, with $\alpha h_z\ne0$.
Here the triples denote $(q,r,\mu_C(T))$, 
To obtain the remaining specialization,
start with $(6,7,1)$ and replace
$\Delta$ by $\Delta-\eta z^4$, keeping
$b_2$ fixed. For $\eta\ne0$, we have $\Delta|_L=z^2(y^2-\eta z^2),$ and the $A_1$-singularity at $p_y$
remains outside $Z$. The general member belongs
to the $(4,6,1)$ stratum, while the
special member belongs to $(6,7,1)$.

\end{proof}

\begin{proposition}\label{prop:A2-example-fibers}
Each of $\mathcal T_{2,4,8}^{\quarticfamilylink{A2}{(A.2)},\mathrm I},
 \mathcal T_{2,4,8}^{\quarticfamilylink{A2}{(A.2)},\mathrm{II}},
 \mathcal T_{2,6,6}^{\quarticfamilylink{A2}{(A.2)},\mathrm I},
 \mathcal T_{2,6,6}^{\quarticfamilylink{A2}{(A.2)},\mathrm{IV}},
 \mathcal T_{2,6,7}^{\quarticfamilylink{A2}{(A.2)},\mathrm I},
\mathcal T_{2,6,7}^{\quarticfamilylink{A2}{(A.2)},\mathrm{II}}$
contains a normal quartic $S:=\{F=0\}\subset\mathbb P^3$ whose cyclic
quartic cover
$Y:=\{u^4=F\}\subset\mathbb P^4_{u,w,x,y,z}$
is irrational.
\end{proposition}

\begin{proof}
By \eqref{eq:A2-square-form}, the pair $(b_2,\Delta)$ uniquely determines $F$. Thus we specify $F$ by giving $\Delta$ and $b_2$. Note $Y$ is birational with the double solid
\[
X=\{v^2=\Delta(x,y,z)+4u^4\}\subset \mathbb P(1,1,1,1,2).
\]
We will show $X$ is irrational. For $(q,r)\in\{(4,8),(6,6)\}$, take
\[
\Delta=y^2z^2-\delta_{q,4}y^4+x(x^3+y^3+z^3),
\qquad\text{with}\qquad b_2=\pm\delta_{q,4}yz.
\]
Both branch curves $V(\Delta)$ are smooth, and by Lemma~\ref{lem:A2-singularity-origins} and
\eqref{eq.A.2-blow-up}, we have an example of a quartic surface in each of the families
$\mathcal T_{2,4,8}^{\quarticfamilylink{A2}{(A.2)},\mathrm I}$,
$\mathcal T_{2,4,8}^{\quarticfamilylink{A2}{(A.2)},\mathrm{II}}$, and
$\mathcal T_{2,6,6}^{\quarticfamilylink{A2}{(A.2)},\mathrm I}$.
In each case, the surface $V(\Delta+4u^4)$ is also
smooth, so $X$, and hence $Y$, is irrational by
\cite[Theorem~1.2]{CPS19}. For the family $\mathcal T_{2,6,6}^{\quarticfamilylink{A2}{(A.2)},\mathrm{IV}}$, take
\begin{equation}\label{eq:A2-test-66}
\Delta=(y^2-x^2)z^2+x^3z+x^2y^2
\qquad\text{and}\qquad b_2=y^2+z^2.
\end{equation}
The curve $V(\Delta)$ has exactly two nodes,
at $\bar p_y,\bar p_z$, where $b_2=1$.
Thus both marked points are $A_1$ points of $T$
outside $Z$, giving numerical data $(6,6,2)$.
On the chart $D_{\mathbb P^3}(x)$ put $R=y^2-1$ and $H=1+4y^2-4y^4$.
The equation of $X$ becomes
$v^2=Rz^2+z+y^2+4u^4$.
Multiplying by $4R$ gives
\begin{equation}\label{eq:A2-test-conic}
(2Rz+1)^2-R(2v)^2=H-16Ru^4.
\end{equation}
Hence, the conic bundle over $\mathbb{C}(y,u)$ has Brauer class $(R,H-16Ru^4)$.
The polynomials $R$ and $H$ are coprime and square free
of degree two and four, respectively. 
Moreover, the cover $u^4=H/(16R)$ has branching indices $(B_4,B_2)=(6,1)$.
Proposition~\ref{prop:fractional-kummer-irrationality}
therefore applies to show that $X$, and hence $Y$, is irrational.

For the two stratum at $(6,7)$, we consider
\begin{equation}\label{eq:A2-test-67}
\Delta=y^2z^2+xz^3+x^2(y^2-z^2),
\qquad \text{ and } \qquad b_2=y^2\pm yz.
\end{equation}
By Lemma~\ref{lem:A2-singularity-origins},
the two choices of sign give surfaces in
$\mathcal T_{2,6,7}^{\quarticfamilylink{A2}{(A.2)},\mathrm I}$ and
$\mathcal T_{2,6,7}^{\quarticfamilylink{A2}{(A.2)},\mathrm{II}}$,
respectively. Interchanging $x$ and $z$ identifies
the discriminant in \eqref{eq:A2-test-67}
with that in \eqref{eq:A2-test-66}.
Hence the corresponding quartic double solid
$X=\{v^2=\Delta+4u^4\}$ is irrational.
\end{proof}

\subsection{Quartics of type \textup{\quarticfamilylink{A3}{(A.3)}}}

\label{sec:A3}
In this subsection, we study when the pair
$(\mathbb P^3,S)$ is of cluster type for
$S\in\mathcal N^{\quarticfamilylink{A3}{(A.3)}}$.
We show that the cusp strata
$\mathcal T_{2,3,r}^{\quarticfamilylink{A3}{(A.3)}}$ are smooth and
irreducible for $7\le r\le10$ and that
$(\mathbb P^3,S)$ is of cluster type for
$r\ge11$ (Proposition~\ref{prop:A3-low-strata}).
We also provide an example with a
$T_{2,3,10}$ cusp for which $(\mathbb P^3,S)$
is not of cluster type
(Proposition~\ref{prop:irrational-memberA.3}).
The quartic equation appearing in the proof of
Proposition~\ref{prop:A3-low-strata}.(iii)
is related to the type II degenerations of
$K3$ surfaces studied by Jones
\cite[\S 3.7]{Jones26} (see also \cite{Shah81}).

\begin{proposition}
\label{prop:A3-low-strata}
There exists a Wall family $\mathcal S_B\to B$ of type
\textup{\quarticfamilylink{A3}{(A.3)}} such that every quartic of this type with a cusp of
type $T_{2,3,r}$ for some $r\ge7$ is projectively equivalent to a
fiber. Let $\mathcal T_{2,3,r}^{\quarticfamilylink{A3}{(A.3)}}\subset B$ denote the
corresponding cusp stratum. Then the following hold.

\begin{enumerate}
\item[(i)] For $7\le r\le10$, the stratum
$\mathcal T_{2,3,r}^{\quarticfamilylink{A3}{(A.3)}}$ is smooth and irreducible.

\item[(ii)] For $7\le r\le10$ and $s>r$, we have
$\mathcal T_{2,3,s}^{\quarticfamilylink{A3}{(A.3)}}
\subset \overline{\mathcal T_{2,3,r}^{\quarticfamilylink{A3}{(A.3)}}}$
where the closure is taken in $B$.
\item[(iii)] For any $S\in\mathcal T_{2,3,r}^{\quarticfamilylink{A3}{(A.3)}}$
with $r\ge11$, the pair $(\pp^3,S)$
is of cluster type.
\end{enumerate}
\end{proposition}

\begin{proof}
Let $S\in\mathcal N^{\quarticfamilylink{A3}{(A.3)}}$
have a cusp of type $T_{2,3,r}$ with $r\geq 7$.
We perform a weighted blow-up of $\pp^3$ as in
Notation~\ref{not:weighted-blowup}.
On the chart $D_{\mathbb P(1,2,3)}(x)$, the curve $C$ is defined by
\begin{equation}\label{A.3equationofnode}
\{[1:y:z]\mid az^2+by^3+cyz+dz+e_0y^2=0\}
\subset\pp(1,2,3)
\quad\text{with}\quad abd\ne0.
\end{equation}
Indeed, $d=0$ would make $S$ singular along
$V(y,z)$, contrary to normality.

By Proposition~\ref{prop:constancy-of-m},
$C$ has an ordinary node.
Along $x=0$, the curve is smooth since
$az^2+by^3=0$ with $yz\ne0$. $D_{E_1}(x)$, its restrictions to $y=0$ and $z=0$
are $z(az+d)$ and $y^2(by+e_0)$, respectively.
All their zeros are smooth points of $C$ as
those away from the origin are simple,
and the origin is smooth since $d\ne0$. Thus the node lies in $D_{E_1}(yz)$. Rescaling $y,z$ and
the defining equation, we may assume that $b=1$ and that
the node is $q=[1:-1:-1/2]$.
Equation~\eqref{A.3equationofnode} at $q$ is singular, thus we have
$c=2-a$, $d=2$, and $e_0=1+\frac a4$.
Consequently, on the chart $D_{E_1}(x)$, the curve $C$ is defined by
\begin{equation}\label{eq:A.3-leading-term}
\left\{[1:y:z]\mid
az^2+y^3+(2-a)yz+2z+
\left(1+\frac a4\right)y^2=0 \right\}
\subset\mathbb P(1,2,3)
\quad\text{where}\quad a\notin\{0,-1\}.
\end{equation}
The discriminant of the quadratic part at $q$
is $4(a+1)$, which is nonzero since $q$ is an ordinary node.

We turn to perform a change of variables to describe $\widetilde{F}|_{x=1}$, by doing so, we will achieve statement (i).
We give $w,x,y,z$ weights $0,1,2,3$, respectively.
We conside the change of variables
\[
w\mapsto w+\alpha x+\xi y+\theta z,\qquad
x\mapsto x+\beta y+\zeta z,\qquad
y\mapsto y+\gamma z. 
\]
Note that we have $6$ free variables. 
Using these free variables we may eliminate the monomials 
$xy^3$, $x^2yz$, $wy^2z$,
$wyz^2$, and $wz^3$, and make the coefficient
of $x^2z^2$ four times that of $y^4$. By doing so, on the chart$D_X(x)$, the strict transform is given by
\begin{equation}\label{equation:strict-transformation-A_3}
\begin{aligned}
\widetilde F|_{x=1}
=&az^2+y^3+(2-a)yz+2z
+\left(1+\frac a4\right)y^2+esz^2\\
&+s^2\left(gy^2z+f(y^2+2z)^2\right)
+s^3(hyz^2+iy^3z)\\
& + s^4(jz^3+ky^2z^2)
+ \ell s^5yz^3+ms^6z^4.
\end{aligned}
\end{equation}

The monomials $\mathfrak{m}$ in the previous equation have homogeneous degree $6+\deg_{\mathfrak{m}}(s)$.
Let $B$ be the open locus of the coefficient tuples $(a,e,g,f,h,i,j,k,\ell,m)$ variables  with $a\notin\{0,-1\}$. Restricting the tautological
Wall family to $B$ gives the required family
$\mathcal S_B\to B$.

We study the equation
\eqref{equation:strict-transformation-A_3} locally on $q$ to evaluate the cusp index
$r$. By equality~\eqref{computationofDuvalsingularities}, the cusp index is determined by the first nonzero entry of
$(e,g,h-2i,j-2k)$, giving $r=7,8,9,10$, respectively,
and $r\ge11$ if all four vanish.
Thus, the strata for $7\leq r\leq 10$ are open subsets of affine linear spaces.
This proves \textup{(i)}.

For the statement \textup{(ii)}, take a point in the stratum of cusp index
$r'>r$. For $r=7,8,9$ or $10$ vary only the coefficient
$e,g,h,j$, respectively. This proves \textup{(ii)}.

For the statement \textup{(iii)}, the assumption $r\ge11$
gives $e=g=0$, $i=\frac h2$, and $j=2k$.
We work on the chart $D_{\pp^3}(z)$. Then we can write
\begin{equation}\label{equationchartA.3}
\begin{aligned}
F|_{z=1}=\;&
aw^2+w(y^3+(2-a)xy)+2x^3+
\left(1+\frac a4\right)x^2y^2\\
&+hxy+\frac h2y^3+f(y^2+2x)^2
+2kx+ky^2+\ell y+m.
\end{aligned}
\end{equation}
Put $u=y^2+2x$.
On the open set $u\ne0$ in this affine chart, set
\[
s:=w-\frac12xy
+\frac{1-\sqrt{a+1}}{2a}uy
+\frac{1}{\sqrt{a+1}}
\left(\frac h2+\frac{\ell}{u}\right).
\]
Then
\begin{equation}\label{toruschartA.3}
F|_{z=1}
=\sqrt{a+1}\,usy
+a\left(
s-\frac{1}{\sqrt{a+1}}
\left(\frac h2+\frac{\ell}{u}\right)
\right)^2
+\frac14u^3+ku+m+fu^2.
\end{equation}
Let us restrict the coordinate functions $(u,s,y)$ to an open set $D_{\mathbb P^3}(zu)$. Since \[ a\left(
s-\frac{1}{\sqrt{a+1}}
\left(\frac h2+\frac{\ell}{u}\right)
\right)^2
+\frac14u^3+ku+m+fu^2\in \mathbb C[u^{\pm},s]\]
Example~\ref{ex:shear-instances}.\textup{(ii)}
gives an open torus embedding   $\mathbb G_m^3\hookrightarrow \pp^3\setminus S$.
Hence $(\pp^3,S)$ is of cluster type by
Lemma~\ref{clustertypeandtorusembedding}.
\end{proof}

We next construct a quartic surface
$S\in\mathcal T_{2,3,10}^{\quarticfamilylink{A3}{(A.3)}}$ such that
$(\pp^3,S)$ is not of cluster type.
Its defining equation is obtained from the
examples studied by Jones by adding monomials
of the highest possible weighted degrees
$10$ and $11$ (see \cite[\S 3.7]{Jones26}).

\begin{proposition}
\label{prop:irrational-memberA.3}
There exists a polynomial $F\in \mathcal{T}_{2,3,10}^{\quarticfamilylink{A3}{(A.3)}}$ 
such that the quartic hypersurface 
$Y:={u^4=F}\subset \pp^4$ is irrational.
\end{proposition}

\begin{proof}
We consider the following polynomial
\begin{equation} \label{eq:A3-F10}
F = w^2z^2+wy^3+wxyz+2x^3z +\frac54x^2y^2+y^2z^2+yz^3 \in \cc[x,y,z,w].
\end{equation}
We apply the weighted blow-up of
Notation~\ref{not:weighted-blowup}. Then, we can write the strict transform in the form
\eqref{equation:strict-transformation-A_3}
to obtain $e=g=h-2i=0$ and $j-2k=-2$. The Jacobian criterion shows that $V(F)$ is smooth away from
$[1:0:0:0]$, hence normal.
Thus $F\in\mathcal T_{2,3,10}^{\quarticfamilylink{A3}{(A.3)}}$. On the chart $D_{\mathbb P^3}(z)$, we set
$t:=2x+y^2$ and $v:=2w+y^3+xy$.
A direct computation gives
\[
4F|_{z=1}
=
v^2+(4-2t^2)y^2+4y+t^3.
\]
Homogenizing in $v,y$ with a new variable $s$,
we obtain a conic bundle over $\mathbb A^2_{t,u}$
birational to $Y$:
\begin{equation}
\label{eq:A3-conic10}
X:=
\left\{
((t,u),[v:y:s])
\;\middle|\;
v^2+(4-2t^2)y^2+4ys+(t^3-4u^4)s^2=0
\right\}
\subset
\mathbb A^2_{t,u}\times\mathbb P^2_{v,y,s}.
\end{equation}
We compute the discriminant curve as in
Lemma~\ref{lem:conic-brauer-discriminant} and get
\[
\Gamma:=
\left\{
(t,u)\mid
16-4(4-2t^2)(t^3-4u^4)=0
\right\}
\subset\mathbb A^2_{t,u}.
\]
The discriminant curve $\Gamma$ is birational to the
smooth projective model of
$u^4=\frac{t^5-2t^3+2}{4(t^2-2)}$.
We verify the hypotheses of
Proposition~\ref{prop:fractional-kummer-irrationality}.
Let us denote the polynomials
$R:=t^5-2t^3+2$ and $H(t):=t^2-2$.
The generic conic has Brauer class
$\alpha=(H,4Hu^4-R)$, so we apply the proposition
with the polynomials $H,-R$ and the constant $c=4$.
First, (i) is satisfied because
$\deg H=2$ is even.
For (ii), the polynomials $R$ and $H$ are square free
and coprime, with $(\deg R,\deg H)=(5,2)$.
Thus the five zeros of $R$ and the two zeros of $H$
are totally ramified, while
$\ord_\infty\!\left(\frac{R(t)}{4H(t)}\right)=-3$
gives ramification index $4$ at infinity. Thus the cover $\rho\colon C\to\mathbb P^1$ has eight
branch points, all with ramification index $4$, so
$(B_4,B_2)=(8,0)$. In particular, $B_4\ge4$, and $3B_4+2B_2=24$.
Therefore
Proposition~\ref{prop:fractional-kummer-irrationality}
applies, and so $Y$ is irrational.
\end{proof}

\subsection{Quartics of type \textup{\quarticfamilylink{A4}{(A.4)}}}
\label{sec:A4}

In this subsection, we study when the pair
$(\mathbb P^3,S)$ is of cluster type for
$S\in\mathcal N^{\quarticfamilylink{A4}{(A.4)}}$.
We classify the irreducible components of the strata
$\mathcal T_{2,3,r}^{\quarticfamilylink{A4}{(A.4)}}$ for $7\le r\le11$
and prove their smoothness
(Lemma~\ref{lem:A.4-low-strata} and
Proposition~\ref{prop:A4-strata}).
We show that $(\mathbb P^3,S)$ is of cluster type
for $r\ge12$ (Proposition~\ref{prop:A4-strata}(iv)),
and provide non-cluster type examples for
$r=10$ and $r=11$
(Proposition~\ref{prop:A4-irrationalexample}).

The quartics considered in
Lemma~\ref{lem:A.4-low-strata}
are type III analogues of the type II degenerations
of quartic surfaces studied by Jones
\cite[\S 3.8]{Jones26}.

\begin{lemma}\label{lem:A.4-low-strata}
There exist an open subset $B\subset\mathbb A^{10}$ and a
Wall family $\mathcal S_B\to B$ of type \textup{\quarticfamilylink{A4}{(A.4)}}
with the following properties.
\begin{enumerate}
\item[(i)] Up to projective equivalence, every quartic in
$\mathcal N^{\quarticfamilylink{A4}{(A.4)}}$ with a cusp is represented by a point of $B$.
\item[(ii)] For $r=7,8,9$, the stratum
$\mathcal T_{2,3,r}^{\quarticfamilylink{A4}{(A.4)}}$ is smooth, irreducible, and
locally closed in $B$.
\end{enumerate}
\end{lemma}

\begin{proof}
First we aim to prove~\noindent\textup{(i)}.
Recall the Wall normal form \eqref{eq:A4-normal}
\[
F=Q^2+wy^3+(2a_1x+b_2)Q+c_2x^2+d_3x+e_4,
\qquad Q=wz+x^2.
\]
We use the weighted blow-up at $p=[1:0:0:0]$ from
Notation~\ref{not:weighted-blowup}. On the chart $D_{\mathbb P^3}(wx)$ we perform weighted blow up via $x=s$, $y=s^2y_1$, and $Q=s^3q$, so that
$z=s^2(sq-1)$. By Proposition~\ref{prop:constancy-of-m},
the node of $C$ lies on $D_{E_1}(x)$. Let $[1:y_0:q_0]$ be a coordinate of node. Perform a change of coordinates $y_1'=y_1-y_0$ and $q'=q-2q_0$ and relabelling it, we may assume we have the node $[1:0:0]$. 

Write $a_1=a_{10}y+a_{01}z$ and
$c_2=c_{20}y^2+c_{11}yz+c_{02}z^2$. On $D_{E_1}(x)$,
the curve $C$ is defined by
\[
q^2+2(a_{10}y_1-a_{01})q+y_1^3
+c_{20}y_1^2-c_{11}y_1+c_{02}=0.
\]
Since the origin is a node, we have
$a_{01}=c_{11}=c_{02}=0$ and $a_{10}^2-c_{20}\ne0$.
Choose $t$ satisfying $4t^2-4a_{10}t+c_{20}=0$ and set
$x'=x+ty$ and $Q'=wz+(x')^2$ and after rescalling and relabelling, we get
\begin{equation}\label{eq:A.4reduced}
F=Q^2+wy^3+(xy+b_2)Q+c_3x+d_4.
\end{equation}
Write $b_{ab},c_{ab},d_{ab}$ for the coefficients of
$y^az^b$ in $b_2,c_3,d_4$, respectively. We may assume
$c_{12}=0$ by setting $w'=w+c_{12}z$ and relabelling.
The change $x'=x-\eta z$, $w'=w+2\eta x-\eta^2z$ fixes $Q$
and replaces $(b_2,c_3)$ by $(b_2+\eta yz,c_3-2\eta y^3)$,
so we may arrange $c_{21}=\tfrac12 b_{11}-2b_{02}$.

Let $B$ parametrize the equations \eqref{eq:A.4reduced}
in $\mathcal N^{\quarticfamilylink{A4}{(A.4)}}$ satisfying these two conditions.
The free coefficients $(b_2,c_{30},c_{03},d_4)$ identify
$B$ with an open subset of $\mathbb A^{10}$, and then equation
\eqref{eq:A.4reduced} defines the family $\mathcal S_B\to B$.
This proves the statement \textup{(i)}.

Now, we turn to prove statement
\noindent\textup{(ii)}.
On the chart $D_X(x)$, the strict transform of a fiber is defined by
\begin{equation}\label{eq:A.4strict-transform}
q^2+y_1q+y_1^3
+s\bigl(b_2(y_1,sq-1)q+c_3(y_1,sq-1)\bigr)
+s^2d_4(y_1,sq-1)=0.
\end{equation}
To expand the equation~\eqref{eq:A.4strict-transform} at the node,
set $y_1'=y_1+b_{02}s$. Since $c_{12}=0$, the coefficients
of the terms linear in $y_1'$ and $q$ are $O(s^2)$.
By setting $y_1'=q=0$, we have
\[
-c_{03}s+d_{04}s^2
+b_{02}\bigl(d_{13}+b_{02}^2-\tfrac12 b_{02}b_{11}\bigr)s^3
+O(s^4).
\]
By Proposition~\ref{prop:constancy-of-m}, the previous equation is of type $T_{2,3,r}$ precisely when the series has order of vanishing  $r-6$. Thus, the strata for $r=7,8,9$ are the open subsets of $B$, $B\cap V(c_{0,3})$, 
and $B\cap V(c_{03},d_{04})$, respectively. 
These open subsets are defined by the non-vanishing of the corresponding coefficients in the series.
\end{proof}

In the following proposition, we use the notation from Lemma~\ref{lem:A.4-low-strata}.
In particular, $B$ is the base constructed in such proposition.

\begin{proposition}
\label{prop:A4-strata}
Let $\mathcal{T}_{2,3,r}^{\quarticfamilylink{A4}{(A.4)}}\subset B$ 
be the stratum with $T_{2,3,r}$ singularities. Then, the following statements hold.
\begin{enumerate}
\item[(i)] The $r=10$ stratum has two irreducible components $\mathcal T_{2,3,10}^{\quarticfamilylink{A4}{(A.4)},I}$ and $
\mathcal T_{2,3,10}^{\quarticfamilylink{A4}{(A.4)},II}$.
\item[(ii)] Only one of these components contains
$\mathcal T_{2,3,11}^{\quarticfamilylink{A4}{(A.4)}}$ in its closure in $B$.
\item[(iii)] For $r\le11$, every irreducible component of
$\mathcal T_{2,3,r}^{\quarticfamilylink{A4}{(A.4)}}$ is smooth.
\item[(iv)] If $S\in\mathcal N^{\quarticfamilylink{A4}{(A.4)}}$ has a cusp of type $T_{2,3,r}$ with $r\geq12$, then
$(\mathbb P^3,S)$ is of cluster type.
\end{enumerate}
\end{proposition}

\begin{proof}
We keep the notation from the proof of Lemma~\ref{lem:A.4-low-strata}. Assume $r\geq9$, so that $c_{03}=d_{04}=0$.
Since $c_{12}=0$, we have $c_3=c_{30}y^3+c_{21}y^2z$.
Define $\rho_3,\rho_4$ and $\rho_5$ by
\begin{equation}\label{equationrho}
\rho_3:=d_{13}+b_{02}^2-\frac12b_{02}b_{11},\quad
\rho_4
:=
d_{22}-b_{02}c_{30}-\frac14b_{11}^2,
\qquad\text{and}\qquad
\rho_5
:=
d_{31}-\frac12b_{20}b_{11}.
\end{equation}
We determine the singularity of the strict transform $\widetilde S$
at the node. We solve the equations $(\tilde F|_{x=1})_{y_1}=(\tilde F|_{x=1})_q=0$ and get
\begin{equation}\label{A.4-term-expansion}
y_1(s)
=
-b_{02}s-2\rho_3s^2+O(s^3)
\quad\text{and}\quad
q(s)
=
\left(\rho_3-\frac12b_{02}b_{11}\right)s^2+O(s^3).
\end{equation}
Substituting \eqref{A.4-term-expansion} into $\tilde F|_{x=1}$ gives the series
$
b_{02}\rho_3s^3
+
\bigl(\rho_3^2+b_{02}^2\rho_4\bigr)s^4
+
O(s^5)$.
By Proposition~\ref{prop:constancy-of-m}, its order is $r-6$. Thus the $r=10$ stratum is the union of the following two loci in $B$:
\begin{equation}\label{twoirreducibleA.2}
\mathcal T_{2,3,10}^{\quarticfamilylink{A4}{(A.4)},\text{I}}
:=
\{c_{03}=d_{04}=\rho_{3}=0,\;
  b_{02}\neq0,\rho_4\neq 0\}\quad\text{and} \quad
\mathcal T_{2,3,10}^{\quarticfamilylink{A4}{(A.4)},\text{II}}
:=
\{c_{03}=d_{04}=b_{02}=0,\;
  \rho_3\neq0\}.
\end{equation}
Suppose that $b_{02}=\rho_3=0$. Then $b_2=y(b_{20}y+b_{11}z)$ and
$c_3=y^2(c_{30}y+c_{21}z)$. By \eqref{equationrho}, we have $d_{13}=0$, and therefore $d_4
=
y^2(d_{40}y^2+d_{31}yz+d_{22}z^2)$. Writing $b_2=y\bar b_1$,
$c_3=y^2\bar c_1$ and
$d_4=y^2\bar d_2$, equation~\eqref{eq:A.4reduced} becomes
\[
F
=
Q^2+y(x+\bar b_1)Q
+y^2\bigl(wy+x\bar c_1+\bar d_2\bigr).
\]
It follows that $\{Q=y=0\}\subset \operatorname{Sing}(S)$,
so $S$ has non-isolated singularities. Thus $b_{02}$ and $\rho_3$ cannot both vanish on the locus
$c_{03}=d_{04}=0$ in $B$. For $r\geq11$, the coefficients of $s^3$ and $s^4$ vanish, so $b_{02}\ne0$ and $\rho_3=\rho_4=0$.
The closure of the second $r=10$ locus lies in $V(b_{02})$.
Varying $\rho_4$ shows that $\mathcal T_{2,3,10}^{\quarticfamilylink{A4}{(A.4)},\text{I}}$ contains $\mathcal T_{2,3,11}^{\quarticfamilylink{A4}{(A.4)}}$ in its closure. Varying $c_{03}$, $d_{04}$, and $\rho_3$, respectively,
also shows that
$\mathcal T_{2,3,11}^{\quarticfamilylink{A4}{(A.4)}}
\subset\overline{\mathcal T_{2,3,r}^{\quarticfamilylink{A4}{(A.4)}}}$
for $r=7,8,9$.

Suppose that $\rho_3=\rho_4=0$. Solving the same equations, we get
\[
y_1
=
-b_{02}s
+b_{02}^2(2c_{30}-b_{20})s^3
+O(s^4)\quad\text{and}\quad
q
=
-\frac12b_{02}b_{11}s^2
-b_{02}^2c_{30}s^3
+O(s^4).
\]
Substitution into $\tilde F|_{x=1}$ gives the polynomial
$b_{02}^3\rho_5s^5+O(s^6)$. Hence the $r=11$ stratum is determined by the equations 
$c_{03}=d_{04}=\rho_3=\rho_4=0$ and $b_{02}\rho_5\ne0$.
Furthermor, for $r \geq 12$, we also have $\rho_5=0$.
We replace $d_{13},d_{22},$ and $d_{31}$
by $\rho_3,\rho_4,$ and $\rho_5$, respectively.
Then, the two loci in~\eqref{twoirreducibleA.2}
and the $r=11$ stratum are open subsets of affine linear subspaces in the $rho$ variables. Together with Lemma~\ref{lem:A.4-low-strata} applied to the cases $r=7,8,9$, completes the argument of \textup{(i)--(iii)}.

For \textup{(iv)}, by Lemma~\ref{lem:A.4-low-strata}.\textup{(i)} we may assume that $S$ is represented by a point of $B$.
On $B$ we have $c_{21}=\frac{b_{11}}2-2b_{02}$, and by \eqref{equationrho} the equalities $\rho_3=\rho_4=\rho_5=0$ give
$d_{13}=\frac{b_{02}b_{11}}2-b_{02}^2$,
$d_{22}=b_{02}c_{30}+\frac{b_{11}^2}{4}$, and
$d_{31}=\frac{b_{20}b_{11}}2$. On $D_{\mathbb P^3}(y)$, set $p:=x+b_{02}z^2$,
$u:=Q+\frac{b_{11}}{2}z+c_{30}$,
$v:=
w-2b_{02}xz+\frac{b_{11}}{2}-b_{02}^2z^3$ and set $g(u):=(u-c_{30})^2+b_{20}(u-c_{30})
+d_{40}-\frac{b_{11}}{2}$. Then
\begin{equation}\label{A.4torusembedding}
\begin{aligned}
F|_{y=1}
={}&
\left(Q+\frac{b_{11}}2z+c_{30}\right)
       (x+b_{02}z^2)+\left(
w-2b_{02}xz+\frac{b_{11}}2-b_{02}^2z^3
\right)\\
&+\left(Q+\frac{b_{11}}2z\right)^2
+b_{20}\left(Q+\frac{b_{11}}2z\right)
+\left(d_{40}-\frac{b_{11}}2\right)\\
={}&up+v+g(u).
\end{aligned}
\end{equation}

Moreover, we have
\[
\begin{aligned}
u-p^2-c_{30}
&=
\left(wz+x^2+\frac{b_{11}}{2}z+c_{30}\right)
-\left(x+b_{02}z^2\right)^2-c_{30}\\
&=
z\left(
w-2b_{02}xz+\frac{b_{11}}{2}-b_{02}^2z^3
\right)=zv.
\end{aligned}
\]
The coordinate functions $(v,z,p)$ identifies $D_{\mathbb P^3}(yv)$ with $\mathbb G_m\times\mathbb A^2$. The equations are of the form considered in
Example~\ref{ex:shear-instances}.\textup{(ii)}. Hence
$\mathbb P^3\setminus S$ contains an open subset isomorphic to
$\mathbb G_m^3$. Therefore $(\mathbb P^3,S)$ is of cluster type by
Lemma~\ref{clustertypeandtorusembedding}.
\end{proof}

\begin{proposition}
\label{prop:A4-irrationalexample}
There exists $F_0 \in \mathcal T_{2,3,10}^{\quarticfamilylink{A4}{(A.4)},\text{II}}$ and $F_1\in\mathcal{T}_{2,3,11}^{\quarticfamilylink{A4}{(A.4)}}$ such that $Y_{\varepsilon}:=\{u^4=F_{\varepsilon}\}$ is irrational for $\epsilon\in \{0,1\}$. 
\end{proposition}
 
\begin{proof}
For $\varepsilon\in\{0,1\}$, set
\begin{equation}
\label{eq:A4-last-strata}
F_\varepsilon
=
Q^2+(xy-\varepsilon z^2)Q+wy^3
+2\varepsilon xy^2z-\varepsilon y^3z
+(1-2\varepsilon)yz^3,\quad \text{with}\quad Q=wz+x^2.
\end{equation}
Comparing with equations~\eqref{eq:A.4reduced} and \eqref{equationrho}, we obtain $c_{03}=d_{04}=0,$ $
b_{02}=-\varepsilon,$ $
\rho_3=(1-\varepsilon)^2,\rho_4=0$ and $\rho_5=-\varepsilon$. For both values of $\varepsilon$, the Jacobian criterion shows that
$V(F_\varepsilon)$ is smooth away from $[1:0:0:0]$, hence normal. Therefore, depending on whether $\varepsilon$ is $0$ or $1$, we have $F_0\in  \mathcal{T}_{2,3,10}^{\quarticfamilylink{A4}{(A.4)},\text{II}}$  or $F_1\in\mathcal T_{2,3,11}^{\quarticfamilylink{A4}{(A.4)}}$. On the chart $D_{\mathbb P^3}(z)$, put \[F_\varepsilon|_{z=1}
=
q^2+(xy+y^3-\varepsilon)q
-x^2y^3+2\varepsilon xy^2-\varepsilon y^3+(1-2\varepsilon)y\quad\text{with}\quad q:=Q|_{z=1}=w+x^2.
\]
Set $v:=2q+xy+y^3-\varepsilon$.
Homogenizing in $v,x$ with a new variable $s$,
we obtain a conic bundle over $\mathbb A^2_{y,u}$
birational to $Y_\varepsilon$:
\[
\{ v^2-y^2(1+4y)x^2
+2y(\varepsilon (1+4y)-y^3)xs
-\bigl(y^6-4y+\varepsilon(2y^3+8y+1)+4u^4\bigr)s^2=0\}\subset \mathbb P^{2}_{v,x,s}\times\mathbb A_{y,u}^2.
\]
Let $\Gamma$ be the closure in $\mathbb P^1_y\times\mathbb P^1_u$
of the discriminant over $y\ne0$. By Lemma~\ref{lem:conic-brauer-discriminant}, on this open set
the discriminant is given by
\[
(1+4y)u^4+y^{1+2\epsilon}P_\epsilon(y)=0,\qquad\text{where}\qquad
P_\varepsilon(y):=
\begin{cases}
y^6-4y-1,&\varepsilon=0,\\
y^4+4y+1,&\varepsilon=1.
\end{cases}\]
Moreover, completing the square in $x$, the
Brauer class of the conic bundle is represented by
\[
\alpha=\left(1+4y,\;-(1+4y)u^4-y^{1+2\varepsilon}P_\varepsilon(y)\right).
\]

The discriminant curve $\Gamma$ is birational to the smooth projective model of $u^4=
-\frac{y^{1+2\varepsilon}P_\varepsilon(y)}{1+4y}.$ We verify the hypotheses of Proposition~\ref{prop:fractional-kummer-irrationality}.
In the notation of that proposition, we have $\alpha=(R(y),cR(y)u^4+H_\varepsilon(y))$ with $R(y):=1+4y,$ $H_\varepsilon(y):=-y^{1+2\varepsilon}P_\varepsilon(y)$ and $ c=-1.$
First, $\deg R=1$ and $\deg H_\varepsilon=7$, so $4\nmid \bigl(\deg R-\deg H_\varepsilon\bigr)=-6.$ Thus condition (i)
of Proposition~\ref{prop:fractional-kummer-irrationality} holds. Moreover, $P_\varepsilon$ is square free and is coprime to
$y(1+4y)$. 
Hence the zeros of $P_\varepsilon$,
the point $y=0$, and the zero $y=-1/4$
of $R$ are all totally ramified. Indeed, $\ord_{y=0}\!
\left(
\frac{H_\varepsilon(y)}{R(y)}
\right)
=1+2\varepsilon\in\{1,3\},$
so the ramification index at $y=0$ is $4$.
Since $\deg P_\varepsilon=6-2\varepsilon$, this gives $B_4=8-2\varepsilon.$ At infinity, $\operatorname{ord}_\infty(H_\varepsilon/R)=-6$,
so the fiber of $\rho\colon C\to\mathbb P^1$ over infinity
consists of two points, each with ramification index $2$.
Thus infinity contributes one branch point to $B_2$, and
$(B_4,B_2)=(8-2\varepsilon,1)$.
In particular, $B_4\geq4,$ and $
3B_4+2B_2\geq 20$. Proposition~\ref{prop:fractional-kummer-irrationality} therefore applies, and $Y_\varepsilon$ is
irrational for $\varepsilon=0,1$.
\end{proof}

\subsection{Cluster type condition for $(\mathbb P^3,S)$ when $S$ is normal}
\label{sec:normal-classification}

In this subsection, we complete the proof of Theorem~\ref{thm:intro-A2} and Theorem~\ref{them:mainclusternormal} by combining the results of the previous subsections.
\begin{proof}[Proof of Theorems~\ref{thm:intro-A2}~ and \ref{them:mainclusternormal}]
By Propositions~\ref{prop:A2-degenerations},~\ref{prop:A3-low-strata}, and~\ref{prop:A4-strata}, the equisingular strata considered below are
smooth. Moreover, by
Proposition~\ref{prop:simultaneous-dlt}, every irreducible equisingular
family of pairs $(\mathbb P^3,S)$ considered below admits a simultaneous $\mathbb Q$-factorial
dlt modification \'etale locally over the base. Hence, by Corollary~\ref{cor:p3-families}, the cluster type
property is constant on each connected component. We may therefore apply
Proposition~\ref{prop:irrational-covers-in-families}. If a stratum
specializes to a member whose cyclic quartic cover is irrational,
then $(\mathbb P^3,S)$ is not of cluster type for every member $S$
of that stratum.

For type \textup{\quarticfamilylink{A2}{(A.2)}}, Proposition~\ref{prop:A2-example-fibers} provides a member with an irrational cyclic quartic cover in each of the six strata listed in Proposition~\ref{prop:A2-degenerations}.(3).
By Proposition~\ref{prop:A2-degenerations}, every connected component considered there specializes to one of these six strata.
The preceding argument therefore shows that $(\mathbb P^3,S)$ is not of cluster type for every $S$ in any of these components.
The remaining \textup{\quarticfamilylink{A2}{(A.2)}} strata are of cluster type by Proposition~\ref{thm:worse-nodal-normal-form}. In particular, $(\mathbb P^3,S)$ is not of cluster type
if $r\le6$, and, by Corollary~\ref{eq:A2-index-bound}, is of
cluster type if $q+r\ge13$ and $(q,r)\ne(6,7)$.

For type \textup{\quarticfamilylink{A3}{(A.3)}}, Proposition~\ref{prop:irrational-memberA.3}
gives a member of $\mathcal T_{2,3,10}^{\quarticfamilylink{A3}{(A.3)}}$ whose associated cyclic quartic
cover is irrational. By Proposition~\ref{prop:A3-low-strata}, $\mathcal T_{2,3,7}^{\quarticfamilylink{A3}{(A.3)}},
\mathcal T_{2,3,8}^{\quarticfamilylink{A3}{(A.3)}}$ and $
\mathcal T_{2,3,9}^{\quarticfamilylink{A3}{(A.3)}}$
specialize to $\mathcal T_{2,3,10}^{\quarticfamilylink{A3}{(A.3)}}$. Hence
$(\mathbb P^3,S)$ is not of cluster type for $S\in\mathcal T_{2,3,r}^{\quarticfamilylink{A3}{(A.3)}},$ with $ 7\le r\le10.$

For type \textup{\quarticfamilylink{A4}{(A.4)}}, Proposition~\ref{prop:A4-irrationalexample}
gives members of $\mathcal T_{2,3,11}^{\quarticfamilylink{A4}{(A.4)}}$ and
$\mathcal T_{2,3,10}^{\quarticfamilylink{A4}{(A.4)},\text{II}}$ whose associated cyclic quartic covers of $\mathbb P^3$ are
irrational. By Proposition~\ref{prop:A4-strata}, every other
irreducible component with $7\le r\le11$ specializes to
$\mathcal T_{2,3,11}^{\quarticfamilylink{A4}{(A.4)}}$. The remaining component \(T_{2,3,10}^{\quarticfamilylink{A4}{(A.4)},\text{II}}\) contains a member whose associated cyclic quartic cover of $\mathbb P^3$ is irrational by Proposition~\ref{prop:A4-irrationalexample}. Therefore $(\mathbb P^3,S)$ is not of cluster type for $S\in\mathcal T_{2,3,r}^{\quarticfamilylink{A4}{(A.4)}}$ with $ 7\le r\le11.$

It remains to consider the complementary loci.
For the \textup{\quarticfamilylink{A2}{(A.2)}} family,
Proposition~\ref{thm:worse-nodal-normal-form} shows
that $(\mathbb P^3,S)$ is of cluster type whenever $T$ has an $A_m$ singularity with $m\ge2$ over $L$.
By \eqref{eq:A2-index-bound}, such a singularity must
occur if $q+r\ge13$ and $(q,r)\ne(6,7)$.
For the \textup{\quarticfamilylink{A3}{(A.3)}} and \textup{\quarticfamilylink{A4}{(A.4)}} families,
Propositions~\ref{prop:A3-low-strata}.\textup{(iii)}
and~\ref{prop:A4-strata}.\textup{(iv)} give the cluster type condition
for $r\ge11$ and $r\ge12$, respectively.
This proves the theorem.
\end{proof}

\section{Non-normal quartic surfaces with degenerate cusp singularities}
\label{sec:nonormal}
We now classify the possible non-normal loci of slc quartic surfaces of coregularity zero in $\pp^3$.

\begin{proposition}\label{prop:reduciblenonnormallocus}
Let $S \subset \mathbb{P}^3$ be an irreducible non-normal quartic surface
such that the pair $(\mathbb{P}^3, S)$ has coregularity zero. Then the non-normal locus
$\Sigma$ of $S$ is one of the following:
\begin{enumerate} 
\item a line; 
\item a smooth plane conic;
\item two intersecting lines;
\item a twisted cubic;
\item a smooth conic and a line not contained in its plane, meeting at one point;
\item three concurrent non-coplanar lines; or
\item a chain of three lines in which the first and the last lines are not coplanar.
\end{enumerate}

\end{proposition}

\begin{proof}
By \cite[Proposition~3.3]{Duc24}, the reduced non-normal locus
$\Sigma$ has degree at most $3$ and is either one of (1), (2), (4) or a degeneration
of one of them.

Let $C$ be a component of $\Sigma$ contained in a plane $\Pi$. Since $S$ is slc,
$S$ is generically double along $C$. Since $S$ is irreducible, $\Pi\not\subset S$,
so $\Pi\cap S$ is a plane quartic curve containing $C$ with multiplicity at least $2$.
Hence $2\deg C\le 4$, that is, $\deg C\le 2$.

Suppose that $\Sigma$ is reducible with two irreducible components. If $\Sigma$ is a
union of two lines, then the lines either intersect or are skew. If they intersect,
we are in case (3). If they are skew, then the pair $(\mathbb P^3,S)$ has
coregularity one by \cite[Proposition~3.3]{Duc24}, a contradiction.

If $\Sigma$ is the union of a line $L$ and a smooth conic $C\subset\Pi$, then $L$
and $C$ intersect. Indeed, by the degree count above, $\Pi\cap S=2C$ as divisors
on $\Pi$. In particular, we have $\Pi\cap S=C$ set theoretically. The line $L$ is not
contained in $\Pi$, since otherwise $L\subset\Pi\cap S=C$. Hence $L$ meets $\Pi$
in exactly one point $p$, and since $p\in L\subset S$, we have $p\in\Pi\cap S=C$.
This is case (5).

Three coplanar double lines $L_1,L_2,L_3\subset\Pi$
are impossible, since $4=\deg S|_\Pi
\ge \deg(2L_1+2L_2+2L_3)=6.$

Next, suppose that $L_1,L_2,L_3$ are pairwise
skew. They lie in the same ruling of a smooth
quadric surface
$Q\cong\mathbb P^1\times\mathbb P^1$
(see, e.g.,\cite[\S 2.12]{Har92}).
If $Q\not\subset S$, then $S|_Q$ has bidegree
$(4,4)$. Since $S|_Q\ge 2L_1+2L_2+2L_3$
and the right hand side has bidegree $(6,0)$,
there would be an effective divisor of
bidegree $(-2,4)$, a contradiction.
Thus $Q\subset S$, contrary to the
irreducibility of $S$.

Finally, suppose that $L_1$ and $L_2$ meet,
and let $\Pi$ be the plane spanned by them.
Since $S$ has degree four, $S|_\Pi=2L_1+2L_2.$
In particular, $S\cap\Pi=L_1\cup L_2$
set-theoretically. Since $L_3\subset S$
and $L_3\not\subset\Pi$, the point
$L_3\cap\Pi$ lies on $L_1\cup L_2$.
Hence $L_3$ meets at least one of
$L_1$ and $L_2$.
\end{proof}

In cases (4) and (6), the pairs $(\mathbb{P}^3, S)$ are known to be of cluster type~\cite[Theorem~6.1]{EFM26}. Moreover, apart from degenerate cusp singularities, these surfaces have only Du Val singularities.

\begin{proposition}[{\cite[Proposition~3.3]{Duc24}}]
\label{uniquedegeneratecuspsurfaces}
Let $S\subset\mathbb P^3$ be an irreducible
non-normal slc quartic surface of coregularity zero.
Then all singularities of $S$ away from its
reduced non-normal locus are Du Val.
\end{proposition}

\subsection{Quartic surfaces with irreducible non--normal locus}
\label{subsec:irreducible-non-normal}

We now describe normal forms for non-normal slc quartic surfaces whose  non--normal locus is a double curve.
\begin{proposition}[cf.~{\cite[Proposition~3.3 and \S5.4--5.5]{Duc24}}]\label{prop:B1B2}
Let $S\subset\mathbb{P}^3$ be an irreducible non-normal slc quartic surface whose double curve $C$ is either a line or a smooth conic. Then, after a projective change of coordinates, the equation defining $S$ has one of the following two forms.
\begin{alignat}{2}
  &\quarticfamilytarget{B1}\textup{\quarticfamilylink{B1}{(B.1)}}\qquad & F &:= a_2w^2+b_2\,wx+c_2x^2+d_3w+e_3x+f_4, \label{eq:B1-normal}\\
  &\quarticfamilytarget{B2}\textup{\quarticfamilylink{B2}{(B.2)}}\qquad & F &:= Q^2+a_1\,wQ+b_2\,w^2. \label{eq:B2-normal}
\end{alignat}
In case \textup{\quarticfamilylink{B1}{(B.1)}}, $C$ is the line $\{y=z=0\}$, and $a_2,b_2,c_2,d_3,e_3,f_4$ are homogeneous polynomials in $y,z$ of degree indicated by the subscript. In case \textup{\quarticfamilylink{B2}{(B.2)}}, $C$ is the smooth conic $\{w=Q=0\}$, where $Q$ is a nondegenerate quadratic form in $x,y,z$, and $a_1,b_2$ are homogeneous polynomials in $w,x,y,z$ of degree indicated by the subscript.
\end{proposition}

We first record some basic properties of non-normal slc quartic surfaces of types \textup{\quarticfamilylink{B1}{(B.1)}} and \textup{\quarticfamilylink{B2}{(B.2)}}.

\begin{proposition}\label{prop:conductor-normalization}
Let $S\subset \mathbb P^3$ be a non-normal slc quartic surface whose
reduced double curve $C$ is either a line or a smooth conic. Assume that
$S$ has multiplicity two at every point of $C$. Let
\[
\pi_C\colon W:=\operatorname{Bl}_C\mathbb P^3\rightarrow\mathbb P^3
\]
be the blow-up of $C$, with exceptional divisor $E$, and let
$\widetilde S\subset W$ be the strict transform of $S$. Put $\nu:=\pi_C|_{\widetilde S}$ and $
\widetilde C:=E|_{\widetilde S}.$ Then the following statements hold.
\begin{enumerate}
\item
The morphism $\nu\colon\widetilde S\to S$ is the normalization of $S$. Moreover, the morphism
$\widetilde C\to C$ is finite of degree $2$.

\item 
A degenerate cusp of type $T_{2,q,\infty}$ yields an
$A_{q-3}$-singularity on $\widetilde S$.

\item
The normalization is crepant  birational, i.e., $K_{\widetilde S}+\widetilde C=\nu^*K_S\sim0.$

\item
The surface $S$ has at most two degenerate cusps along $C$.
\end{enumerate}
\end{proposition}

\begin{proof}
Let $\bar\nu\colon\bar S\to S$ be the normalization
and let $\bar C\subset\bar S$ be the conductor. Since $S$ is slc, its conductor subscheme on $S$
is the reduced curve $C$, with ideal $\mathcal I_{C/S}$. Moreover, $\bar C\to C$ is a finite morphism of degree two
(see \cite[5.1--5.2]{K13}).

We first show that $\bar S$ is Gorenstein.
The normalization is smooth over $C$ away from the
degenerate cusps by Proposition~\ref{prop:slc_local_models}. Away from $C$, it is an isomorphism, and $S$ is
Gorenstein because it is a hypersurface. It remains to consider a degenerate cusp of type
$T_{2,q,\infty}$, locally analytically given by
$x^2+y^q-xyz=0$. Adjoining $u=x/y$ gives the finite birational morphism $\{u^2-uz+y^{q-2}=0\}.$
This has an $A_{q-3}$-singularity and is therefore
normal, so it is the normalization near the cusp. Thus $\bar S$ is Gorenstein.
By \cite[Theorem~1.4]{BS16}, blowing up
the conductor ideal $\mathcal I_{C/S}$
yields the normalization.
Hence $\widetilde S\cong\bar S$ and
$\widetilde C=\bar C$.
This proves \textup{(1)} and \textup{(2)}.

Furthermore, we have $K_{\widetilde S}+\widetilde C=\nu^*K_S$ by~\cite[Equation 5.7.3]{K13}.
Since $S\subset\mathbb P^3$ is a quartic,
$K_S\sim0$, proving \textup{(3)}.

At a degenerate cusp, the conductor on the normalization
is given by $y=0$ and $u(u-z)=0$, so it has an ordinary node. Away from the degenerate cusps, the conductor
is smooth. Thus $\widetilde C$ is a reduced
nodal curve. Thus $\widetilde C$ is a reduced nodal curve whose nodes
correspond bijectively to the degenerate cusps of $S$. Adjunction gives
$\omega_{\widetilde C}\cong\mathcal O_{\widetilde C}$,
hence $p_a(\widetilde C)=1$.
Let $N$ be the number of its nodes, $r$ the number
of its irreducible components, and $g_i$ their
geometric genera. Then
$1=p_a(\widetilde C)
  =\sum_{i=1}^r g_i+N-r+1$
(see \cite[Chapter~IV, Exercise~1.8]{Har77}).
Therefore $N\le r$.
Since $\widetilde C\to C$ is finite of degree two
and $C$ is irreducible, $r\le2$.
Hence $N\le2$, proving \textup{(4)}.
\end{proof}

Suppose that the double curve $C$ is a line. Then the exceptional
divisor of the blow-up along $C$ is
$E=\mathbb P(N_{C/\mathbb P^3})
\cong \mathbb P^1\times\mathbb P^1$.
By Proposition~\ref{prop:conductor-normalization}, the conductor
satisfies
$\widetilde C=\widetilde S\cap E
\in |\mathcal O_E(2,2)|$.
We distinguish two cases according to the geometry of $\widetilde C$.
If $\widetilde C$ is the union of two curves of bidegree $(1,1)$,
we say that $S$ is of type \textup{\quarticfamilylink{B1-alpha}{(B.1;$\alpha$)}};
otherwise, we say that $S$ is of type \textup{\quarticfamilylink{B1-beta}{(B.1;$\beta$)}}.
Together with the conic non--normal locus case \textup{\quarticfamilylink{B2-alpha}{(B.2;$\alpha$)}},
these are referred to as the three \emph{Urabe families}.

We call equations~\eqref{eq:B1-normal} and~\eqref{eq:B2-normal}
the \emph{Urabe normal forms} for non-normal slc quartic surfaces.
In case \textup{\quarticfamilylink{B1}{(B.1)}}, the coefficients
$a_2,b_2,c_2,d_3,e_3$ and $f_4$ involve
$3,3,3,4,4$ and $5$ parameters, respectively.
Hence equation~\eqref{eq:B1-normal} is parametrized by
an affine space $\mathcal U^{\quarticfamilylink{B1}{(B.1)}}$ of dimension $22$.
Similarly, the coefficient space $\mathcal U^{\quarticfamilylink{B2}{(B.2)}}$
associated with \eqref{eq:B2-normal} has dimension $14$.

\begin{definition}
\label{def:Urabe-coefficient-spaces}
\em
For the Urabe normal forms \textup{\quarticfamilylink{B1}{(B.1)}} and \textup{\quarticfamilylink{B2}{(B.2)}}, let $S_t=V(F_t)$ for $t\in \mathcal{I}^{\quarticfamilylink{B1}{(B.i)}}$ for $i\in \{\quarticfamilylink{B1}{1},\quarticfamilylink{B2}{2}\}$. Let $C$ denote the fixed line and smooth conic, respectively.
For each parameter $t$, let $C_t$ denote the reduced
non-normal locus of $S_t$.

We define the corresponding \emph{Urabe coefficient spaces}
as follows:
\[
\begin{aligned}
\quarticfamilytarget{B1-alpha}\mathcal I^{\quarticfamilylink{B1-alpha}{(B.1;\alpha)}}
&:=
\left\{
t\in\mathcal U^{\quarticfamilylink{B1}{(B.1)}}
\ \middle|\
\begin{array}{l}
S_t\text{ is irreducible and slc, and its conductor}\\
\widetilde C_t\text{ is the union of two curves of bidegree }(1,1)
\end{array}
\right\},\\
\quarticfamilytarget{B1-beta}\mathcal I^{\quarticfamilylink{B1-beta}{(B.1;\beta)}}
&:=
\left\{
t\in\mathcal U^{\quarticfamilylink{B1}{(B.1)}}
\ \middle|\
\begin{array}{l}
S_t\text{ is irreducible and slc, and its conductor}\\
\widetilde C_t\text{ is not the union of two curves of bidegree }(1,1)
\end{array}
\right\},\\
\quarticfamilytarget{B2-alpha}\mathcal I^{\quarticfamilylink{B2-alpha}{(B.2;\alpha)}}
&:=
\left\{
t\in\mathcal U^{\quarticfamilylink{B2}{(B.2)}}
\ \middle|\
\begin{array}{l}
S_t\text{ is irreducible and slc},\\
C_t=C,\quad
\operatorname{mult}_p(S_t)=2\quad\text{for every }p\in C
\end{array}
\right\}.
\end{aligned}
\]
\end{definition}

We set $\mathcal I^{\quarticfamilylink{B1}{(B.1)}}
:=
\mathcal I^{\quarticfamilylink{B1-alpha}{(B.1;\alpha)}}
\cup
\mathcal I^{\quarticfamilylink{B1-beta}{(B.1;\beta)}}$.
Irreducibility and the multiplicity condition are open.
The slc condition is open by inversion of adjunction
and semicontinuity of log canonical thresholds
(see \cite[Example~9.5.41]{Laz04b}). Thus $\mathcal I^{\quarticfamilylink{B2-alpha}{(B.2;\alpha)}}$ is open in
$\mathcal U^{\quarticfamilylink{B2}{(B.2)}}$.
Moreover, Proposition~\ref{Prop.B.1.betacondition}
shows that $\mathcal I^{\quarticfamilylink{B1-beta}{(B.1;\beta)}}$ is smooth.

\begin{definition}
\label{def:marked-Urabe-family}
\em Let $i$ denote $\quarticfamilylink{B1-alpha}{1;\alpha}$ or $\quarticfamilylink{B1-beta}{1;\beta}$ or $\quarticfamilylink{B2-alpha}{2;\alpha}$. Let
$\mathcal I^{\quarticfamilylink{B1-alpha}{(B.i)}}$ be the corresponding Urabe coefficient space. Let
\[
\mathcal S^{(i)}
:=
\{(x,t)\mid F_t(x)=0\}
\subset
\mathbb P^3\times\mathcal I^{\quarticfamilylink{B1-alpha}{(B.i)}}
\]
be the \emph{tautological family}. Its fiber over
$t\in\mathcal I^{\quarticfamilylink{B1-alpha}{(B.i)}}$ is $S_t=\{F_t=0\}$.

The double curve $C$ is fixed by the Urabe normal form, and we define the \emph{relative double curve} as $\mathcal C^{(i)}
:=
C\times\mathcal I^{\quarticfamilylink{B1-alpha}{(B.i)}}
\subset
\mathcal S^{(i)}.$

\begin{enumerate}
\item
For a morphism $f\colon B\rightarrow\mathcal I^{\quarticfamilylink{B1-alpha}{(B.i;\alpha)}}$, set $\mathcal S_B
:=
\mathcal S^{(i)}\times_{\mathcal I^{\quarticfamilylink{B1-alpha}{(B.i)}}}B$ and $
\mathcal C_B
:=
\mathcal C^{(i)}\times_{\mathcal I^{\quarticfamilylink{B1-alpha}{(B.i)}}}B
\cong C\times B.$
We call $(\mathcal S_B,\mathcal C_B)\rightarrow B$
an \emph{Urabe family of type}
\textup{\quarticfamilylink{B1-alpha}{(B.1;$\alpha$)}}, \textup{\quarticfamilylink{B1-beta}{(B.1;$\beta$)}}
\textup{\quarticfamilylink{B2-alpha}{(B.2;$\alpha$)}}.

\item
Let $q_1\geq q_2\geq3$. An \emph{equisingular
$\mathcal T^{\quarticfamilylink{B1-alpha}{(B.i)}}_{q_1,q_2}$-family}
is an Urabe family $(\mathcal S_B,\mathcal C_B)\rightarrow B$ with sections $\sigma_1$ and $\sigma_2\colon B\rightarrow\mathcal C_B$ such that, for every $b\in B$, the points
$\sigma_1(b)$ and $\sigma_2(b)$ are the degenerate cusps of
$S_b$, and
\[
(S_b;\sigma_j(b))\cong T_{2,q_j,\infty},\quad\text{with}\quad j=1,2.
\]
Similarly, for $q\geq3$, an equisingular
$\mathcal T^{\quarticfamilylink{B1-alpha}{(B.i)}}_{q}$-family is defined by requiring a
unique degenerate cusp of type $T_{2,q,\infty}$ on each fiber.
\end{enumerate}
\end{definition}
If $\mathcal{T}_{q_1,q_2}^{\quarticfamilylink{B1-alpha}{(B.i)}}$ has $m$ irreducible components, we denote them by $\mathcal T^{\quarticfamilylink{B1-alpha}{(B.i)},I}_{q_1,q_2},
\ldots,
\mathcal T^{\quarticfamilylink{B1-alpha}{(B.i)},m}_{q_1,q_2}.$ By convention, if $F \in \mathcal{I}^{\quarticfamilylink{B1-alpha}{(B.i)}}$ (respectively, $F \in \mathcal{T}_{q_1,q_2}^{\quarticfamilylink{B1-alpha}{(B.i)}}$), then for the quartic surface $S := \{F=0\} \subset \mathbb{P}^3$, we may also write $S \in \mathcal{I}^{\quarticfamilylink{B1-alpha}{(B.i)}}$ (respectively, $S \in \mathcal{T}_{q_1,q_2}^{\quarticfamilylink{B1-alpha}{(B.i)}}$). Next, we determine the range of possible values of $q_1$ and $q_2$.

\begin{proposition}
\label{prop:non-normal-quartic-lattices}
Let $S\subset\mathbb P^3$ be an irreducible non-normal
slc quartic surface with a degenerate cusp of type $T_{2,q,\infty}$.
Assume that its reduced double curve $C$ is either a line or a
smooth conic and that $S$ has multiplicity two at every point of $C$.
Then $q\le 10$ if $C$ is a line
$q\le 7 $ if $C$ is a smooth conic.
\end{proposition}

\begin{proof}
Let $Y\to\widetilde S$ be the minimal resolution, let $H$ be the pullback of
$\mathcal O_S(1)$, and put $\Lambda:=\langle K_Y,H\rangle^{\perp}\subseteq\operatorname{Pic}(Y).$
Then by Urabe's classification, when $C$ is a line, the surface $Y$ is isomorphic to the
blow up of $\mathbb P^2$ at nine points $p_1,\ldots,p_9$ (see~\cite[Table~3]{Duc24}).
Let $h$ denote the pullback of the class of a line in $\mathbb P^2$,
and let $e_i$ denote the exceptional classes over $p_i$. In this basis, gives the polarization $H=4h-2e_1-\sum_{i=2}^9e_i$ and $K_Y=-3h+\sum_{i=1}^9e_i$ (see~\cite[Table~3]{Duc24}).

We now identify the lattice $\Lambda$. For
$r=ah-\sum_i b_ie_i$ the conditions $r\cdot K_Y=r\cdot H=0$ imply $b_1=a$ and
$\sum_{i\geq2}b_i=2a$, and then $r^{2}=-\sum_{i\geq2}b_i^{2}$, so
$r\mapsto(b_2,\dots,b_9)$ identifies $\Lambda$ with $D_8(-1)$. When $C$ is a
smooth conic, a similar argument identifies 
$\Lambda=K_Y^{\perp}\cong D_5(-1)$ (see e.g.\ \cite[Theorem~25.4]{Manin86}).

By Proposition~\ref{prop:conductor-normalization}, the exceptional curves over a
$T_{2,q,\infty}$-point span a copy of $A_{q-3}(-1)$
in $\Lambda$. Hence $q-3\leq\operatorname{rk}\Lambda$.

In the \quarticfamilylink{B1}{(B.1)} Urabe family, we have that $\Lambda\cong D_8(-1)$.
The case $q-3=8$ is impossible, since a full rank
inclusion $A_8\subseteq D_8$ would imply that $\frac{\det A_8}{\det D_8}=\frac{9}{4}$
is the square of an integer. Thus $q\leq10$.
Similarly, in the conic case, we have
$\Lambda\cong D_5(-1)$, and the case $q-3=5$
is excluded since
$\det A_5/\det D_5=6/4$ is not the square
of an integer. Hence $q\leq7$.
\end{proof}

We next show that marked Urabe families with fixed degenerate cusp data admit
simultaneous dlt modifications.

\begin{proposition}
\label{prop:simultaneous-urabe-dlt}
Let $B\subset\mathcal I^{\quarticfamilylink{B1-alpha}{(B.i)}}$ be a smooth irreducible
marked equisingular stratum of type
$\mathcal T^{\quarticfamilylink{B1-alpha}{(B.i)}}_q$ or
$\mathcal T^{\quarticfamilylink{B1-alpha}{(B.i)}}_{q_1,q_2}$.
Let $\mathcal S_B\subset\mathbb P^3\times B$
be the tautological Urabe family.
Then, \'etale locally on $B$, the pair
$(\mathbb P^3_B,\mathcal S_B)$ admits a simultaneous
$\mathbb Q$-factorial dlt modification.
\end{proposition}
\begin{proof}
Let $\pi\colon\mathcal X_1
:=\operatorname{Bl}_{\mathcal C_B}(\mathbb P^3_B)
\longrightarrow\mathbb P^3_B$
be the blow-up morphism, with exceptional divisor $E_B$,
and let $\widetilde{\mathcal S}_B$ be the strict transform
of $\mathcal S_B$.
Since $\mathcal C_B=C\times B$, both $\mathcal X_1$
and $E_B$ are smooth over $B$.
Since $\mathcal S_B$ has multiplicity two along
$\mathcal C_B$, the morphism $\pi$ is crepant.

By Proposition~\ref{prop:conductor-normalization},
the restriction $\pi|_{\widetilde{\mathcal S}_B}$
is finite and gives the normalization fiberwise,
with conductor $E_b\cap\widetilde S_b$.
Since $\widetilde{\mathcal S}_B$ is flat over $B$,
this restriction is a simultaneous normalization.
In particular, each marked degenerate cusp section $\sigma_j$ lifts
uniquely to a section $\tilde{\sigma_j}$ of
$\widetilde{\mathcal S}_B\to B$.

Note that the index $q_j$ is constant over $B$.
Moreover, the conductor has a non-degenerate node
along $\tilde{\sigma_j}$. 
Therefore, the calculation of Proposition~\ref{prop:conductor-normalization} applies relatively over $B$.
This gives us the local model $B_{q_j-22}$
along $\tilde{\sigma_j}$.
The pair induced on $\mathcal X_1$ by pull-back is dlt
away from the marked sections.
Therefore, Lemma~\ref{lem:relative-uv-sm-dlt}
gives a projective crepant dlt modification.
We compose this projective crepant dlt modification with $\pi$, getting the required simultaneous
$\qq$-factorial dlt modification of the initial pair, 
\'etale locally over $B$.
\end{proof}

\subsubsection{Quartics of type \quarticfamilylink{B1}{(B.1)}}

In this subsection, we study when the pair
$(\mathbb P^3,S)$ is of cluster type for
$S$ in the \textup{\quarticfamilylink{B1}{(B.1)}} Urabe families.
We show that $(\mathbb P^3,S)$ is of cluster type
if $S\in\mathcal I^{\quarticfamilylink{B1-alpha}{(B.1;\alpha)}}$, or if
$S\in\mathcal I^{\quarticfamilylink{B1-beta}{(B.1;\beta)}}$ has a degenerate cusp
of index $q\ge7$.
We also provide three examples
$S_1\in\mathcal T_{4,4}^{\quarticfamilylink{B1-beta}{(B.1;\beta)}}$,  $S_2\in\mathcal T_6^{\quarticfamilylink{B1-beta}{(B.1;\beta)},\mathrm{II}},$ $
S_3\in\mathcal T_{6,3}^{\quarticfamilylink{B1-beta}{(B.1;\beta)}}$ and 
such that none of the pairs $(\mathbb P^3,S_i)$
is of cluster type
(Proposition~\ref{prop:three-irrational-models}).

Alves da Silva, Figueroa, and the first author
studied whether a non-normal quartic surface
whose double curve is a line can define a
cluster type pair (see \cite[Example~7.2]{EFM26}).
Their example belongs to the \textup{\quarticfamilylink{B1-alpha}{(B.1;$\alpha$)}}
family. The following proposition extends their
construction to the entire family.

\begin{proposition}
\label{prop:11-splitting-cluster}
Let $S \in \mathcal{I}^{\quarticfamilylink{B1-alpha}{(B.1;\alpha)}}$. 
Then the pair $(\mathbb{P}^3, S)$ is of cluster type.
\end{proposition}

\begin{proof}
Write the Urabe normal form~\eqref{eq:B1-normal} as
\[
F
=
b_{2,2}(w,x,y,z)+d_3(y,z)w+e_3(y,z)x+f_4(y,z),
\]
where $b_{2,2}(w,x,y,z)
=
a_2(y,z)w^2+b_2(y,z)wx+c_2(y,z)x^2.$ Let $C=(y=z=0)$. On the exceptional divisor $E\cong\mathbb P^1_{w,x}\times\mathbb P^1_{y,z}$ 
of $\operatorname{Bl}_C\mathbb P^3$, the conductor is defined by $$\tilde{C}:=\{b_{2,2}=0\}\subset \mathbb P^1_{w,x}\times\mathbb P^1_{y,z}. $$
Moreover, since $S\in\mathcal I^{\quarticfamilylink{B1-alpha}{(B.1;\alpha)}}$, the conductor is the union of two curves of bidegree $(1,1)$. Since the conductor cover is finite, the two
$(1,1)$-curves are smooth. The local slc models
show that they meet transversely at two points.
We may therefore choose coordinates on the two
factors of $E$ so that
\[
b_{2,2}=(wy-xz)(wy-\lambda xz),
\quad\text{with}\quad \lambda\in\mathbb C^*\setminus\{1\}.
\]

We now work on the affine chart $D_{\mathbb P^3}(yz)$. After dehomogenizing by $z$, set $u:=wy-x$ and $v:=wy-\lambda x.$
Since $x=\frac{u-v}{\lambda-1}$ and $w=\frac{\lambda u-v}{(\lambda-1)y},$
 in these coordinates the equation of $S$ has the form
\[
F=uv+a(y)u+b(y)v+c(y)
\quad\text{with}\quad
a,b,c\in\mathbb C[y^{\pm1}].
\]
Now make the coordinate change $s:=u+b(y)$ and $ t:=v+a(y).$
Then, the equation becomes $F
=st+g(y)$ where $
g(y):=c(y)-a(y)b(y)\in\mathbb C[y^{\pm1}].$ This gives a torus embedding as in Example~\ref{ex:shear-instances}. This is a dense open subset of $\mathbb P^3\setminus S$.
Hence $(\mathbb P^3,S)$ is of cluster type by Lemma~\ref{clustertypeandtorusembedding}.
\end{proof}

\begin{proposition}\label{Prop.B.1.betacondition}
The Urabe family $\mathcal{I}^{\quarticfamilylink{B1-beta}{(B.1;\beta)}}$
is an open subset of $\mathcal{U}^{\quarticfamilylink{B1}{(B.1)}}$.
\end{proposition}

\begin{proof}
By Definition~\ref{def:Urabe-coefficient-spaces}, $I^{\quarticfamilylink{B1-beta}{(B.1;\beta)}}$ is the complement of $I^{\quarticfamilylink{B1-alpha}{(B.1;\alpha)}}$, so it is suffice to show $I^{\quarticfamilylink{B1-alpha}{(B,1;\alpha)}}$ is defined a closed subspace of $\mathcal{U}^{\quarticfamilylink{B1}{(B,1)}}$. Indeed, this locus is defined as a closed subspaces of $|\mathcal{O}_{\mathbb P^1\times \mathbb P^1}(2,2)|$ defined by the image of the multiplication map $|\mathcal{O}_{\mathbb P^1\times \mathbb P^1}(1,1)|\times |\mathcal{O}_{\mathbb P^1\times \mathbb P^1}(1,1)|$. By a dimension counting, it is a proper subset of $|\mathcal{O}_{\mathbb P^1\times\mathbb P^1}(2,2)|$. 
\end{proof}

\begin{proposition}\label{prop:B1beta-marked-family}
There exists a marked Urabe family $\mathcal S_B\to B$
of type \textup{\quarticfamilylink{B1-beta}{(B.1;$\beta$)}} representing every
quartic of this type with a marked degenerate cusp,
up to projective equivalence.
\end{proposition}

\begin{proof}
 The conductor has an ordinary node above $p_1$. After changing coordinates on the two factors of $E$, its equation is
\begin{equation}\label{eq:newformB.1.beta}
 (y^2+z^2)w_1^2+\lambda yzw_1x_1+\epsilon z^2x_1^2=0,
 \qquad D:=\lambda^2-4\epsilon\ne0.
\end{equation}
Here $\lambda$ is allowed to vanish, and $\epsilon=0$ is allowed. Indeed, the equation is
\[
 (y^2+a_1yz+a_2z^2)w_1^2
 +(b_1yz+b_2z^2)w_1x_1+c_2z^2x_1^2=0.
\]
The node condition is $b_1^2-4c_2\ne0$. The changes $y':= y+tz$ and $x_1':= x_1+uw_1$ remove $yzw_1^2$ and $z^2w_1x_1$ by a linear system with determinant $4c_2-b_1^2$. The coefficient of $z^2w_1^2$ is then nonzero, otherwise the conductor splits into two divisors of mutidegree $(1,1)$. By rescaling, we get an equation of the form~\eqref{eq:newformB.1.beta}.

The conductor is irreducible when $\epsilon\ne0$ and its discriminant is $z^2(Dy^2-4\epsilon z^2)$. Replacing $y,z$ by $y+L(w,x)$ and $z+M(w,x)$ for $L, M\in\mathbb C[w,x]$ removes the cubic terms $yw^3,zw^3,yw^2x$ and $zw^2x$. We obtain
\begin{equation}\label{eq:B1beta-global}
 F={}w^2(y^2+z^2)+\lambda xyzw+\epsilon x^2z^2
       +x^2w(\alpha y+\gamma z)+x^3(\beta y+\delta z)+aw^4+bw^3x+cw^2x^2+dwx^3+ex^4.
\end{equation}
Let $B$ be the locus of coefficient tuples in \eqref{eq:B1beta-global} satisfying $\lambda^2-4\epsilon\ne0$. Then, $B$ satisfies desired property.
\end{proof}

\begin{proposition}\label{prop:line-conductorstrata}
For the Urabe family of type \textup{\quarticfamilylink{B1-beta}{(B.1;$\beta$)}}, the following statements hold.
\begin{enumerate}
\item The locus with a marked $T_{2,6,\infty}$-singularity has two irreducible components.
\item For $6\geq q_1$, we have $4\geq q_2$. If $q_2=4$, then $q_1=4$.
\item For $q\le6$ and $6\ge q_1\ge q_2$, the strata $\mathcal T_q^{\quarticfamilylink{B1-beta}{(B.1;\beta)}}$ and $\mathcal T_{q_1,q_2}^{\quarticfamilylink{B1-beta}{(B.1;\beta)}}$ are smooth.
\item If $S\in\mathcal I^{\quarticfamilylink{B1-beta}{(B.1;\beta)}}$ has a degenerate cusp of type $T_{2,q,\infty}$ with $q\ge 7$, then $(\mathbb P^3,S)$ is of cluster type.
\end{enumerate}
\end{proposition}

\begin{proof}
We first determine the cusp type at $p_1=[0:0:1:0]$ of the surface $S:=\{F=0\}\subset \mathbb P^3_{w,x,y,z}$, where $F$ is as in \eqref{eq:B1beta-global}. Let $D:=\lambda^2-4\epsilon\neq0$ and on the open set $D_{\operatorname{Bl}_C(\mathbb P^3)}(y)$, we set $x=s, w=sw_1$ and $y=1$ so that the strict transform of $S$ is defined by
\[
\widetilde F|_{y=1} ={}w_1^2(1+z^2)+\lambda w_1z+\epsilon z^2 +s\bigl(\alpha w_1+\gamma zw_1+\beta+\delta z\bigr)+s^2\bigl(aw_1^4+bw_1^3+cw_1^2+dw_1+e\bigr).
\]
Since $\partial_s\widetilde F(0,0,0)|_{y=1}=\beta,$ we obtain a $T_{2,3,\infty}$-singularity when $\beta\neq0$. Let $\beta=0$, and put
\[
 \sigma:=\lambda\delta-2\epsilon\alpha \qquad\text{and}\qquad \tau:=2\delta-\lambda\alpha.
\]
Since the quadratic part $w_1^2+\lambda w_1z+\epsilon z^2$ is nondegenerate ($D\neq0$), we solve the equation $\partial_{w_1}\widetilde F|_{y=1}=\partial_z\widetilde F|_{y=1}=0$ by
\[
 w_1=-\frac{\sigma}{D}s+\left(\frac{2 \epsilon d}{D}+\frac{\gamma(2\lambda\sigma-D\delta)}{D^2}\right)s^2+O(s^3)\quad\text{and} \quad z=\frac{\tau}{D}s-\left(\frac{\lambda d}{D}+\frac{\gamma(2\lambda\tau+D\alpha)}{D^2}\right)s^2+O(s^3).
\]
Then
\begin{equation}\label{eq:B1beta-rhos}
 \widetilde F(s,w_1(s),z(s))|_{y=1} = \rho_2s^2-\frac{\sigma}{D}\rho_3s^3+\left( \frac{\sigma^2}{D^4}\rho_4 +\frac{\lambda\gamma\sigma}{D^2}\rho_3 +\frac{\epsilon}{D}\rho_3^2 \right)s^4+O(s^5),
\end{equation}
where
\begin{equation}\label{rhob_1beta}
 \rho_2:=e+\frac{\epsilon\alpha^2-\lambda\alpha\delta+\delta^2}{D}, \qquad \rho_3:=d+\frac{\gamma\tau}{D}\quad\text{and}\quad \rho_4:= D^2c+D\gamma^2+\tau^2.
\end{equation}
By Proposition~\ref{prop:conductor-normalization}, if $F$ has singularity of type $T_{2,q_1,\infty}$ at $p_1$ then the vanishing order in \eqref{eq:B1beta-rhos} at $p_1$ is $q_{1}-2$. For the strata with a unique degenerate cusp, we impose $\epsilon\ne0$. The strata of indices four and five and the first stratum of index six are given by the following conditions.
\[
\mathcal{T}_{4}^{\quarticfamilylink{B1-beta}{(B.1;\beta)}}:=\{\varepsilon\neq0,\beta=0,\rho_2\neq 0\},\quad
 \mathcal{T}_{5}^{\quarticfamilylink{B1-beta}{(B.1;\beta)}}:=\{\varepsilon\neq0,\beta=\rho_2=0,\rho_3\neq0,\sigma\neq 0\}
\]
\[
\quad \text{and}\quad
\mathcal{T}_{6}^{\quarticfamilylink{B1-beta}{(B.1;\beta)},\mathrm{I}}:=\{\varepsilon\neq0,\beta=\rho_2=\rho_3=0,\rho_4\ \neq 0,\sigma\neq 0\}.
\]
Suppose instead that $\sigma=0$. Then the coefficient of $s^3$ vanishes and that of $s^4$ is $\frac{\epsilon}{D}\rho_3^2$, so we obtain the irreducible stratum
\[
 \mathcal T_{6}^{\quarticfamilylink{B1-beta}{(B.1;\beta)},\mathrm{II}} :=\{\beta=\sigma=\rho_2=0,\ \epsilon\neq0,\ \rho_3\neq0\}.
\]
The two components of the marked $T_{2,6,\infty}$-locus are obtained by allowing $\epsilon=0$ in the first component.By equation~\eqref{eq:newformB.1.beta}, a second degenerate cusp occurs at $p_2=[0:0:0:1]$ precisely when $\epsilon=0$. In this case $D=\lambda^2\neq0$, so $\lambda\neq0$ and $\sigma=\lambda\delta$. We consider separately the cases $\delta\neq0$ and $\delta=0$. If $\delta\neq0$, the same calculation as in (i) gives $q_2=3$. This shows that $\mathcal T_{6,3}^{\quarticfamilylink{B1-beta}{(B.1;\beta)}}$ lies in the closure of $\mathcal T_6^{\quarticfamilylink{B1-beta}{(B.1;\beta)},\mathrm{I}}$.

Now suppose that $\delta=0$. If $\epsilon=0$ and $F$ is irreducible, then $e\neq0$. Otherwise, equation~\eqref{eq:B1beta-global} is divisible by $w$. We have $\delta=\epsilon=0$ and $\rho_2=e\neq0$, so equation~\eqref{eq:B1beta-rhos} gives $(q_1,q_2)=(4,4)$. On $\mathcal T_{6}^{\quarticfamilylink{B1-beta}{(B.1;\beta)},\text{II}}$ we have $\sigma=\rho_2=0$. If moreover $\epsilon=0$, then $\delta=0$ and $e=\rho_2=0$, so every such member is reducible. Consequently, no irreducible surface $S\in \mathcal T_{6}^{\quarticfamilylink{B1-beta}{(B.1;\beta)},\text{II}}$ of degenerate to the normal surface with a second degenerate cusp at $p_2$. In particular, the component $\mathcal T_{6}^{\quarticfamilylink{B1-beta}{(B.1;\beta)},\text{II}}$ admits no specialization within the irreducible quartic locus to a stratum with two degenerate cusps. Varying $\beta,\rho_2,\rho_3$, and $\epsilon$ in \eqref{eq:B1beta-global}--\eqref{rhob_1beta} shows that every equisingular stratum considered above with $q_1\le6$, other than $\mathcal T_{4,4}^{\quarticfamilylink{B1-beta}{(B.1;\beta)}}$ and $\mathcal T_6^{\quarticfamilylink{B1-beta}{(B.1;\beta)},\text{II}}$, contains $\mathcal T_{6,3}^{\quarticfamilylink{B1-beta}{(B.1;\beta)}}$ in its closure. Smoothness of the each strata follows from the Jacobian criterion. This proves (ii) and (iii).

Finally, to prove (iv), suppose that $q_1\geq 7$. If $\sigma=0$, then $\epsilon\neq0$, and $\rho_2=\rho_3=0$ means that $\epsilon z^2+\delta xz+ex^2$ has a double root and that $\partial_wF$ vanishes along it, so $S$ is singular along the line $\{w=2\epsilon z+\delta x=0\}\neq C$, which is excluded. Hence $\sigma \neq 0$ and $\rho_2=\rho_3=\rho_4=0$, that is, we have
\[
 e=\frac{D\alpha^2-\tau^2}{4D},\qquad d=-\frac{\gamma\tau}{D},\qquad c=-\frac{D\gamma^2+\tau^2}{D^2}.
\]
Substituting the previous relations into~\eqref{eq:B1beta-global} gives us
\[
 F|_{w=1}= y^2+z^2+\lambda xyz+\epsilon x^2z^2+\alpha x^2y+\gamma x^2z+\delta x^3z+a+bx+cx^2+dx^3+ex^4.
\]
We set $u:=2y+\lambda xz+\alpha x^2.$ 
Using $4(y^2+\lambda xyz+\alpha x^2y) = u^2-x^2(\lambda z+\alpha x)^2,$ 
we obtain
\[
 \begin{aligned}
 4F={}&u^2-Dx^2z^2+2\tau x^3z-\frac{\tau^2}{D}x^4 +4\gamma x^2z-\frac{4\gamma^2}{D}x^2-\frac{4\gamma\tau}{D}x^3\\
 &\quad+4z^2+4a+4bx -\frac{4\tau^2}{D^2}x^2.
 \end{aligned}
\]
Thus, setting $r:=z-\frac{\tau}{D}x-\frac{2\gamma}{D}$, the terms on the first line combine to give $-Dr^2x^2+u^2$. By eliminating $z$ from the second line, we have
\[
 4F=u^2-Dr^2x^2+g(r)x+h(r)\quad\text{with}\quad g(r)\quad\text{and}\quad h(r)\in\mathbb C[r].
\]
By Example~\ref{ex:shear-instances}, this gives a torus embedding $\mathbb G_m^3\hookrightarrow \mathbb P^3\setminus S$. Hence $(\mathbb P^3,S)$ is of cluster type. This proves (iv).
\end{proof}

\begin{proposition}
\label{prop:three-irrational-models}
There are quartics $F_{1}\in\mathcal T_{4,4}^{\quarticfamilylink{B1-beta}{(B.1;\beta)}},
F_{2}\in\mathcal T_{6,3}^{\quarticfamilylink{B1-beta}{(B.1;\beta)}}$ and $
F_{3}\in\mathcal T_{6}^{\quarticfamilylink{B1-beta}{(B.1;
\beta)},\text{II}}$
such that $Y_i:=\{u^4=F_i\}\subset \mathbb P^4_{u,w,x,y,z}$ is irrational for all $i=1,2,3$.
\end{proposition}

\begin{proof}
Consider the quartic $F$ in~\eqref{eq:B1beta-global} with $\lambda=a=1$
and $\alpha=\beta=b=c=0$.
Let $F_i$ be obtained by setting
$(\epsilon,\gamma,\delta,d,e)$ to $(0,1,0,0,1),
(0,0,1,0,-1),$ and $
(-3/4,1,0,1,0),$
respectively. Each $S_i:=V(F_i)$ is irreducible and smooth off $V(w,x)$, with slc
singularities along that line by inversion of adjunction.
By Proposition~\ref{prop:line-conductorstrata}, the corresponding surfaces
$S_1$, $S_2$, and $S_3$ belong to the indicated strata.

We write \eqref{eq:B1beta-global} as
\[
F=\left(wy+\frac12 xz\right)^2+Az^2+Bz+C,
\]
where $A=w^2+\left(\epsilon-\frac14\right)x^2,
B=x^2(\gamma w+\delta x),$ and $
C=w^4+dwx^3+ex^4.$ By Lemma~\ref{lem:conic-brauer-discriminant}, we have that $Y_i=(u^4=F_i)$ is birational to
a conic bundle with Brauer class
\[
\alpha=\left(A,\;Au^4-AC+\frac{B^2}{4}\right),
\]
whose discriminant curve is
given by $Au^4=AC-\frac{B^2}{4}.$
The discriminant curve is
birational to the smooth projective model of
\[
u^4=\frac{4-y^2+3y^4-y^6}{4-y^2},
\qquad
u^4=\frac{4-y^2-4y^4}{4-y^2},
\qquad
u^4=\frac{4y^5+y^4-4y^3+4y^2-4}{4(y^2-1)}.
\]
We verify the hypotheses of
Proposition~\ref{prop:fractional-kummer-irrationality}. In the notation
of that proposition, we have $\alpha=(R(y),R(y)u^4+H(y))$ with $R=A$ and $ H=\frac{B^2}{4}-AC,$
where $w=1$. In each case, $R$ and $H$ are square free and coprime, and
$\deg R=2$. Since $\deg R$ is even, condition (i) of Proposition~\ref{prop:fractional-kummer-irrationality} holds.  Let $B_4$ and $B_2$ denote the numbers
of branch points of ramification indices
$4$ and $2$, respectively, for
$\rho:C_i\to\mathbb P^1$.
The degrees of the numerators above are $6,4$ and $5$, and the degree of
the denominator is $2$ in each case. These give $(B_4,B_2)=(8,0),(6,1)$ and $(8,0)$
respectively.
Hence $B_4\geq4$ and
$3B_4+2B_2\in\{20,24\}$, so the ramification
conditions of Proposition~\ref{prop:fractional-kummer-irrationality} are satisfied.
Therefore $Y_1$, $Y_2$, and $Y_3$ are
irrational. 
\end{proof}

\subsubsection{Quartics of type \quarticfamilylink{B2-alpha}{(B.2;$\alpha$)}} 
In this subsection, we study when the pair
$(\mathbb P^3,S)$ is of cluster type for
$S\in\mathcal I^{\quarticfamilylink{B2-alpha}{(B.2;\alpha)}}$.
We show that $(\mathbb P^3,S)$ is of cluster type
for $S\in\mathcal T_{q_1,q_2}^{\quarticfamilylink{B2-alpha}{(B.2;\alpha)}}$
with $q_1\ge4$ (Proposition~\ref{prop:B2-low-cusp-strata}),
and provide an example
$S\in\mathcal T_{3,3}^{\quarticfamilylink{B2-alpha}{(B.2;\alpha)}}$
that is not of cluster type
(Proposition~\ref{prop:two-cusp-irrational}).
The quartics in \eqref{eq:B2-marked-form}
are type III analogues of the type II degenerations
studied by Jones \cite[\S 3.3]{Jones26}.
\begin{proposition}
\label{prop:B2-low-cusp-strata}
There exist an open subset $B\subset\mathbb A^6$ and a marked
Urabe family $\mathcal S_B\to B$ such that every marked quartic
of type \textup{\quarticfamilylink{B2-alpha}{(B.2;$\alpha$)}} is projectively equivalent to a
marked fiber of this family.

\begin{enumerate}
\item[(i)] The strata $\mathcal T_3^{\quarticfamilylink{B2-alpha}{(B.2;\alpha)}}$ and
$\mathcal T_{3,3}^{\quarticfamilylink{B2-alpha}{(B.2;\alpha)}}$ in $B$
are smooth and irreducible.

\item[(ii)] We have
$\mathcal T_{3,3}^{\quarticfamilylink{B2-alpha}{(B.2;\alpha)}}
\subset\overline{\mathcal T_3^{\quarticfamilylink{B2-alpha}{(B.2;\alpha)}}}$,
where the closure is taken in $B$.

\item[(iii)] If $S\in\mathcal I^{\quarticfamilylink{B2-alpha}{(B.2;\alpha)}}$ has a degenerate cusp
of type $T_{2,q,\infty}$ with $q\ge4$, then $(\mathbb P^3,S)$ is
of cluster type.
\end{enumerate}
\end{proposition}

\begin{proof}
Starting from \eqref{eq:B2-normal}, write $F=A^2+w^2R$, where
\[
A=Q+w(\ell_xx+\ell_yy+\ell_zz)+\kappa_0w^2,
\qquad R=b_2-\tfrac14a_1^2.
\]
The change of coordinates
$(x',y',z')=(x+\ell_zw,y-\ell_yw/2,z+\ell_xw)$ fixes $C$.
After relabeling, we have $A=Q+\kappa w^2$ with $Q=xz-y^2$.
On the open subset of $\operatorname{Bl}_C\mathbb P^3$ where
$q=Q/w$ is regular, the strict transform is given by
\begin{equation}\label{B.2strict}
\widetilde F=(q+\kappa w)^2+R=0.
\end{equation}

We parametrize the curve $C$ by $[0:s^2:st:t^2]$.
Then, the cusp appears at the point $p_1=[0:0:0:1]$.
By Proposition~\ref{prop:conductor-normalization},
the cusp is a double zero of $R|_C$. 
We may perform a change of coordinates on $C$, fixing $p-1$, and rescale the variable $w$ to obtain $R|_C=s^2(\tau s^2-t^2)=(\tau x^2-y^2)|_C$. 
Therefore, we have 
$R=\tau x^2-y^2+cQ+wL$. 
We replace $\kappa$ by $\kappa+c/2$
and $L$ by $L-(\kappa c+c^2/4)w$ to obtain the following equality:
\begin{equation}\label{eq:B2-marked-form}
F=(Q+\kappa w^2)^2+w^2(\tau x^2-y^2)
  -w^3(\alpha x+\beta y+\gamma z+\delta w).
\end{equation}
Let $B\subset\mathbb{A}^6_{\kappa,\tau,\alpha,\beta,\gamma,\delta}$
be the open locus where $V(F)\in\mathcal I^{\quarticfamilylink{B2-alpha}{(B.2;\alpha)}}$.
The corresponding family, marked at $p_1$, has the required properties. 

The point $p_1$ is always a double zero of $R|_C$. 
Moreover, when $\tau=0$, the point $p_2=[0:1:0:0]$
is also a double zero of $R|_C$.
The strict transform is smooth over $p_1$
exactly when $\gamma \neq 0$. Moreover, the strict transform is smooth over $p_2$ precisely when
$\alpha \neq 0$ if $\tau=0$. Indeed, the respective linear terms are $-\gamma w$ and $-\alpha w$. 
By Proposition~\ref{prop:conductor-normalization}.\textup{(2)}, this gives
\begin{equation}\label{taugammaalpha}
\mathcal T_3^{\quarticfamilylink{B2-alpha}{(B.2;\alpha)}}=B\cap\{\tau\ne0,\ \gamma\ne0\},
\qquad
\mathcal T_{3,3}^{\quarticfamilylink{B2-alpha}{(B.2;\alpha)}}
=B\cap\{\tau=0,\ \alpha\gamma\ne0\}.
\end{equation}

For (iii), mark the cusp of index $q\ge4$. The strict transform
is singular over $p_1$, so $\gamma=0$. On $D_{\mathbb P^3}(w)$,
setting $w=1$ gives
\[
F|_{w=1}=(xz-y^2+\kappa)^2-\left(y+\frac{\beta}{2}\right)^2
+\tau x^2-\alpha x-\delta+\frac{\beta^2}{4}.
\]
Example~\ref{ex:shear-instances}.\textup{(i)} gives an open embedding
$\mathbb{G}_m^3\hookrightarrow\mathbb{P}^3\setminus S$.
Therefore, the pair $(\pp^3,S)$ is of cluster type by Lemma~\ref{clustertypeandtorusembedding}.
\end{proof}

\begin{proposition}
\label{prop:two-cusp-irrational}
There is a quartic polynomial 
$F\in \mathcal{T}_{3,3}^{\quarticfamilylink{B2-alpha}{(B.2; \alpha)}}$ such that $Y:=\{u^4=F\}\subset \mathbb P^4_{u,w,x,y,z}$ is irrational.
\end{proposition}

\begin{proof}
Let \begin{equation}\label{eq:B.2example}F:=Q^2+wyQ+w^3(x+z)\qquad \text{with}\qquad Q:=xz-y^2 \end{equation} With $Q=wq$, the strict transform
$q(q+y)+w(x+z)=0$ is smooth at both conductor nodes since
$x+z\ne0$ there, so both cusps have type $T_{2,3,\infty}$. The surface $S:=V(F)$ is irreducible and smooth off $V(w,Q)$. Hence $S$ has slc
singularities.
We show that $Y$ is irrational. On the chart $D_{\mathbb P^3}(w)$,
we have $xz-y^2=Q$ and $x+z+Qy=(u^4-Q^2)$. 
Eliminating $z$ and homogenizing with a new variable $v$,
we obtain the conic bundle
\begin{equation}\label{conicbundleB.2}
\left\{
([x:y:v],(Q,u))\ \middle|\
(2y+Qx)^2+(4-Q^2)x^2
-4(u^4-Q^2)xv+4Qv^2=0
\right\}
\subset \mathbb P^2_{x,y,v}\times\mathbb A^2_{Q,u}.
\end{equation}

By Lemma~\ref{lem:conic-brauer-discriminant} the Brauer class is represented as
$\alpha:=\bigl(4-Q^2,\ (u^4-Q^2)^2-Q(4-Q^2)\bigr)$ and the discriminant curve $\Gamma$ can be computed as
\[\Gamma:=\{(Q,u)\mid (u^4-Q^2)^2-Q(4-Q^2)=0\}\subset \mathbb A^2_{Q,u}. \]

Let $r:=u^4-Q^2$ and let $C$ be the smooth projective model of
$\Gamma$. Then $C$ is a quartic cover of the elliptic curve $E:=\{[Q:r:z]\mid r^2z=Q(4z^2-Q^2)\}\subset  \mathbb P^2_{Q,r,z},$ with $p\colon C\longrightarrow E$ defined by
$u^4=f:=Q^2+r$.

On $E$, the function $f$ has four simple zeros, namely $P_0:=[0:0:1]$ and the three points with
$Q^3+Q^2-4=0$, $r=-Q^2$, and a pole of order $4$ at
the point $O$ at infinity. In particular, $f$ is not a square in $\mathbb C(E)$, so $\Gamma$ is irreducible.
Hence $p$ is totally ramified at four
points and unramified elsewhere, so Riemann--Hurwitz gives $g(C)=7.$
Moreover $p^{-1}(O)=\{O_1,\ldots,O_4\}$ and $(u)_\infty=O_1+\cdots+O_4,$
so $u$ defines a degree four map $u\colon C\to\mathbb P^1$.

We consider the quotient of $C$ by the involution
given by $u\mapsto -u$. Let $D$ be such involution.
If we write $s=u^2$, then $D$ is the smooth plane quartic given by $(s^2-Q^2)^2+Q^3z-4Qz^3=0$. 
We claim that the curve $C$ is not hyperelliptic.
Indeed, its hyperelliptic involution would commute with $u\mapsto -u$. If the two involutions coincide, then we would have $D\cong \mathbb{P}^1$.
Otherwise, the hyperelliptic involution would descend to a non-trivial involution of $D$.
The quotient of $D$ by the hyperelliptic involution would be a rational curve. 
This contradicts the fact that a smooth plane quartic is non-hyperelliptic.  Since $g(C)=7$, the curve $C$ is not
birational to a plane quintic. If $C$ was trigonal, the degree
$3$ map and $u\colon C\to\mathbb{P}^1$ would give $g(C)\le (3-1)(4-1)=6$ by the Castelnuovo--Severi inequality, a contradiction.

By the residue formula~\eqref{residueformula}
the residue of $\alpha$ along $\Gamma$ is represented by
$4-Q^2$. On $C$, writing $r=u^4-Q^2$, we have $r^2=Q(4-Q^2).$ Thus $4-Q^2$ and $Q$ represent the
same class in $\mathbb C(C)^\times/\mathbb C(C)^{\times2}$. Hence the double cover $\tau\colon\widetilde{C}\rightarrow C$ is obtained by
adjoining $\sqrt Q$.

We claim that the double cover $\tau\colon\widetilde{C}\rightarrow C$ is connected and \'etale.
Since $\mathbb C(C)=\mathbb C(E)(u)$ with $u^4=f$
is a cyclic extension of degree $4$, its unique quadratic
subextension is $\mathbb C(E)(u^2)=\mathbb C(E)(\sqrt f).$ Thus, if $Q$ were a square in $\mathbb C(C)$, then either $Q$
would be a square in $\mathbb C(E)$, or $Q/f$ would be a square in
$\mathbb C(E)$. The first is impossible since $\operatorname{div}_E(Q)=2P_0-2O$ and $P_0-O$ is not principal, while the second is impossible since
$Q/f$ has a simple zero at $P_0$. Hence $Q$ is not a square in
$\mathbb C(C)$, so the double cover is connected.

The double cover $\tau$ is \'etale because the quadratic extension obtained by adjoining $\sqrt Q$ is
ramified precisely at the points $P$ where $\operatorname{ord}_P(Q)$ is odd. Since $P_0$ is a zero of $f$, it is
totally ramified under $p$; let $\widetilde P_0$ be the point of $C$ over it. Then
\[
\operatorname{div}_C(Q)
=
p^*\operatorname{div}_E(Q)
=
8\widetilde P_0-2(O_1+\cdots+O_4),
\]
and all coefficients are even. Hence the double cover is \'etale.
S
On $\mathbb P^1_Q\times\mathbb P^1_u$, the only possible
ramification divisors besides $\Gamma$ are $Q=\pm2$,
$u=\{\infty\}$, and $Q=\infty$. Along $Q=\pm2$, the second
entry of $\alpha$ specializes to $(u^4-4)^2$.
Along $u=\{\infty\}$, the second entry has a pole of order $8$,
while along $Q=\infty$, the two entries have poles of
orders $2$ and $4$, respectively.
Thus \eqref{residueformula} gives trivial residues
away from $\Gamma$.
By Lemma~\ref{lem:standardization-after-resolution},
$Y$ is birational to a standard conic bundle over a
smooth projective rational surface whose discriminant
curve is isomorphic to $C$.
Since $C$ is neither hyperelliptic, trigonal, nor
birational to a plane quintic,
Proposition~\ref{Thm:Shokurov} implies that
$Y$ is irrational.

\end{proof}

\subsection{Quartic surfaces with reducible non--normal locus}
\label{subsec:reducible-non-normal}

Recall from Proposition~\ref{prop:reduciblenonnormallocus} that the remaining
cases with reducible non-normal locus are
the following three configurations:

\begin{itemize}
  \item[\quarticfamilylink{B2-beta}{(B.2;$\beta$)}] \quarticfamilytarget{B2-beta}the union of two distinct lines meeting in a point;
  \item[\quarticfamilylink{B3-beta}{(B.3;$\beta$)}] \quarticfamilytarget{B3-beta}the union of a smooth conic $C$ and a line $L$ not contained
        in the plane spanned by $C$, such that $L$ meets $C$ in exactly one point; and 
  \item[\quarticfamilylink{B4-beta}{(B.4;$\beta$)}] \quarticfamilytarget{B4-beta}a chain of three lines $L_1 \cup L_2 \cup L_3$, that is,
        $L_1 \cap L_2 \neq \emptyset$ and $L_2 \cap L_3 \neq \emptyset$, in which
        the lines $L_1$ and $L_3$ are not coplanar.
\end{itemize}
 For each $i\in\{\quarticfamilylink{B2-beta}{2},\quarticfamilylink{B3-beta}{3},\quarticfamilylink{B4-beta}{4}\}$, let $\mathcal I^{\quarticfamilylink{B2-beta}{(B.i;\beta)}}$ be the
family of irreducible slc quartic surfaces in $\mathbb P^3$
with reduced non-normal locus of type $\textup{\quarticfamilylink{B2-beta}{(B.$i$;$\beta$)}}$
and multiplicity at most two at every point. In what follows, we study each of these families in turn.
\begin{lemma}\label{lem:lc-quadratic}
Let $(X,S)$ be a log canonical pair with $X$
a smooth threefold and $S$ reduced.
If a local equation of $S$ has weighted order four
with respect to the weights $(2,1,1)$,
then its weighted initial form is reduced.
\end{lemma}

\begin{proof}
The weighted blow-up with weights $(2,1,1)$
is crepant for $(X,S)$.
By adjunction to the exceptional divisor is reduced so as the weighted initial form is.
\end{proof}
\begin{proposition}\label{prop:ClusterB.2.beta}
Let $S\in\mathcal I^{\quarticfamilylink{B2-beta}{(B.2;\beta)}}$ and assume that
$(\mathbb P^3,S)$ has coregularity zero. If $S$ has a point of
multiplicity three or at least two distinct degenerate cusps,
then $(\mathbb P^3,S)$ is of cluster type.
\end{proposition}

\begin{proof}
The multiplicity three case follows from Theorem~\ref{introthm3}.
We may therefore assume that all points of $S$ have multiplicity
at most two. Choose coordinates $[w:x:y:z]$ in which the double lines are
$L_1=V(x,y)$ and $L_2=V(x,z)$, meeting at $p=[1:0:0:0]$.
Since $F\in(x,y)^2\cap(x,z)^2$, we can write
\[
 F=\alpha y^2z^2+xyz\ell_1+x^2q_2,
\]
where $\ell_1$ is linear and $q_2$ is quadratic. Irreducibility
gives $\alpha\ne0$, and $\operatorname{mult}_p S=2$ implies that
the coefficient of $x^2w^2$ is nonzero. Rescaling $F$ and $x$,
we obtain
\begin{equation}\label{eq:B.2.beta}
 F=x^2w^2+\lambda xyzw+y^2z^2+xyz(ay+bz)
   +x^2(cy^2+dz^2+eyz+fxy+gxz+hx^2).
\end{equation}
We consider the open chart $D_{\mathbb{P}^3}(w)$.
We give weights $(2,1,1)$ to the variables $(x,y,z)$.
In these weights, the weighted initial form of the previous equation is $x^2+\lambda xyz+y^2z^2$. 
Therefore, we have $k:=1-\lambda^2/4\neq 0$ by Lemma~\ref{lem:lc-quadratic}. Completing the square in $x$ leaves the monomial $y^2z^2$ times a unit.
Therefore, the singularity at $p$ is of type  $T_{2,\infty,\infty}$.

The computation in the proof of
Proposition~\ref{prop:conductor-normalization} shows that
$L_1\setminus\{p\}$ has a degenerate cusp precisely when
$(\lambda w+b)^2-4(w^2+d)$ has a double root. 
This is equivalent to imposing that 
$b^2-4kd=0$. Analogously, $L_2\setminus \{p\}$ has a degenerate cusp precisely when $a^2-4kc=0$.
We can relabel the lines so that
a second degenerate cusp lies on $L_1$. Therefore, we get
$ky^2+by+d=t(y)^2$.

On the chart $D_{\mathbb P^3}(x)$ with $t\ne0$, put $w'=w+\frac{\lambda yz}{2}$ and $ z'=tz+\frac{ay^2+ey+g}{2t}.$
Then
\[
 F|_{x=1}=(w'+\mathrm{i}z')(w'-\mathrm{i}z')+G(t)
 \quad\text{for some }G\in\mathbb C[t^{\pm1}].
\]
Example~\ref{ex:shear-instances}\textup{(ii)} gives an open
embedding $\mathbb G_m^3\hookrightarrow\mathbb P^3\setminus S$.
Hence $(\mathbb P^3,S)$ is of cluster type.
\end{proof}

Set $k=1-\lambda^2/4,
 \delta_1=a^2-4kc$ and $
 \delta_2=b^2-4kd.$
Let $B=D(k\delta_1\delta_2)$ in the affine space $\mathbb A^9$
with coordinates $(\lambda,a,b,c,d,e,f,g,h)\subset \mathcal{I}^{\quarticfamilylink{B2-beta}{(B.2;\beta)}}$, and let
$\mathcal S\subset\mathbb P^3\times B$ be the family defined by
\eqref{eq:B.2.beta}. 
By Lemma~\ref{lem:lc-quadratic} and
Proposition~\ref{prop:ClusterB.2.beta}, if $S_t\in \mathcal{S}_B$ then 
$S$ of is irreducible slc quartics with a unique degenerate cusp $T_{2,\infty,\infty}$. For the next statement, we prove this $B$ admit simlteneous $\mathbb Q$-factorial dlt modification.

\begin{proposition}\label{prop:two-line-dlt}
Let $S\in\mathcal I^{\quarticfamilylink{B2-beta}{(B.2;\beta)}}$ have a unique degenerate cusp,
of type $T_{2,\infty,\infty}$.
Then $S$ lies in a smooth irreducible family admitting a
simultaneous $\mathbb Q$-factorial dlt modification étale locally
on the base.
\end{proposition}

\begin{proof}
Let $\mathcal S_B\subset\mathbb P^3\times B$ be the tautological family.
We blow up $L_1\times B$ and then the strict transform of $L_2\times B$. 
We write $\pi\colon \mathbb{P}^3\times B$ for the resulting morphism.
Let $E_1$ and $E_2$ be the exceptional divisors.
Let $\mathcal{S}_{\mathcal{Y}}$ be the strict transform of $\mathcal{S}$. 
The morphism 
$\pi \colon (\mathcal{Y},E_1+E_2+\mathcal{S}_{\mathcal{Y}})\rightarrow (\pp^3,B)$ is crepant. 
Moreover, the variety $\mathcal{Y}$ is smooth over $B$.
Away from the exceptional divisors, the fibers
are plt by inversion of adjunction and Proposition~\ref{uniquedegeneratecuspsurfaces}. Along each exceptional divisor $E_1$ and $E_2$, we are in the situation of Lemma~\ref{lem:relative-uv-sm-dlt}. 
Thus, the family admits a simultaneous $\qq$-factorial dlt modification \'etale over the base. 
\end{proof}

\begin{proposition}\label{prop:two-line-irrational}
There exists a quartic surface $S\in\mathcal I^{\quarticfamilylink{B2-beta}{(B.2;\beta)}}$
with exactly one degenerate cusp, of type $T_{2,\infty,\infty}$,
such that $(\mathbb P^3,S)$ is not of cluster type.
\end{proposition}

\begin{proof}
Let
$F=x^2w^2+(y^2-x^2)z^2+x^3z+x^2y^2+x^4 \in \cc[x,y,z,w]$.
Using the notation of Proposition~\ref{prop:ClusterB.2.beta}, we have
$k=1$, $\delta_1=-4$, and $\delta_2=4$ so that $S:=\{F=0\}\subset\mathbb P^3_{w,x,y,z}$
belongs to the $\quarticfamilylink{B2-beta}{(B.2;\beta)}$ family.

Consider the cyclic quartic cover 
$Y:=\{u^4=F\}\subset \mathbb{P}^4_{w,x,y,z,u}$ branched along $S$. 
Over the open chart $D_{\mathbb P^3}(x)$, the quartic cover has equation 
\[
w^2+(y^2-1)z^2+z+y^2+1-u^4=0.
\]
We set $R=y^2-1$ and $H=-y^4+\frac{5}{4}$ and homogenizing in $w$ and $z$ with a new variable $v$ gives
\[
w^2+(y^2-1)z^2+zv+(y^2+1-u^4)v^2=0.
\]
This is a conic bundle over  $\mathbb A_{y,u}$ with the Brauer class
$\alpha=(R,Ru^4+H)$. Moreover, the discriminant curve of the conic bundle is
birational to the smooth projective model $C$ of $u^4=\frac{4y^4-5}{4y^2-4}$. 
The polynomiasl $R$ and $H$
are square free and also coprime.
The morphism $C\to \mathbb{P}^1_y$ has six branch points 
of index four and one branch point having index two.
Thus, we conclude that $B_4=6$ and $B_2=1$ for this finite morphism. As a consequence $3B_4+2B_2=20$.
Since the degree of $R$ is even, all the hypotheses
of Proposition~\ref{prop:fractional-kummer-irrationality}
are satisfied. We conclude
that $Y$ is irrational, therefore
the pair $(\pp^3,S)$ is not of cluster type.
\end{proof}

\begin{proposition}\label{prop:3betaclustertype}
Let $S\in\mathcal I^{\quarticfamilylink{B3-beta}{(B.3;\beta)}}$. Then
$(\mathbb P^3,S)$ is of cluster type.
\end{proposition}

\begin{proof}
Choose coordinates $[w:x:y:z]$ in which
$C=V(w,Q)$, $Q=xz-y^2$, and $L=V(x,y)$.
Since $S$ is irreducible and double along $C$ and $L$,
its equation can be written as
\[
 F=Q^2+wQ(ax+by)+w^2(cx^2+dxy+ey^2).
\]
The curves meet at $P=[0:0:0:1]$. On the open chart $D_{\mathbb{P}^3}(z)$, set
$q=x-y^2$. For the weighted blow-up at $P$ with weights
$(2,1,1)$ on $(q,y,w)$, the weighted initial form of $F$ is
$q^2+bqyw+ey^2w^2$. Thus, we have $b^2-4e \neq 0$ by
Lemma~\ref{lem:lc-quadratic}.

Now, we work on the chart $D_{\mathbb P^3}(x)$ with $w\ne0$. Set $p=(z-y^2)/w$.
Substituting $z=y^2+wp$ gives
\[
 F|_{x=1}=w^2f(y,p),\qquad\text{with}\qquad
 f(y,p)=p^2+byp+ey^2+ap+dy+c.
\]
Since $b^2-4e\ne0$, we can choose affine coordinates $(s,t)$
in which $f=st-\eta$ for some $\eta\in\mathbb C$.
Applying Lemma~\ref{prop:shear-torus-chart} on
$\mathbb G_m\times\mathbb A^2_{s,t}$ gives an open embedding
$\mathbb G_m^3\hookrightarrow\mathbb P^3\setminus S$.
Hence $(\mathbb P^3,S)$ is of cluster type.
\end{proof}

\begin{proposition}\label{B4betaclustertype}
Let $S\in\mathcal{I}^{\quarticfamilylink{B4-beta}{(B.4;\beta)}}$.
Then $(\mathbb P^3,S)$ is of cluster type.
\end{proposition}

\begin{proof}
By \cite[Lemma 2.4]{skewline}, we may assume that
the three lines are $L_1=V(w,z)$, $
L_2=V(w,x)$ and $
L_3=V(x,y).$ Since $F$ vanishes to multiplicity two along each line,
we have $F\in(w,z)^2\cap(w,x)^2\cap(x,y)^2.$ Hence we may write
\[
F=ax^2z^2+bx^2zw+cx^2w^2
+dxyzw+exyw^2+hy^2w^2.
\]

On the open chart $D_{\mathbb P^3}(wx)$, set $u=y/x$ and $F|_{w=1}=x^2q(u,z)$, where
\[
q(u,z)=hu^2+duz+az^2+eu+bz+c.
\]
Since $S$ is irreducible, the polynomial $q$ is irreducible. If the quadratic part of $q$ is nondegenerate,
we may choose affine coordinates $(r,s)$
such that $q=rs-\eta$ for some
$\eta\in\mathbb C^\times$.
Otherwise, irreducibility implies that
$q$ can be completed to a
coordinate system $(r,q)$ of $\mathbb A^2$.

In the first case, we recover $s$ as $(q+\eta)/r$ on the open chart $D(r)$.
Thus, in both cases, the coordinates $(x,r,q)$ identify
$D(xrq)$ with $\mathbb G_m^3$.
Since $F|_{w=1}=x^2q$, this torus is an open subset of
$\mathbb P^3\setminus S$.
By Lemma~\ref{clustertypeandtorusembedding},
the pair $(\mathbb P^3,S)$ is of cluster type.
\end{proof}

\subsection{Cluster type condition for $(\mathbb P^3,S)$ when $S$ is irreducible and non-normal}

We now complete the proof of Theorem~\ref{them:mainclusternonnormal}.
\begin{proof}[Proof of Theorem~\ref{them:mainclusternonnormal}]
By Propositions~\ref{prop:line-conductorstrata}
and~\ref{prop:B2-low-cusp-strata}.(i),
the equisingular strata used below are smooth.
Moreover, Proposition~\ref{prop:simultaneous-urabe-dlt}
provides \'etale local simultaneous
$\mathbb Q$-factorial dlt modifications for the log pairs $(\pp^3,S)$ considered below.
Hence, by Corollary~\ref{cor:p3-families}, the cluster type
property is constant on each connected component
of these strata. 
Suppose first that the double curve is a line.
If $S\in \mathcal{I}^{\quarticfamilylink{B1-alpha}{(B.1;\alpha)}}$ then $(\mathbb P^3,S)$ is of cluster type
by Proposition~\ref{prop:11-splitting-cluster}.
If $S\in \mathcal{I}^{\quarticfamilylink{B1-beta}{(B.1;\beta)}}$ then
Proposition~\ref{prop:three-irrational-models} gives, in each of the strata, a member whose
associated cyclic quartic cover is irrational
$\mathcal{T}_{4,4}^{\quarticfamilylink{B1-beta}{(B.1;\beta)}},
\mathcal{T}_{6}^{\quarticfamilylink{B1-beta}{(B.1;\beta)},II}$ and
$\mathcal T_{6,3}^{\quarticfamilylink{B1-beta}{(B.1;\beta)}}$.
By Proposition~\ref{prop:line-conductorstrata},
all the remaining $\quarticfamilylink{B1-beta}{(B.1;\beta)}$ equisingular strata with
$q\le 6$ specialize to these strata.
Thus $(\mathbb P^3,S)$ is not of cluster type
in this range.

Suppose next that the double curve is a smooth conic.
Proposition~\ref{prop:two-cusp-irrational}
gives a member with $S
\in \mathcal{T}_{3,3}^{\quarticfamilylink{B2-alpha}{(B.2;\alpha)}}$
whose associated cyclic quartic cover is irrational.
By Proposition~\ref{prop:B2-low-cusp-strata},
the stratum $\mathcal{T}_{3}^{\quarticfamilylink{B2-alpha}{(B.2;\alpha)}}$ contains
$\mathcal{T}_{3,3}^{\quarticfamilylink{B2-alpha}{(B.2;\alpha)}}$ in its closure.
Hence $(\mathbb P^3,S)$ is not of cluster type
for the strata with $q=3$.

It remains to consider the complementary ranges.

For the \textup{\quarticfamilylink{B1-beta}{(B.1;$\beta$)}} family, $(\mathbb P^3,S)$
is of cluster type when $q\ge7$ by
Proposition~\ref{prop:line-conductorstrata}.\textup{(4)}.
For the \textup{\quarticfamilylink{B2-alpha}{(B.2;$\alpha$)}} family, $(\mathbb P^3,S)$
is of cluster type when $q\ge4$ by
Proposition~\ref{prop:B2-low-cusp-strata}.\textup{(iii)}.

Finally, suppose that the non-normal locus is reducible.
For the \textup{\quarticfamilylink{B2-beta}{(B.2;$\beta$)}} family, if $S$ has at
least two degenerate cusps, then $(\mathbb P^3,S)$ is of
cluster type by Proposition~\ref{prop:ClusterB.2.beta}.
Otherwise, Proposition~\ref{prop:two-line-dlt} shows that the
\textup{\quarticfamilylink{B2-beta}{(B.2;$\beta$)}} surfaces for which the cusp is the only
zero dimensional log canonical center of $(\mathbb P^3,S)$ and
has type $T_{2,\infty,\infty}$ form a smooth affine family.
All log pairs $(\mathbb P^3,S)$ in this family admit a simultaneous
$\mathbb Q$-factorial dlt modification. Moreover,
Proposition~\ref{prop:two-line-irrational} shows that this family
contains a surface whose associated cyclic quartic cover, branched over it,
is irrational. The families \textup{\quarticfamilylink{B3-beta}{(B.3;$\beta$)}} and \textup{\quarticfamilylink{B4-beta}{(B.4;$\beta$)}}
are of cluster type by the Propositions~\ref{prop:3betaclustertype} and~\ref{B4betaclustertype}, and the
case \textup{\quarticfamilylink{B3-alpha}{(B.3;$\alpha$)}} follows from \cite[Theorem~6.1]{EFM26}. This completes the irreducible non-normal case.
\end{proof}

\section{Cusp singularities and cluster type degenerations of \texorpdfstring{$\pp^3$}{P3}}
\label{sec:cusp-sing-cluster-type-deg}

In this section, we prove that there are precisely $92$ cusp singularities $T_{p,q,r}$ appearing in quartic surfaces $S\subset \pp^3$ for which the pair $(\pp^3,S)$ is of cluster type. 
Then, we use these quartic surfaces and produce some new examples of klt cluster type degenerations of $\pp^3$.

\begin{proposition}\label{prop:101+92}
There are exactly $101$ analytic type of cusp singularities $T_{p,q,r}$ that occur at the
unique non-klt center $o$ of a normal log canonical quartic surface
$S\subset\mathbb P^3$ of coregularity zero. Of these, exactly $92$
are realized by quartics $S$ for which $(\mathbb P^3,S)$ is of cluster type.
\end{proposition}

\begin{proof}
Let $p,q$, and $r$ be the indices of the cusp singularity. Suppose $p\ge3$. Write $S=V(wf_3+f_4)$, where $f_i$ is homogeneous of degree $i$ (see~\cite[Equation (4.1.1)]{Deg90}). On the exceptional divisor $E_1$ of the blow-up at $o$, set $C=V_{E_1}(f_3)$ and $D=V_{E_1}(f_4)$. By ~\cite[Lemma~4.2 and Section~4.3.2]{Deg90}  if there are cusp $(S;o)\cong T_{p,q,r}$ then there are $A_{p-4},A_{q-4}$ and $A_{r-4}$-singularities on the strict transformation $\tilde{S}$ of $S$ in the blow-up. Let $C_i$ be an irreducible component of $C$. Then by B\'ezout's theorem the sum of the local intersection number between $D$ and $C_i$ is at most $4\deg C_i$. 

When $p=q=3$, then the curve $C$ has a unique singular point, at which the local intersection number of $C$ and $D$ is at most $12$. Hence $4\leq r\leq  15$. In the case that $p=3<q$, the curve $C$ consists of two irreducible components $C_1$ and $C_2$. If $q=4$, then we have $r-4\leq 8$ and hence $r\leq 12$ (see~\cite[Section~4.7]{Deg90}). If $q\geq 5$, then there is a component $C_i$ with a node of index $r$
has multiplicity $r-4$ also has multiplicity at least $1$
at the node of index $q$
\cite[Sections~4.3.2 and~4.7(b)]{Deg90}.
Its total intersection number is at most $8$ by~\cite[Section~4.7(a)]{Deg90}. So, we have $(r-4)+1\le8$ and $r\le11$. If $q\ge8$, both $q-4$ and $r-4$ must occur on the conic, since
the line meets $D$ at both nodes and has total intersection number $4$
with $D$. Thus, we conclude that $(q-4)+(r-4)\le8$, so $q=r=8$.

If $4\leq p$ then $C$ is a triangle, i.e., a cycle of three rational curves intersecting transversally. If $p=4$ then the above argument gives $4\leq q\leq r\leq 8$. If $p\geq 5$, then since each local intersection number is bounded by $4$ we have $(r-4)+1\leq 4$. By \cite[Section 4.8]{Deg90} all the values $5\leq q\leq r\leq 7$ are possible. This sums up to $65$ possible local analytic types $(S;o)\cong T_{p,q,r}$. By \cite[Section~4.8]{Deg90}, all triples satisfying
$5\le p\le q\le r\le7$ are possible as a cusp singularity of the log canonical quartic surface.

Now suppose $p=2$ and $q\ge4$.
Then $S$ belongs to the family \textup{\quarticfamilylink{A2}{(A.2)}}.
In this case, we can use Notation~\ref{not:T}.
By Lemma~\ref{lem:A2-singularity-origins},
if $b_2(p_z)=0$, then $r\le8$.
Thus, to establish the upper bound on $r$,
we may assume that $b_2(p_z)\ne0$.
In this case, the curve $D$ has an $A_{r-5}$-singularity
at $\bar p_z$. If $D$ is irreducible then $p_a(D)=3$ and $\left\lfloor\frac{r-4}{2}\right\rfloor=\delta_{p_z}(D)\le3,$
so $r\le11$. If $q\geq 5$, then $C$ is reducible. Suppose $r\ge12$. By Lemma~\ref{lem:A2-singularity-origins},
we have $b_2(\bar p_z)\ne0$ and $D$ has an $A_{r-5}$-singularity at
$\bar p_z$. By B\'ezout theorem, we conclude that 
\[
4\le \frac{r-4}{2}=I_{\bar p_z}(D_1,D_2)
\le \deg D_1\deg D_2\le4.
\]
Hence $r=12$, and $D_1,D_2$ are conics meeting only at $\bar p_z$.
Since $I_{\bar p_z}(D,L)=2$, the line $L$ meets each conic
transversally there. The two remaining intersections of $L$ with $D$
are distinct, so $q=4$. Now suppose $q\ge8$.
By Corollary~\ref{eq:A2-index-bound} and
Lemma~\ref{lem:A2-singularity-origins},
the images of $p_y$ and $p_z$ on $D$ have
singularities of types $A_{q-5}$ and $A_{r-5}$,
respectively. Since $q-5,r-5\ge3$, the classification of reduced plane quartics forces both points to be of type $A_3$ (see~\cite[Tables~1 and~2(a)]{BDLL25}).
Hence $q=r=8$.

 We now construct quartic surfaces with the prescribed cusp types. In the variables of $x,y,z,w$ of the projective space $\mathbb P^3$ we can consider the homogeneous polynomial $4F=(2wx+b_2)^2-\Delta$ from \eqref{eq:A2-square-form}. For $(q,r)=(4,12)$ and $(8,8)$, we can choose $\Delta=(yz-x^2)^2+y^4$ and $
\Delta=(yz-x^2)(yz-2x^2)$ respectively. The two conics meet with multiplicity four at
$\bar p_z$ in the first case and with multiplicity two at
each of $\bar p_y$ and $\bar p_z$ in the second. 
In the remaining cases, choose a plane quartic $D$ of type $A_{r-5}$
at $\bar p_z$ and smooth elsewhere, where $A_0$ denotes a smooth point. For $q=4$, take $L$ to meet $D$ with multiplicity $2$ at $\bar p_z$
and transversely at two further points. For $5\le q\le7$, take $L$ to meet $D$ with multiplicity $2$ at
both $\bar p_z$ and a smooth point $\bar p_y\ne\bar p_z$.

If $q=3$, then Proposition~\ref{Theorem:q+rleq18} gives $r\le15$. In family \textup{\quarticfamilylink{A3}{(A.3)}}, set $j=2k$ in
\eqref{equation:strict-transformation-A_3} and take
$(k,f)=(2,0),(0,-\tfrac12),(0,0)$ to obtain $r=13,14,15$,
respectively.

Counting the cusp types appearing on quartic surfaces
in the families \textup{\quarticfamilylink{A2}{(A.2)}}--\textup{\quarticfamilylink{A4}{(A.4)}}
for which $(\mathbb P^3,S)$ is of cluster type,
together with the $65$ types above,
we obtain a total of $92$ distinct analytic types
of cusp singularities.
\end{proof}

We now construct new examples of klt toric degenerations
of $\mathbb P^3$ with a relative cluster type boundary of $\mathbb P^3$. 

\begin{example}
{\em 
Consider the weighted projective space $\pp:=\pp(1,1,1,1,2,3,3,5)$ with variables $w,x,y,z$, $u,v,p,q$. 
Consider the family $\mathcal X\subset\mathbb P\times\mathbb A^1$
which is the scheme theoretic closure of the complete
intersection over 
$\mathbb A^1\setminus\{0\}$ defined by
\begin{equation}\label{equationofF}t^3u=wx+z^2,\quad t^2v=zu+xy^2,\quad t^2p=wu-zy^2,\quad\text{and}\quad tq=up-y^2v.\end{equation}
For $t\ne0$, the fiber $\mathcal X_t$ is isomorphic
to $\mathbb P^3$, since the variables $u,v,p,q$
are uniquely determined by $w,x,y,z$. The special fiber $\mathcal X_0$ is the Zariski closure
of the image of the torus $\mathbb G_m^3$
under the embedding
\[
(w,p,q)\longmapsto
[\,w:-w^3q^2:1:-w^2q:-wq:-wpq:p:q\,].
\]
Thus, the variety $\mathcal{X}_0$ is a toric variety. This toric variety $\mathcal{X}_0$ does not appear in the literature. Indeed, $\mathcal X_0$ is not $\mathbb Q$-factorial
and has a point of local Gorenstein index $3$.
The previously known examples are either
$\mathbb Q$-factorial or non-$\mathbb Q$-factorial
with Gorenstein index $1$ or $5$
(see \cite[Theorem~5.1]{Dev22} and \cite{DVLT26}).

We argue that $\mathcal{X}_0$ is not $\qq$-factorial at the torus invariant point $o$ where all variables except $p$ vanishes. 
Moreover, its local Gorenstein index at such point equals $3$. 
Consider the cyclic cover of $D_{\mathcal X_0}(p)$ induced by the cyclic morphism from $\mathbb A^7$ to its ambient space. The defining relations give $u=y^2v, z=wv$ and $x=-wv^2.$
Substituting into $wq=-y^4u$ gives $F=wq+y^6v=0$, with no further relations. Thus the cyclic cover of $D_{\mathcal{X}_0}(p)$ is isomorphic to $Y=V(F)\subset\mathbb A^4_{w,y,q,v}$. The quotient is given by $\mu_3$ acting
with weights $(1,1,2,0)$.

The open set $wyqv\ne 0$ is a torus with coordinates $w,y,v$. Its complement has four prime divisors $V(w,y), V(w,v), V(q,y)$ and $V(q,v)$. These divisors are $\mu_3$-equivariant. Hence $\mathcal{F}$ is not $\mathbb Q$-factorial around $o$. The
canonical form $\omega=\operatorname{Res}_Y
\frac{dw\wedge dy\wedge dq\wedge dv}{wq+y^6v}$
has $\mu_3$-weight $1$.
Hence the local Gorenstein index of $\mathcal X_0$
at $o$ is $3$.

Define a boundary divisor $\mathcal{B}$ on $\mathcal{X}$ by $F=yv+t x^2(y^2+z^2)+t^4x^3z$. We finish the example by arguing that the pair $(\mathcal{X},\mathcal{B})$ is of cluster type over $\mathbb{A}^1$. On the open chart $D_{\mathcal X}(y)$, set $y=1$.
On the open subset where $xzF\ne0$, we have $v=F-tx^2(1+z^2)-t^4x^3z$ and $u=\frac{t^2v-x}{z}$, 
thus 
\[
w=\frac{t^3u-z^2}{x},\qquad
p=\frac{wv-tu}{z},\qquad
q=\frac{tpv-u^2}{z}.
\]
Therefore, the coordinate functions $(x,z,F)$ give an isomorphism $D_{\mathcal X}(yxzF)
\cong\mathbb A^1_t\times\mathbb G_m^3$ over $\mathbb A^1$.
}
\end{example} 

\begin{example}
{\em 
Consider the weighted projective space 
$\pp:=\pp(1,1,1,1,2,3,3)$ with variables $w,x,y,z,u,v,p$. Consider the family 
$\mathcal{X}\subset \pp\times \mathbb{A}^1$ to be the defined by the equations 
\[
t^2u=wx+yz,\quad  tv=wu-z^3,\quad tp=yu+xz^2,\quad xv+zp=tu^2,\quad  wp-yv=tuz^2.
\]
For $t\neq 0$, the fiber $\mathcal{X}_t$ is isomorphic
to $\mathbb{P}^3$, since the variables $u,v,p$
are uniquely determined by $w,x,y,z$.
We argue that the toric special fiber $\mathcal{X}_0$ is novel in the literature.
Indeed, we will see that it is not $\mathbb{Q}$-factorial and has a point
of local Gorenstein index $3$. 

Let $o_u$ and $o_p$ be the points where all variables vanishes except the one in the subindex. 
We argue that $\mathcal{X}_0$ is not $\qq$-factorial at $o_u$ and has local index $3$ at $o_p$. Consider the open subsets $D_{\mathcal X_0}(u)$ and
$D_{\mathcal{X}_{0}}(p)$ and consider the cyclic cover of $D_{\mathcal{X}}(p)$ induced by the cyclic map from $\mathbb A^6$ to its ambient space. The cover of $\mathcal{X}_0$ is defined by
$Y_u=V(xv+zp)\subset\mathbb A^4_{x,z,v,p}$ and $Y_p=V(yu+x^3v^2)\subset\mathbb A^4_{x,y,u,v}.$
Since $Y_u$ is an ordinary double point with a toric quotient action $\mu_2$. Hence $\mathcal X_0$ is not
$\mathbb Q$-factorial at $o_u$.

On $Y_p$, the canonical form $\omega=\operatorname{Res}_{Y_p}
\frac{dx\wedge dy\wedge du\wedge dv}{yu+x^3v^2}$
has $\mu_3$-weight $1$. Hence the local Gorenstein index of $\mathcal{X}_0$ at $o_p$ is $3$.

Finally, we argue that this family $\mathcal{X}\rightarrow \mathbb{A}^1$ is of cluster type over the base. 
Set $a=x+tw$ and $b=v+tux+t^2uw=v+tua$.
We define the boundary $\mathcal{B}=V_{\mathcal X}(F)$ by 
\[
F=tu^2-ab
 =zp-twv-tux^2-2t^2uwx-t^3uw^2.
\]
We argue that the pair $(\mathcal X,\mathcal{B})$ is of relative
cluster type over $\mathbb A^1$. On the open chart $D_{\mathcal{X}}(z)$, put $z=1$. The defining equations reduce to $uw=1+tv$,
$y=t^2u-wx$, and $p=tu^2-xv$. On $D_{\mathcal{X}}(uaF)$, we have
\[
v=\frac{tu^2-F}{a}-tua,\qquad
w=\frac{1+tv}{u},\qquad\text{and}\qquad x=a-tw.
\]
The remaining coordinates $y,p$ are then uniquely
determined by the above relations.
Thus the coordinate functions $(u,a,F)$ give an isomorphism 
$D_{\mathcal F}(zuaF)\cong\mathbb A^1_t\times\mathbb G_m^3$.
}
\end{example} 

\appendix

\section{Use of AI disclosure}
\label{appendix-A}

In this appendix, we explain how the tools of Artificial Intelligence (AI), more precisely ChatGPT,\footnote{The precise models and names changed during the completion of this work, so it does not seem necessary to specify.} were used to understand the cluster type boundaries on $\mathbb{P}^3$. 
Below, we explain in a chronological way which ideas came from the authors and which ideas came from AI.

The idea of using a modified version of invariance of plurigenera to prove the constancy of the cluster type condition on equisingular families of pairs is due to the authors. 
This idea extends earlier work of the first
author with Ji and with Zúñiga. This idea and the work of Ducat, Urabe, and Wall simplified the task of understanding which pairs $(\pp^3,S)$, with $S$ a quartic surface, are cluster type. Indeed, now the task was reduced to consider a finite set of pairs $\{(\pp^3,S_i)\mid i\in J\}$. The authors used AI in order to test different obstructions for these pairs to be of cluster type. AI proposed using Lemma~\ref{lem:clustertypeobstruction} and Shokurov's irrationality criterion on the obtained quartic threefold. This strategy required some careful computations in each case. At this stage, the set $J$ of equisingular strata 
could have more than $100$ different cases and the task was not feasible. Then, the authors had the idea of using the fact that in suitable families of threefolds over curves, one irrational fiber implies that a very general fiber is irrational. AI proposed a precise reference for this fact in the works of de Fernex and Fusi. Using this idea and previous work due to Alves da Silva, Figueroa, and the first author, the authors reduced the task to a subset $J_0\subsetneq J$ of $14$ families.
These are the $14$ families which are presented in the sketch
of the proof in Subsection~\ref{strat}. Finally, the authors used AI to help construct, in each of these
families, a quartic with an irrational cyclic quartic cover. The authors then modified some of the examples given by AI to simplify the exposition and wrote the present article.

\bibliographystyle{habbvr}
\bibliography{bib}

\vspace{0.5cm}
\end{document}